\documentclass[a4paper,12pt,twoside]{article}
\usepackage{amsmath,amssymb,ifthen}
\usepackage[amsmath,thmmarks]{ntheorem}
\usepackage{mathrsfs}
\usepackage{paper}
\usepackage[all]{xy}
\usepackage{url}
\DeclareMathOperator{\Spf}{Spf}
\DeclareMathOperator{\Spa}{Spa}
\DeclareMathOperator{\rec}{rec}
\DeclareMathOperator{\spp}{sp}
\DeclareMathOperator{\cosp}{cosp}
\DeclareMathOperator{\supp}{supp}
\DeclareMathOperator{\height}{height}
\DeclareMathOperator{\Gys}{Gys}
\DeclareMathOperator{\Ind}{Ind}
\DeclareMathOperator{\Nr}{Nr}
\DeclareMathOperator{\Nrd}{Nrd}
\DeclareMathOperator{\JL}{JL}
\DeclareMathOperator{\SL}{SL}
\DeclareMathOperator{\Iw}{Iw}
\DeclareMathOperator{\cInd}{c-Ind}

\DeclareMathOperator{\Art}{Art}

\DeclareMathOperator{\algInd}{algInd}

\newcommand{\ITw}{\mathrm{IT}_{\mathrm{w}}}

\SelectTips{cm}{11}
\title{Geometric approach to the explicit local Langlands correspondence}
\author{Yoichi Mieda}

\begin{document}
\maketitle

\begin{firstfootnote}
 Graduate School of Mathematical Sciences, The University of Tokyo, 3--8--1 Komaba, Meguro-ku, Tokyo, 153--8914, Japan

 E-mail address: \texttt{mieda@ms.u-tokyo.ac.jp}

\end{firstfootnote}

\begin{abstract}
 We propose a geometric strategy of giving explicit description of the Langlands parameter of
 an irreducible supercuspidal representation of $\GL(n)$ over a non-archimedean local field.
 The key is to compare the cohomology of an affinoid in the Lubin-Tate space at infinite level
 and that of the reduction of its formal model.
 As examples, we treat the cases of depth $0$ supercuspidal representations and simple supercuspidal
 representations.
\end{abstract}

\section{Introduction}\label{sec:intro}
Let $F$ be a non-archimedean local field. 
By the local Langlands correspondence for $\GL_n(F)$ (\cite{MR1876802}, \cite{MR1228127}),
irreducible supercuspidal representations of $\GL_n(F)$ are parameterized by $n$-dimensional
irreducible smooth representations of the Weil group $W_F$ of $F$.
On the other hand, irreducible supercuspidal representations of $\GL_n(F)$ are
completely classified by the theory of types \cite{MR1204652}.
However, in general it is very difficult to determine the parameter corresponding to
each supercuspidal representation explicitly. 
This problem, called ``the explicit local Langlands correspondence'',
is extensively studied by Bushnell and Henniart;
see \cite{MR2138141}, \cite{MR2148193}, \cite{MR2679700}, \cite{MR2275640} and references therein.
Recently, Imai and Tsushima \cite{2015arXiv150902960I} also gave a result in this direction,
which is not covered by the works cited above.
The methods in all of these works are purely algebraic, and sometimes involve very detailed
and complicated computations.

The aim of this paper is to propose a new geometric method for the explicit local Langlands correspondence.
It is well-known that the local Langlands correspondence for $\GL_n(F)$ has a nice geometric realization;
it appears in the middle degree $\ell$-adic cohomology of the Lubin-Tate tower.
This fact is called the non-abelian Lubin-Tate theory
after Carayol \cite{MR1044827}, and has been proved by Harris-Taylor \cite{MR1876802} and Boyer \cite{MR1719811}.
Therefore, it is quite natural to expect that we can understand the local Langlands correspondence
by studying the geometry of the Lubin-Tate tower.
Apart from many results in the case $n=2$, the first result in this direction was given by Yoshida \cite{MR2676163}.
He constructed a semistable model of the first layer of the Lubin-Tate tower, and discovered
the Deligne-Lusztig variety for $\GL_n$ in the reduction of it.
By using this result, he concluded that a part of the non-abelian Lubin-Tate theory
boils down to the Deligne-Lusztig theory \cite{MR0393266}.
More recently, the Lubin-Tate space at infinite level,
namely the projective limit of the Lubin-Tate tower, turns out to be much simpler than individual layers.
This idea goes back to Faltings \cite{MR1936369} and Fargues \cite{MR2441312}, and now is formulated
by the theory of perfectoid spaces \cite{Scholze-perfectoid};
see \cite{Scholze-Weinstein} and \cite{Weinstein-stabred} for detail.
Based on this theory, Weinstein and Boyarchenko \cite{MR3402698} constructed a family of affinoids in
the Lubin-Tate space at infinite level.
They also constructed formal models of the affinoids,
and observed that a part of the local Langlands correspondence is encoded in the $\ell$-adic cohomology
of their reduction.
Following this work, some other authors also obtain similar results;
see \cite{Imai-Tsushima-tame}, \cite{Imai-Tsushima-wild} and \cite{Tokimoto-thesis}.
The common strategy in these works is to combine an explicit computation of the reduction
of the formal models and a previously known explicit description of the local Langlands correspondence.
The idea of this article is to reverse this process; we combine a computation of the reduction
and the non-abelian Lubin-Tate theory to obtain a result in the explicit local Langlands correspondence.
To carry it out, we need to relate the cohomology of the reduction with that of the Lubin-Tate tower itself,
which is not included in the previous works.
The main ingredient of this article is to give a useful condition which ensures the existence of such a relation.

Here is a key technical result for us:

\begin{thm}[Definition \ref{defn:specialization-map-infinite-level}, Proposition \ref{prop:specialization-map-functoriality}, Theorem \ref{thm:injectivity-criterion-infinite-level}]\label{thm:specialization-injectivity-intro}
 Let $k$ be a complete algebraically closed non-archimedean field, $k^\circ$ its valuation ring,
 and $\kappa$ the residue field of $k^\circ$.
 Assume that the characteristic of $\kappa$ is positive.
 Fix a non-zero element $\varpi$ in the maximal ideal of $k^\circ$.
 Let $A$ be a $\varpi$-adically complete flat $k^\circ$-algebra endowed with an action of
 a profinite group $K_0$. We write $\widetilde{A}$ for the integral closure of $A$ in $A[1/\varpi]$.
 We assume the following two conditions (see Section 3 for detail):
 \begin{itemize}
  \item The space $\Spa(A[1/\varpi],\widetilde{A})$ is the projective limit of a tower $\{X_K\}_K$
	of smooth affinoids over $k$ indexed by open normal subgroups $K$ of $K_0$.
  \item The reduction $\Spec (A\otimes_{k^\circ}\kappa)$ of $\Spf A$ is the perfection of an affine scheme
	$Y$ of finite type over $\kappa$.
 \end{itemize}
 \begin{enumerate}
  \item For a prime number $\ell$ invertible in $\kappa$,
	we have a map 
	\[
	\spp^*\colon H^i_c(Y,\overline{\Q}_\ell)\to \varinjlim_K H^i_c(X_K,\overline{\Q}_\ell),
	\]
	which we call the specialization map.
	It is functorial with respect to automorphisms of $A$ of finite level 
	(i.e., automorphisms coming from those of the tower $\{X_K\}_K$).
  \item Assume moreover that $Y$ is pure-dimensional and smooth over $\kappa$.
	Let $\{Z_K\}_K$ be a tower of rigid spaces which contains $\{X_K\}_K$ as a tower of open rigid subspaces.
	For a subspace $V\subset H^i_c(Y,\overline{\Q}_\ell)$ such that the composite
	\[
	 V\hookrightarrow H^i_c(Y,\overline{\Q}_\ell)\to H^i(Y,\overline{\Q}_\ell)
	\]
	is injective, the composite
	\[
	 V\hookrightarrow H^i_c(Y,\overline{\Q}_\ell)\xrightarrow{\spp^*} \varinjlim_K H^i_c(X_K,\overline{\Q}_\ell)\to \varinjlim_K H^i(Z_K,\overline{\Q}_\ell)
	\]
	is also injective.
 \end{enumerate}
\end{thm}
In the application, $\{Z_K\}$ becomes the Lubin-Tate tower, $\Spa(A[1/\varpi],\widetilde{A})$ an affinoid
in the Lubin-Tate space at infinite level, and $\Spf A$ a formal model of it.
We remark that, unlike the case of finite type, the specialization map $\spp^*$ is not necessarily surjective,
even if $Y$ is smooth (see Example \ref{exa:non-reduced}).

To demonstrate how Theorem \ref{thm:specialization-injectivity-intro} is used to give an explicit description
of the local Langlands correspondence, we choose two classes of irreducible supercuspidal representations
as examples; one is depth $0$ supercuspidal representations and the other is simple supercuspidal
representations. Recall that the latter are supercuspidal representations with minimal positive depth;
see \cite{MR2730575}, \cite{MR3164986}, \cite{Knightly-Li}.
By geometric method, we will obtain the following result:

\begin{thm}\label{thm:explicit-LLC-intro}
 Let $\pi$ be either a depth $0$ supercuspidal representation or a simple supercuspidal representation.
 Then, we can describe the Langlands parameter of $\pi$ explicitly.
\end{thm}
For a more precise statement, see Theorems \ref{thm:LLC-JLJC-depth-0} and \ref{thm:LLC-JLJC-ssc}.
We can also determine the image of $\pi$ under the local Jacquet-Langlands correspondence
between $\GL_n(F)$ and $D^\times$, where $D$ is the central division algebra over $F$ with invariant $1/n$.
Note that the resulting description in Theorem \ref{thm:explicit-LLC-intro} is not new,
as in Remarks \ref{rem:history-depth-0}, \ref{rem:history-ssc}.
However our proof is totally different from the previous ones.
We do not need to investigate involved representations closely.
In fact, we do not have to compute neither characters nor epsilon factors.
The author expects that the same method can be applied to supercuspidal representations with
larger depth. This will be considered in our future works.

The idea to see the map $H^i_c(Y,\overline{\Q}_\ell)\to H^i(Y,\overline{\Q}_\ell)$ in
Theorem \ref{thm:specialization-injectivity-intro} (ii) emerged from discussion with Takahiro Tsushima.
After we wrote this paper, we are informed that Tsushima is also working on a similar problem
in the case where $n=2$ and $F$ is an odd equal characteristic local field.
He also pays attention to the map $H^i_c(Y,\overline{\Q}_\ell)\to H^i(Y,\overline{\Q}_\ell)$, but
his basic strategy is to find a family of affinoids with good reduction
in the Lubin-Tate spaces at suitable finite levels, which is different from ours.

The outline of this paper is as follows. In Section \ref{sec:finite-level}, we give an analogous result
as Theorem \ref{thm:specialization-injectivity-intro} in the finite level setting by using
the theory of formal nearby cycles. In Section \ref{sec:infinite-level},
we deduce Theorem \ref{thm:specialization-injectivity-intro} from the finite level case.
To enable this step, we should allow bad reduction in Section \ref{sec:finite-level},
even if the reduction appearing in Theorem \ref{thm:specialization-injectivity-intro}
is the perfection of a smooth affine scheme over $\kappa$.
In Section \ref{sec:LT-tower}, we recall the definition of the Lubin-Tate tower and the statement
of the non-abelian Lubin-Tate theory. We also verify that Theorem \ref{thm:specialization-injectivity-intro} is
applicable to an affinoid in the Lubin-Tate space at infinite level.
In Section \ref{sec:depth-0}, we explicitly describe the Langlands parameter of a depth $0$ supercuspidal
representation. We use a formal model of an affinoid whose reduction becomes the perfection of
the Deligne-Lusztig variety for $\GL_n$.
This formal model can be seen as the infinite level version of a piece of the model constructed 
by Yoshida \cite{MR2676163}.
In Section \ref{sec:ssc}, we give a description of the Langlands parameter of a simple supercuspidal
representation. Here we use results in \cite{Imai-Tsushima-tame} and \cite{Imai-Tsushima-wild}.
Appendix \ref{sec:ring-theory} contains some results for algebras over the valuation ring $k^\circ$
which do not satisfy any finiteness condition.
The results on the cohomology of the Artin-Schreier sheaves in Appendix \ref{sec:Artin-Schreier} are
used in Section \ref{sec:ssc}.

\medbreak
\noindent{\bfseries Acknowledgment}\quad
The author is grateful to Takahiro Tsushima for inspiring discussions.
He also learned many techniques to treat the $\ell$-adic cohomology of the Artin-Schreier varieties from him.
The author also thanks Kazuki Tokimoto for motivating discussions.
This work was supported by JSPS KAKENHI Grant Number 15H03605.

\section{Specialization map: the case of finite level}\label{sec:finite-level}
Let $k$ be a complete algebraically closed non-archimedean field. We write $k^\circ$ for
the valuation ring of $k$ and $\mathfrak{m}$ the maximal ideal of $k^\circ$.
The residue field $k^\circ/\mathfrak{m}$ of $k^\circ$ is denoted by $\kappa$.
We fix a prime number $\ell$ invertible in $\kappa$ and an integer $m\ge 1$, and put $\Lambda=\Z/\ell^m\Z$.

Let $\mathcal{X}$ be a quasi-compact admissible formal scheme over $\Spf k^\circ$ in the sense of \cite{MR1202394};
namely, $\mathcal{X}$ is of topologically finite type and flat over $\Spf k^\circ$.
We put $\mathcal{X}_s=\mathcal{X}\otimes_{k^\circ}\kappa$, which is a scheme of finite type
over $\kappa$. The closed immersion $\mathcal{X}_{\mathrm{red}}\hookrightarrow \mathcal{X}_s$
induces an isomorphism between \'etale sites. In the following we identify the \'etale sites
$(\mathcal{X}_{\mathrm{red}})_\et$ and $\mathcal{X}_{s,\et}$ by this isomorphism.

As in \cite[\S 1.9]{MR1734903} and \cite[\S 2]{MR1620114}, we can attach to $\mathcal{X}$
an adic space $d(\mathcal{X})$ of finite type over $\Spa(k,k^\circ)$.
Moreover, we have a morphism of sites 
$\lambda_{\mathcal{X}}\colon d(\mathcal{X})_\et\to (\mathcal{X}_{\mathrm{red}})_\et=\mathcal{X}_{s,\et}$.
We denote the right derived functor of $\lambda_{\mathcal{X}}$ by $R\Psi_{\!\mathcal{X}}$ (or simply by $R\Psi$),
and call it the formal nearby cycle functor.
By adjointness, we have a natural map $\spp^*\colon \Lambda\to R\Psi_{\!\mathcal{X}}\Lambda$.

Let us assume that $\mathcal{X}$ is pseudo-compactifiable in the sense of \cite[Definition 4.24]{formalnearby}.
Then, by \cite[Corollary 4.29]{formalnearby}, there exists a functorial isomorphism
$R\Gamma_c(\mathcal{X}_s,R\Psi_{\!\mathcal{X}}\Lambda)\xrightarrow{\cong}R\Gamma_c(d(\mathcal{X}),\Lambda)$.
Hence the map $\spp^*\colon \Lambda\to R\Psi_{\!\mathcal{X}}\Lambda$ induces a morphism
\[
 R\Gamma_c(\mathcal{X}_s,\Lambda)\to R\Gamma_c(\mathcal{X}_s,R\Psi_{\!\mathcal{X}}\Lambda)\xrightarrow{\cong}R\Gamma_c(d(\mathcal{X}),\Lambda),
\]
which is also denoted by $\spp^*$.

In the following, assume that $d(\mathcal{X})$ is smooth over $\Spa(k,k^\circ)$.
Since $\mathcal{X}$ is admissible, the dimensions of $\mathcal{X}_s$ and $d(\mathcal{X})$ are equal, 
for which we write $d$.
Then, we have the trace map $\Tr_{d(\mathcal{X})}\colon R\Gamma_c(d(\mathcal{X}),\Lambda)\to\Lambda(-d)[-2d]$.
By composing the isomorphism above, we obtain a map $R\Gamma_c(\mathcal{X}_s,R\Psi_{\!\mathcal{X}}\Lambda)\to \Lambda(-d)[-2d]$.
If we denote the structure morphism $\mathcal{X}_s\to \Spec \kappa$ by $f_s$, 
$R\Gamma_c(\mathcal{X}_s,R\Psi_{\!\mathcal{X}}\Lambda)$ is identical to 
$Rf_{s!}R\Psi_{\!\mathcal{X}}\Lambda$.
Hence, by adjointness we have the map $\cosp^*\colon R\Psi_{\!\mathcal{X}}\Lambda\to Rf_s^!\Lambda(-d)[-2d]$.

The main result in this section is as follows:

\begin{thm}\label{thm:Gysin-specialization}
 The composite of $\Lambda\xrightarrow{\spp^*}R\Psi_{\!\mathcal{X}}\Lambda\xrightarrow{\cosp^*} Rf_s^!\Lambda(-d)[-2d]$ is
 equal to the Gysin map $\Gys_{f_s}\colon \Lambda\to Rf_s^!\Lambda(-d)[-2d]$ with respect to $f_s$,
 that is, the adjoint of the trace map $\Tr_{f_s}\colon Rf_{s!}\Lambda\to \Lambda(-d)[-2d]$.
\end{thm}

\begin{lem}\label{lem:reduction}
 To prove Theorem \ref{thm:Gysin-specialization}, we may assume that $\mathcal{X}$ is affine.
\end{lem}

\begin{prf}
 Since $\mathcal{X}_s$ is $d$-dimensional, we have $H^i_c(\mathcal{X}_s,\Lambda(d))=0$ for $i>2d$.
 Therefore, by adjointness, it suffices to show that the composite of
 \[
  H^{2d}_c(\mathcal{X}_s,\Lambda(d))\xrightarrow{\spp^*} H^{2d}_c(\mathcal{X}_s,R\Psi_{\!\mathcal{X}}\Lambda(d))\xrightarrow{\cong}H^{2d}_c(d(\mathcal{X}),\Lambda(d))\xrightarrow{\Tr_{d(\mathcal{X})}}\Lambda
 \]
 is equal to $\Tr_{\mathcal{X}_s}\colon H^{2d}_c(\mathcal{X}_s,\Lambda(d))\to \Lambda$.

 Let $\mathcal{U}$ be an open formal subscheme of $\mathcal{X}$ such that $\dim (\mathcal{X}_s\setminus \mathcal{U}_s)<d$. We prove that Theorem \ref{thm:Gysin-specialization} for $\mathcal{U}$ implies that for $\mathcal{X}$.
 Consider the following diagram:
 \[
  \xymatrix{%
 H^{2d}_c(\mathcal{U}_s,\Lambda(d))\ar[r]^-{\spp^*}\ar[d]& H^{2d}_c(\mathcal{U}_s,R\Psi_{\mathcal{U}}\Lambda(d))\ar[r]^-{\cong}\ar[d]& H^{2d}_c(d(\mathcal{U}),\Lambda(d))\ar[r]^-{\Tr_{d(\mathcal{U})}}\ar[d]& \Lambda\ar@{=}[d]\\
  H^{2d}_c(\mathcal{X}_s,\Lambda(d))\ar[r]^-{\spp^*}& H^{2d}_c(\mathcal{X}_s,R\Psi_{\!\mathcal{X}}\Lambda(d))\ar[r]^-{\cong}& H^{2d}_c(d(\mathcal{X}),\Lambda(d))\ar[r]^-{\Tr_{d(\mathcal{X})}}& \Lambda\lefteqn{.}
 }
 \]
 Since the restriction of $\spp^*\colon \Lambda\to R\Psi_{\!\mathcal{X}}\Lambda$ equals
 $\spp^*\colon \Lambda\to R\Psi_{\mathcal{U}}\Lambda$, the left rectangle is commutative.
 By \cite[Lemma 4.28 (iii)]{formalnearby}, the middle rectangle commutes.
 As $d(\mathcal{U})\hookrightarrow d(\mathcal{X})$ is an open immersion, the right rectangle
 is also commutative.
 Therefore, if Theorem \ref{thm:Gysin-specialization} holds for $\mathcal{U}$,
 then the composite of
 \[
  H^{2d}_c(\mathcal{U}_s,\Lambda(d))\to H^{2d}_c(\mathcal{X}_s,\Lambda(d))\xrightarrow{\spp^*} H^{2d}_c(\mathcal{X}_s,R\Psi_{\!\mathcal{X}}\Lambda(d))\xrightarrow{\cong}H^{2d}_c(d(\mathcal{X}),\Lambda(d))\xrightarrow{\Tr_{d(\mathcal{X})}}\Lambda
 \]
 coincides with $\Tr_{\mathcal{U}_s}\colon H^{2d}_c(\mathcal{U}_s,\Lambda(d))\to \Lambda$, which is equal to
 the composite of
 \[
  H^{2d}_c(\mathcal{U}_s,\Lambda(d))\to H^{2d}_c(\mathcal{X}_s,\Lambda(d))\xrightarrow{\Tr_{\mathcal{X}_s}} \Lambda.
 \]
 On the other hand, by the exact sequence
 \[
  H^{2d}_c(\mathcal{U}_s,\Lambda(d))\to H^{2d}_c(\mathcal{X}_s,\Lambda(d))\to H^{2d}_c(\mathcal{X}_s\setminus \mathcal{U}_s,\Lambda(d))
 \]
 and the assumption $\dim (\mathcal{X}_s\setminus \mathcal{U}_s)<d$, the map
 $H^{2d}_c(\mathcal{U}_s,\Lambda(d))\to H^{2d}_c(\mathcal{X}_s,\Lambda(d))$ is surjective.
 Hence we conclude that the composite of
 \[
  H^{2d}_c(\mathcal{X}_s,\Lambda(d))\xrightarrow{\spp^*} H^{2d}_c(\mathcal{X}_s,R\Psi_{\!\mathcal{X}}\Lambda(d))\xrightarrow{\cong}H^{2d}_c(d(\mathcal{X}),\Lambda(d))\xrightarrow{\Tr_{d(\mathcal{X})}}\Lambda
 \]
 is equal to $\Tr_{\mathcal{X}_s}\colon H^{2d}_c(\mathcal{X}_s,\Lambda(d))\to \Lambda$.

 By this observation, we may replace $\mathcal{X}$ by its open formal subscheme $\mathcal{U}$
 satisfying $\dim (\mathcal{X}_s\setminus \mathcal{U}_s)<d$.
 First we remove from $\mathcal{X}$ the intersections of irreducible components of $\mathcal{X}_s$.
 Then, by considering each $d$-dimensional connected component, we may assume that $\mathcal{X}_s$ is irreducible.
 Finally, by taking an arbitrary non-empty affine open formal subscheme of $\mathcal{X}$, 
 we may reduce to the case where $\mathcal{X}$ is affine. 
\end{prf}

Later we assume that $\mathcal{X}$ is affine. 
Fix a non-zero element $\varpi\in\mathfrak{m}$.

\begin{lem}\label{lem:Elkik}
 Let $\mathcal{X}$ be an affine formal scheme of topologically finite type and flat over $\Spf k^\circ$.
 Assume that $d(\mathcal{X})$ is $d$-dimensional and smooth over $\Spa(k,k^\circ)$.
 Then there exists an affine scheme $X$ over $\Spec k^\circ$ satisfying the following conditions:
 \begin{itemize}
  \item $X$ is of finite presentation and flat over $k^\circ$,
  \item the generic fiber of $X$ is $d$-dimensional and smooth over $k$, and
  \item the $\varpi$-adic formal completion of $X$ is isomorphic to $\mathcal{X}$ over $\Spf k^\circ$.
 \end{itemize}
\end{lem}

\begin{prf}
 Write $\mathcal{X}=\Spf A$. Since $d(\mathcal{X})$ is smooth over $\Spa(k,k^\circ)$,
 $A$ is formally smooth outside $V(\varpi)$ in the sense of \cite[p.~581]{MR0345966};
 see \cite[Proposition 3.3.2]{MR2435647} (in \cite[\S 3.3]{MR2435647},
 $k$ is assumed to be a discrete valuation field of characteristic $0$,
 but it does not play a role in the proof of \cite[Proposition 3.3.2]{MR2435647}).
 Further, by \cite[Propositions 1.1 (c), 1.3]{MR1202394}, $A$ can be written
 in the form $k^\circ\langle T_1,\ldots,T_n\rangle/I$, where $I$ is a finitely presented ideal of 
 $k^\circ\langle T_1,\ldots,T_n\rangle$.
 Therefore, by \cite[Th\'eor\`eme 7 and Remarque 2 (c)]{MR0345966}, there exists 
 a finitely generated $k^\circ$-algebra $B$ such that $B\otimes_{k^\circ}k$ is smooth over $k$
 and the $\varpi$-adic completion of $B$ is isomorphic to $A$.
 Put $B_{\text{$\varpi$-tors}}=\{b\in B \mid \varpi^mb=0 \text{ for some $m\ge 0$}\}$.
 Since $A$ is $\varpi$-torsion free, \cite[Chapter 0, Proposition 7.4.5]{Fujiwara-Kato} tells us that
 the $\varpi$-adic completion of $B/B_{\text{$\varpi$-tors}}$ is again isomorphic to $A$.
 Hence, replacing $B$ by $B/B_{\text{$\varpi$-tors}}$, we may suppose that $B$ is flat over $k^\circ$.
 Then, by \cite[Corollaire 3.4.7]{MR0308104}, $B$ is a finitely presented $k^\circ$-algebra.
 
 We put $X'=\Spec B$, and denote the structure morphism $X'\to \Spec k^\circ$ by $g'$.
 We write $X'_s$ (resp.\ $X'_\eta$) for the special (resp.\ generic) fiber of $g'$.
 By Chevalley's semi-continuity theorem \cite[IV, Th\'eor\`eme 13.1.3]{EGA},
 the subset $Y=\{x\in X'\mid \dim_x g'^{-1}(g'(x))\ge d+1\}$ is closed.
 On the other hand, since $X'_s=\mathcal{X}_s$ is $d$-dimensional,
 $Y$ is contained in $X'_\eta$.
 Therefore $Y$ is equal to the disjoint union of all connected components of $X'_\eta$ whose dimensions are
 greater than $d$. In particular, $Y$ is open in $X'$.
 Put $X=X'\setminus Y$, which is an open and closed subscheme of $X'$ containing $X'_s$.
 We will see that $X$ satisfies the conditions in the lemma. Clearly, $X$ is an affine scheme
 of finite presentation and flat over $k^\circ$, and has a smooth generic fiber $X_\eta$
 whose dimension is at most $d$.
 Since $X$ is an open subscheme of $X'$ containing $X'_s$, 
 the $\varpi$-adic completion of $X$ coincides with that of $X'$. 
 Therefore the $\varpi$-adic completion of $X$ is isomorphic to $\mathcal{X}$.
 Hence we have an open immersion $d(\mathcal{X})\hookrightarrow X_\eta^{\mathrm{ad}}$, where 
 $X_\eta^{\mathrm{ad}}$ denotes the adic space associated to $X_\eta$. As $d(\mathcal{X})$ is $d$-dimensional,
 we conclude that $X_\eta$ is $d$-dimensional. This completes the proof.
\end{prf}

Fix an affine scheme $X$ as in Lemma \ref{lem:Elkik} and denote the structure morphism $X\to \Spec k^\circ$
by $g$.
We put $S=\Spec k^\circ$ and write $s$ (resp.\ $\eta$) for the closed (resp.\ generic) point of $S$.
Consider the following diagram whose rectangles are cartesian:
\[
 \xymatrix{
 \mathcal{X}_s\ar@{=}[r]\ar[d]^{f_s}& X_s\ar[r]^-{i}\ar[d]^{g_s}& X\ar[d]^{g}& X_\eta\ar[l]_-{j}\ar[d]^{g_\eta}\\
 s\ar@{=}[r]& s\ar[r]^-{i}& S& \eta\ar[l]_-{j}\lefteqn{.}
 }
\]
As usual, we put $R\psi_X\Lambda=i^*Rj_*\Lambda$.
By adjointness, we have a natural map $\spp^*\colon \Lambda\to R\psi_X\Lambda$.
On the other hand, since $X_\eta$ is $d$-dimensional, we have the trace map
\[
 Rg_{s!}R\psi_X\Lambda\to R\psi_S Rg_{\eta!}\Lambda\xrightarrow{\Tr_{X_\eta}}R\psi_S\Lambda(-d)[-2d]=\Lambda(-d)[-2d].
\]
The map $R\psi_X\Lambda\to Rg_s^!\Lambda(-d)[-2d]$ obtained by adjointness is denoted by $\cosp^*$.

\begin{lem}\label{lem:comparison}
 There exists a natural isomorphism $R\psi_X\Lambda\to R\Psi_{\!\mathcal{X}}\Lambda$,
 which makes the following diagram commute:
 \[
 \xymatrix{%
 \Lambda\ar[r]^-{\spp^*}\ar@{=}[d]& R\psi_X\Lambda\ar[r]^-{\cosp^*}\ar[d]^-{\cong}& Rg_s^!\Lambda(-d)[-2d]\ar@{=}[d]\\
 \Lambda\ar[r]^-{\spp^*}& R\Psi_{\!\mathcal{X}}\Lambda\ar[r]^-{\cosp^*}& Rf_s^!\Lambda(-d)[-2d]\lefteqn{.}
 }
 \]
\end{lem}

\begin{prf}
 The construction of an isomorphism $R\psi_X\Lambda\to R\Psi_{\!\mathcal{X}}\Lambda$ is due to
 \cite[Theorem 3.5.13]{MR1734903}. We shall recall it briefly.
 Consider the following diagram of sites (see \cite[3.5.12]{MR1734903}):
 \[
  \xymatrix{%
 d(\mathcal{X})_{\et}\ar[r]^-{e}\ar[d]^-{\lambda}& X_{\eta,\et}\ar[d]^-{j}\\
 X_{s,\et}\ar[r]^{i}& X_{\et}\lefteqn{.}
 }
 \]
 As in \cite[Theorem 3.5.13]{MR1734903}, one can construct a natural morphism of functors
 $\phi\colon \lambda^{-1}\circ i^{-1}\to e^{-1}\circ j^{-1}$.
 For every sheaf $L$ on $X_\et$, a map
 \begin{align*}
  R\psi_Xj^*L=i^*Rj_*j^*L&\xrightarrow{\adj}i^*Rj_*Re_*e^*j^*L\xrightarrow{\phi} i^*Ri_*R\lambda_*e^*j^*L\\
  &\xrightarrow{\adj}R\lambda_*e^*j^*L=R\Psi_{\!\mathcal{X}}e^*j^*L
 \end{align*}
 is induced.
 If we put $L=\Lambda$, we get the map $R\psi_X\Lambda\to R\Psi_{\!\mathcal{X}}\Lambda$,
 which is in fact an isomorphism.

 Let us prove the commutativity of the left rectangle. For a sheaf $L$ on $X_\et$,
 we have the following commutative diagram (see the subsequent Lemma \ref{lem:site-general} (iii)):
 \[
 \xymatrix{%
 L\ar[r]^-{\adj}\ar@{=}[d] & R(j\circ e)_*(j\circ e)^*L\ar[r]^-{\phi}& R(i\circ \lambda)_*(j\circ e)^*L\ar@{=}[d]\\
 L\ar[r]^-{\adj} & R(i\circ \lambda)_*(i\circ \lambda)^*L\ar[r]^-{\phi}& R(i\circ \lambda)_*(j\circ e)^*L\lefteqn{.}
 }
 \]
 By the adjointness of $i^*$ and $Ri_*$, we obtain the commutative diagram
 \[
 \xymatrix{%
 i^*L\ar[r]^-{\adj}\ar@{=}[d] & i^*Rj_*Re_*e^*j^*L\ar[r]^-{\phi}& i^*Ri_*R\lambda_*e^*j^*L\ar[r]^-{\adj}& R\lambda_*e^*j^*L\ar@{=}[d]\\
 i^*L\ar[r]^-{\adj} & R\lambda_*\lambda^*i^*L\ar[rr]^-{\phi}&& R\lambda_*e^*j^*L\lefteqn{.}
 }
 \]
 Now put $L=\Lambda$. By definition, the composite of the top row is equal to
 that of $\Lambda\xrightarrow{\spp^*}R\psi_X\Lambda\to R\Psi_{\!\mathcal{X}}\Lambda$.
 On the other hand, by Lemma \ref{lem:site-general} (ii), 
 the map $\lambda^*i^*\Lambda\xrightarrow{\phi}e^*j^*\Lambda$ is the identity on the constant sheaf $\Lambda$.
 Hence the composite of the bottom row equals $\spp^*\colon \Lambda\to R\Psi_{\!\mathcal{X}}\Lambda$.
 This concludes the commutativity of the left rectangle.

 Next we consider the right rectangle. 
 By exactly the same method as in the proof of \cite[Proposition 4.42]{formalnearby},
 we can check that the rectangle in the following diagram is commutative
 (recall that $X^{\mathrm{ad}}_\eta$ denotes the adic space over $\Spa(k,k^\circ)$ associated to $X_\eta$):
 \[
  \xymatrix{%
 H^{2d}_c(\mathcal{X}_s,R\Psi_{\!\mathcal{X}}\Lambda(d))\ar[d]^-{\cong}&& 
 H^{2d}_c(X_s,R\psi_{X}\Lambda(d))\ar[ll]_-{\cong}\ar[d]\\
 H^{2d}_c(d(\mathcal{X}),\Lambda(d))\ar[r]\ar[rd]_-{\Tr_{d(\mathcal{X})}}& H^{2d}_c(X_\eta^{\mathrm{ad}},\Lambda(d))\ar[d]^-{\Tr_{X_\eta^{\mathrm{ad}}}}& 
 H^{2d}_c(X_\eta,\Lambda(d))\ar[l]^-{\cong}\ar[ld]^-{\Tr_{X_\eta}}\\
 &\Lambda\lefteqn{.}
 }
 \]
 The lower left triangle clearly commutes. The lower right triangle is commutative by 
 \cite[Proposition 2.2]{adicLTF}. Hence we obtain the following commutative diagram:
 \[
  \xymatrix{%
 R\Gamma_c(X_s,R\psi_{X}\Lambda(d)[2d])\ar[r]\ar[d]^-{\cong}& R\Gamma_c(X_\eta,\Lambda(d)[2d])\ar[rr]^-{\Tr_{X_\eta}}&&\Lambda\ar@{=}[d]\\
 R\Gamma_c(\mathcal{X}_s,R\Psi_{\!\mathcal{X}}\Lambda(d)[2d])\ar[r]^-{\cong}& R\Gamma_c(d(\mathcal{X}),\Lambda(d)[2d])\ar[rr]^-{\Tr_{d(\mathcal{X})}}&&\Lambda\lefteqn{.}
 }
 \]
 By adjointness, we conclude that the diagram
 \[
 \xymatrix{%
 R\psi_X\Lambda\ar[r]^-{\cosp^*}\ar[d]^-{\cong}& Rg_s^*\Lambda(-d)[-2d]\ar@{=}[d]\\
 R\Psi_{\!\mathcal{X}}\Lambda\ar[r]^-{\cosp^*}& Rf_s^*\Lambda(-d)[-2d]
 }
 \]
 commutes, as desired.
\end{prf}

In the proof of the lemma above, we have used the following general fact, whose proof is immediate:

\begin{lem}\label{lem:site-general}
 Let $\mathcal{C}$, $\mathcal{C}'$ be sites, $f,g\colon \mathcal{C}\to \mathcal{C}'$ morphisms
 of sites, and $\phi\colon f^{-1}\to g^{-1}$ a morphism of functors.
 \begin{enumerate}
  \item The morphism $\phi$ induces morphisms of functors $f^*\to g^*$ and $g_*\to f_*$,
	where $(f^*,f_*)$ (resp.\ $(g^*,g_*)$) is the morphism of toposes induced by $f$ (resp.\ $g$).
  \item For a constant sheaf $L$ on $\mathcal{C}'$, the morphism $\phi\colon f^*L\to g^*L$ is the identity.
  \item The composite of $\id\xrightarrow{\adj}Rf_*f^*\xrightarrow{\phi}Rf_*g^*$ equals
	that of $\id\xrightarrow{\adj}Rg_*g^*\xrightarrow{\phi}Rf_*g^*$.
 \end{enumerate}
\end{lem}

By Lemma \ref{lem:comparison}, Theorem \ref{thm:Gysin-specialization} is reduced to the following lemma:

\begin{lem}\label{lem:Gysin-specialization-algebraic}
 Let $g\colon X\to S=\Spec k^\circ$ be as in Lemma \ref{lem:Elkik}.
 Then, the composite of $\Lambda\xrightarrow{\spp^*}R\psi_X\Lambda\xrightarrow{\cosp^*} Rg_s^!\Lambda(-d)[-2d]$
 is equal to the Gysin map $\Gys_{g_s}\colon \Lambda\to Rg_s^!\Lambda(-d)[-2d]$ with respect to $g_s$.
\end{lem}

\begin{prf}
 By adjointness, it suffices to show that the composite of
 \begin{equation*}
  Rg_{s!}\Lambda(d)[2d]\xrightarrow{\spp^*}Rg_{s!}R\psi_X\Lambda(d)[2d]\to R\psi_S Rg_{\eta!}\Lambda(d)[2d]\xrightarrow{\Tr_{X_\eta}=\Tr_{g_\eta}}R\psi_S\Lambda\tag{$*$}
 \end{equation*}
 is equal to that of $Rg_{s!}\Lambda(d)[2d]\xrightarrow{\Tr_{g_s}}\Lambda \xrightarrow{\cong} R\psi_S\Lambda$.
 As $g\colon X\to S$ is of finite presentation and flat, we have the trace map
 $\Tr_g\colon Rg_!\Lambda(d)[2d]\to \Lambda$ (see \cite[Expos\'e XVIII, Th\'eor\`eme 2.9]{SGA4}). 
 Since the trace map is compatible with base change, it induces the following commutative diagram:
 \[
  \xymatrix{%
 Rg_!\Lambda(d)[2d]\ar[r]^-{\adj}\ar[d]^-{\Tr_g}& Rg_!Rj_*j^*\Lambda(d)[2d]\ar[r]& Rj_*Rg_{\eta!}\Lambda(d)[2d]\ar[d]^-{Rj_*\Tr_{g_\eta}}\\
 \Lambda\ar[r]^-{\adj}& Rj_*j^*\Lambda\ar@{=}[r]& Rj_*\Lambda\lefteqn{.}
 }
 \]
 By taking $i^*$, we obtain the commutativity of the lower rectangle of the following diagram:
 \[
  \xymatrix{%
 Rg_{s!}\Lambda(d)[2d]\ar[r]^-{\spp^*}\ar[d]^-{\cong}& Rg_{s!}R\psi_X\Lambda(d)[2d]\ar[d]^-{\cong}\\
 i^*Rg_!\Lambda(d)[2d]\ar[r]^-{\adj}\ar[dd]^-{i^*\Tr_g}& i^*Rg_!Rj_*j^*\Lambda(d)[2d]\ar[d]\\
 & R\psi_S Rg_{\eta !}\Lambda(d)[2d]\ar[d]^-{R\psi_S\Tr_{g_{\eta}}}\\
 \Lambda\ar[r]& R\psi_S\Lambda\lefteqn{.}
 }
 \]
 The upper rectangle is clearly commutative. Hence the composite of $(*)$ is equal to that of
 \[
  Rg_{s!}\Lambda(d)[2d]\xrightarrow{\cong}i^*Rg_!\Lambda(d)[2d]\xrightarrow{i^*\Tr_g}\Lambda\xrightarrow{\cong}R\psi_S\Lambda.
 \]
 Since the trace map is compatible with base change, the composite of the first two maps equals $\Tr_{g_s}$.
 This concludes the proof.
\end{prf}

\begin{cor}\label{cor:Gysin-specialization}
 Let $Y$ be a purely $d$-dimensional separated smooth scheme of finite type over $\kappa$.
 Assume that a finite surjective morphism $\pi\colon \mathcal{X}_s\to Y$ over $\kappa$ is given.
 Take a decomposition $Y=\coprod_j Y_j$ into connected components, and let $\delta_j\ge 1$ be the generic degree
 of $\pi\vert_{Y_j}\colon \pi^{-1}(Y_j)\to Y_j$. We denote by $\deg \pi$ the locally constant function on $Y$
 that equals $\delta_j$ on $Y_j$. It induces an endomorphism $\Lambda\xrightarrow{\times\deg\pi}\Lambda$
 of a sheaf over $Y$, and thus that of $H^i_c(Y,\Lambda)$. Concretely, it is described as follows:
 \[
  \times\deg \pi\colon H^i_c(Y,\Lambda)=\bigoplus_j H^i_c(Y_j,\Lambda)\xrightarrow{\bigoplus_j (\times\delta_j)}\bigoplus_j H^i_c(Y_j,\Lambda)=H^i_c(Y,\Lambda).
 \]
 Then, we have the following commutative diagram:
 \[
 \xymatrix{%
 H^i_c(Y,\Lambda)\ar[r]^-{\pi^*}\ar[d]_-{\times\deg\pi}& H^i_c(\mathcal{X}_s,\Lambda)\ar[d]^-{\Gys_{f_s}}\ar[r]^-{\spp^*}& H^i_c(\mathcal{X}_s,R\Psi_{\!\mathcal{X}}\Lambda)\ar[r]^-{\cong}\ar@{=}[d]& H^i_c(d(\mathcal{X}),\Lambda)\ar@{=}[d]\\
 H^i_c(Y,\Lambda)\ar[d]&H^i_c(\mathcal{X}_s,Rf_s^!\Lambda(-d)[-2d])\ar[l]_-{\pi_*}\ar[d]& H^i_c(\mathcal{X}_s,R\Psi_{\!\mathcal{X}}\Lambda)\ar[l]_-{\cosp^*}\ar[d]\ar[r]^-{\cong}& H^i_c(d(\mathcal{X}),\Lambda)\ar[d]\\
  H^i(Y,\Lambda)&H^i(\mathcal{X}_s,Rf_s^!\Lambda(-d)[-2d])\ar[l]_-{\pi_*}& H^i(\mathcal{X}_s,R\Psi_{\!\mathcal{X}}\Lambda)\ar[l]_-{\cosp^*}\ar@{=}[r]& H^i(d(\mathcal{X}),\Lambda)\lefteqn{.}
 }
 \]
 Here $\pi_*$ are the maps induced from 
 \[
 \pi_*Rf_s^!\Lambda(-d)[-2d]=R\pi_!R\pi^!Rh^!\Lambda(-d)[-2d]\xleftarrow[\cong]{\Gys_h}R\pi_!R\pi^!\Lambda\xrightarrow{\adj}\Lambda,
 \]
 where $h\colon Y\to \Spec \kappa$ denotes the structure map of $Y$.
\end{cor}

\begin{prf}
 By Theorem \ref{thm:Gysin-specialization}, the upper middle rectangle commutes.
 To see the commutativity of the upper left rectangle, it suffices to prove that the composite of
 \begin{align*}
 \Lambda\xrightarrow{\adj} \pi_*\pi^*\Lambda&=R\pi_!\Lambda\xrightarrow{\Gys_{f_s}}R\pi_!Rf_s^!\Lambda(-d)[-2d]\\
  &=R\pi_!R\pi^!Rh^!\Lambda(-d)[-2d]\xleftarrow[\cong]{\Gys_h}R\pi_!R\pi^!\Lambda\xrightarrow{\adj}\Lambda\tag{$*$}
 \end{align*}
 is equal to the multiplication by $\deg\pi$.
 Since this problem is local on $Y$, we may assume that $Y$ is connected.
 Further, by shrinking $Y$, we may assume that $\pi\colon \mathcal{X}_s\to Y$ is flat
 (\cite[IV, Th\'eor\`eme 6.9.1]{EGA}).
 Then, as the trace map is compatible with composition, we can easily observe that the composite of $(*)$
 is equal to that of
 \[
  \Lambda\xrightarrow{\adj}\pi_*\pi^*\Lambda=\pi_!\pi^*\Lambda\xrightarrow{\Tr_\pi}\Lambda.
 \]
 Hence the claim follows from \cite[Expos\'e XVIII, Th\'eor\`eme 2.9 (Var 4)]{SGA4}.
 
 The commutativity of the lower right rectangle follows from the construction of the isomorphism
 $H^i_c(\mathcal{X}_s,R\Psi_{\!\mathcal{X}}\Lambda)\to H^i_c(d(\mathcal{X}),\Lambda)$
 (see \cite[Definition 4.27]{formalnearby}) and the proof of \cite[Lemma 4.8]{formalnearby}.
 The other rectangles are obviously commutative.
\end{prf}

Note that $H^i_c(d(\mathcal{X}),\Z_\ell)=\varprojlim_m H^i_c(d(\mathcal{X}),\Z/\ell^m\Z)$
by \cite[Theorem 3.1]{MR1626021};
recall that we are assuming that $d(\mathcal{X})$ is smooth, thus locally algebraic by
\cite[(1.7.7)]{MR1734903}. Hence, by taking projective limit, we obtain the specialization maps
\[
 \spp^*\colon H^i_c(\mathcal{X}_s,\Z_\ell)\to H^i_c(d(\mathcal{X}),\Z_\ell),\quad
 \spp^*\colon H^i_c(\mathcal{X}_s,\overline{\Q}_\ell)\to H^i_c(d(\mathcal{X}),\overline{\Q}_\ell)
\]
between $\ell$-adic cohomology.

\begin{thm}\label{thm:injectivity-criterion-finite-level}
 Let $Y$ and $\pi\colon \mathcal{X}_s\to Y$ be as in Corollary \ref{cor:Gysin-specialization}.
 Let $V$ be a subspace of $H^i_c(Y,\overline{\Q}_\ell)$ such that the composite
 $V\hookrightarrow H^i_c(Y,\overline{\Q}_\ell)\to H^i(Y,\overline{\Q}_\ell)$ is
 an injection.
 Then, for any open immersion $d(\mathcal{X})\hookrightarrow Z$ into an adic space $Z$
 which is locally of finite type, separated and taut over $\Spa(k,k^\circ)$, the composite
 \[
  V\hookrightarrow H^i_c(Y,\overline{\Q}_\ell)\xrightarrow{\pi^*} H^i_c(\mathcal{X}_s,\overline{\Q}_\ell)\xrightarrow{\spp^*} H^i_c(d(\mathcal{X}),\overline{\Q}_\ell)\to H^i_c(Z,\overline{\Q}_\ell)
 \]
 is also an injection.
\end{thm}

\begin{prf}
 By Corollary \ref{cor:Gysin-specialization}, we have the following commutative diagram for each integer $m\ge 1$:
 \[
  \xymatrix{%
 H^i_c(Y,\Z/\ell^m\Z)\ar[rr]^-{\spp^*\circ \pi^*}\ar[d]_-{\times \deg\pi}&& H^i_c(d(\mathcal{X}),\Z/\ell^m\Z)\ar[r]\ar@{=}[d]& H^i_c(Z,\Z/\ell^m\Z)\ar@{=}[d]\\
 H^i_c(Y,\Z/\ell^m\Z)\ar[d]&& H^i_c(d(\mathcal{X}),\Z/\ell^m\Z)\ar[ll]_-{\pi_*\circ\cosp^*}\ar[r]\ar[d]& H^i_c(Z,\Z/\ell^m\Z)\ar[d]\\
 H^i(Y,\Z/\ell^m\Z)&&H^i(d(\mathcal{X}),\Z/\ell^m\Z)\ar[ll]_-{\pi_*\circ\cosp^*}& H^i(Z,\Z/\ell^m\Z)\ar[l]\lefteqn{.}
 }
 \]
 Therefore, the composite $\phi\colon H^i_c(Y,\overline{\Q}_\ell)\xrightarrow{\times\deg\pi}H^i_c(Y,\overline{\Q}_\ell)\to H^i(Y,\overline{\Q}_\ell)$ can be decomposed into
 \[
  H^i_c(Y,\overline{\Q}_\ell)\xrightarrow{\spp^*\circ\pi^*} H^i_c(d(\mathcal{X}),\overline{\Q}_\ell)\xrightarrow{(*)} \Bigl(\varprojlim_m H^i_c(Z,\Z/\ell^m\Z)\Bigr)\otimes_{\Z_\ell}\overline{\Q}_\ell\to H^i(Y,\overline{\Q}_\ell).
 \]
 Clearly $(*)$ is the composite of 
 \[
  H^i_c(d(\mathcal{X}),\overline{\Q}_\ell)\to H^i_c(Z,\overline{\Q}_\ell)\to\Bigl(\varprojlim_m H^i_c(Z,\Z/\ell^m\Z)\Bigr)\otimes_{\Z_\ell}\overline{\Q}_\ell.
 \]
 Hence $\phi$ has a decomposition
 \[
  H^i_c(Y,\overline{\Q}_\ell)\xrightarrow{\spp^*\circ\pi^*} H^i_c(d(\mathcal{X}),\overline{\Q}_\ell)\to
 H^i_c(Z,\overline{\Q}_\ell)\to H^i(Y,\overline{\Q}_\ell).
 \]
 Since $\phi\vert_V\colon V\hookrightarrow H^i_c(Y,\overline{\Q}_\ell)\xrightarrow{\phi} H^i(Y,\overline{\Q}_\ell)$
 is the composite of
 \[
  V\hookrightarrow H^i_c(Y,\overline{\Q}_\ell)\to H^i(Y,\overline{\Q}_\ell)\xrightarrow{\times\deg\pi}H^i(Y,\overline{\Q}_\ell)
 \]
 and $\times\deg\pi$ is an isomorphism, $\phi\vert_V$ is an injection by the assumption on $V$.
 Therefore, the composite
 \[
 V\hookrightarrow H^i_c(Y,\overline{\Q}_\ell)\xrightarrow{\spp^*\circ\pi^*} H^i_c(d(\mathcal{X}),\overline{\Q}_\ell)\to
 H^i_c(Z,\overline{\Q}_\ell)
 \]
 is also an injection, as desired.
\end{prf}

\section{Specialization map: the case of infinite level}\label{sec:infinite-level}
Let $k$, $k^\circ$, $\mathfrak{m}$ and $\kappa$ be as in the previous section.
Assume that the characteristic of $\kappa$ is $p>0$. Fix a non-zero element $\varpi\in\mathfrak{m}$.

In this section, we consider a flat $k^\circ$-algebra $A$ which is $\varpi$-adically complete.
We denote by $\widetilde{A}$ the integral closure of $A$ in $A[1/\varpi]$.
Then, $(A[1/\varpi],\widetilde{A})$ is an affinoid $(k,k^\circ)$-algebra (for the topology on $A[1/\varpi]$,
see \cite[Examples 1.1 (iv)]{MR1207303}).
Suppose that $A$ is equipped with an action over $k^\circ$ of a profinite group $K_0$.
Then, $\widetilde{A}$ is a $K_0$-stable subring of $A[1/\varpi]$.

We will assume two conditions on $A$ and the action of $K_0$.
The first is that the affinoid $(k,k^\circ)$-algebra $(A[1/\varpi],\widetilde{A})$
is obtained as a limit of a tower of smooth affinoid $(k,k^\circ)$-algebras.
The precise statement is as follows:

\begin{assump}\label{assump:Hecke}
 For each open normal subgroup $K$ of $K_0$, we are given a complete f-adic $k$-algebra $B_K$
 of topologically finite type endowed with an action of $K_0/K$.
 For open normal subgroups $K$, $K'$ with $K'\subset K$,
 we are given a continuous $K_0$-equivariant $k$-homomorphism $B_K\to B_{K'}$
 so that $\{B_K\}_{K\lhd K_0}$ becomes an inductive system. 
 Further we impose the following conditions:
 \begin{enumerate}
  \item[(a)] For an open normal subgroup $K$ of $K_0$, the adic space $X_K=\Spa(B_K,B_K^\circ)$ is smooth over $\Spa(k,k^\circ)$.
  \item[(b)] For open normal subgroups $K$, $K'$ of $K_0$ with $K'\subset K$, 
	     the transition map $X_{K'}\to X_K$ is a finite \'etale Galois covering with Galois group $K/K'$.
	     Note that this implies that $X_{K'}\to X_K$ is surjective and
	     $B_K\to (B_{K'})^{K/K'}$ is an isomorphism.
  \item[(c)] There exists a $K_0$-equivariant homomorphism 
	     $\{(B_K,B_K^\circ)\}_{K\lhd K_0}\to (A[1/\varpi],\widetilde{A})$ of 
	     inductive systems of affinoid $(k,k^\circ)$-algebras
	     (here $(A[1/\varpi],\widetilde{A})$ means the constant inductive system),
	     such that
	     \begin{itemize}
	      \item the induced continuous map
		    $\Spa(A[1/\varpi],\widetilde{A})\to \varprojlim_{K\lhd K_0}\Spa(B_K,B_K^\circ)$
		    is a homeomorphism, and
	      \item the induced homomorphism $\varinjlim_K B_K\to A[1/\varpi]$ has dense image.
	     \end{itemize}
 \end{enumerate}
\end{assump}

The second assumption is that the reduction of $A$ is the perfection of a finitely generated $\kappa$-algebra.
More precisely,

\begin{assump}\label{assump:fin-type-perfection}
 There exist a finitely generated $\kappa$-algebra $R$ and an isomorphism
 \[
  \varinjlim_{x\mapsto x^p}R\xrightarrow{\cong} A/\mathfrak{m}A
 \]
 of $\kappa$-algebras. We fix such an $R$ and an isomorphism, and put $Y=\Spec R$.
\end{assump}

Under these assumptions, we will compare two cohomology groups
$H^i_c(Y,\overline{\Q}_\ell)$ and $\varinjlim_{K}H^i_c(X_K,\overline{\Q}_\ell)$.
We begin with basic properties of the inductive system $\{B_K\}_{K\lhd K_0}$.

\begin{lem}\label{lem:limit-surj}
 For each open normal subgroup $K$ of $K_0$, the map
 $\Spa(A[1/\varpi],\widetilde{A})\to \Spa(B_K,B_K^\circ)=X_K$ is surjective.
\end{lem}

\begin{prf}
 By Assumption \ref{assump:Hecke}, it suffices to show the surjectivity
 of the map 
 \[
  \varprojlim_{K'\lhd K_0}X_{K'}=\varprojlim_{K'\lhd K_0, K'\subset K}X_{K'}\to X_K.
 \]
 This follows from the surjectivity and the finiteness of $X_{K'}\to X_K$;
 recall that the filtered projective limit of non-empty finite sets is non-empty.
\end{prf}

\begin{lem}\label{lem:B_K-hom-inj}
 The homomorphism $B_K\to A[1/\varpi]$ is injective.
\end{lem}

\begin{prf}
 First note that $B_K$ is reduced, as $\Spa(B_K,B_K^\circ)$ is assumed to be smooth over $\Spa(k,k^\circ)$. 
 Therefore we have $0=\sqrt{0}=\bigcap_{\mathfrak{n}\subset B_K}\mathfrak{n}$,
 where $\mathfrak{n}$ runs through  maximal ideals of $B_K$
 (recall that $B_K$ is a Jacobson ring; see \cite[Chapter 0, Proposition 9.3.10]{Fujiwara-Kato}).

 We denote the homomorphism $B_K\to A[1/\varpi]$ by $\phi$. 
 Take an arbitrary maximal ideal $\mathfrak{n}$ of $B_K$,
 and let $v_\mathfrak{n}\in\Spa(B_K,B_K^\circ)$ be the classical point corresponding to $\mathfrak{n}$.
 By Lemma \ref{lem:limit-surj}, 
 we can find $\widetilde{v}_\mathfrak{n}\in \Spa(A[1/\varpi],\widetilde{A})$ satisfying
 $v_\mathfrak{n}=\widetilde{v}_\mathfrak{n}\circ\phi$.
 For $a\in \Ker\phi$, we have $v_\mathfrak{n}(a)=\widetilde{v}_\mathfrak{n}(\phi(a))=\widetilde{v}_\mathfrak{n}(0)=0$.
 Hence $\Ker\phi\subset \supp v_\mathfrak{n}=\mathfrak{n}$.
 Therefore we conclude that $\Ker\phi\subset \bigcap_{\mathfrak{n}\subset B_K}\mathfrak{n}=0$,
 namely, $\phi$ is an injection.
\end{prf}

By this lemma, we regard $B_K$ as a $k$-subalgebra of $A[1/\varpi]$.

\begin{defn}
 We put $B'=\bigcup_{K\lhd K_0} B_K\subset A[1/\varpi]$, $A'=A\cap B'$ and $\widetilde{A}'=\widetilde{A}\cap B'$.
 Note that $B'$ is a dense $K_0$-stable $k$-subalgebra of $A[1/\varpi]$, while $A'$ and $\widetilde{A'}$ are 
 $K_0$-stable $k^\circ$-subalgebras of $B'$.
 For an open normal subgroup $K$ of $K_0$, we put $A_K=A\cap B_K$ and $\widetilde{A}_K=\widetilde{A}\cap B_K$.
 These are open $K_0$-stable $k^\circ$-subalgebras of $B_K$.
\end{defn}

\begin{lem}\label{lem:A'-special-fiber}
 The inclusion $A'\hookrightarrow A$ induces isomorphisms $A'/\varpi^mA'\xrightarrow{\cong} A/\varpi^mA$
 for every integer $m\ge 0$ and $A'/\mathfrak{m}A'\xrightarrow{\cong}A/\mathfrak{m}A$.
\end{lem}

\begin{prf}
 Let $m\ge 0$ be an integer.
 Since $A'$ is dense in $A$, we have $A=A'+\varpi^mA$. Thus the map 
 $A'/\varpi^mA'\to A/\varpi^mA$ is surjective. Let us show the injectivity.
 Take an element $x\in A'\cap \varpi^mA$ and write $x=\varpi^my$ with $y\in A$.
 Since $B'$ is a $k$-subalgebra of $A[1/\varpi]$, we have $y=\varpi^{-m}x\in \varpi^{-m}A'\subset B'$.
 Hence $y$ lies in $A\cap B'=A'$ and $x$ lies in $\varpi^mA'$. This means that 
 $A'/\varpi^mA'\to A/\varpi^mA$ is injective.

 In particular, the map $A'/\varpi A'\to A/\varpi A$ is an isomorphism. By taking the base change
 $(-)\otimes_{k^\circ/\varpi k^\circ}k^\circ/\mathfrak{m}$, we conclude that the homomorphism
 $A'/\mathfrak{m}A'\to A/\mathfrak{m}A$ is also an isomorphism.
\end{prf}

\begin{lem}\label{lem:B'-fixed}
 For an open normal subgroup $K$ of $K_0$, we have $(B')^K=B_K$, $(A')^K=A_K$ and 
 $(\widetilde{A}')^K=\widetilde{A}_K$.
\end{lem}

\begin{prf}
 For an open normal subgroup $K'$ of $K_0$ contained in $K$, we have $(B_{K'})^K=B_K$ by assumption.
 By taking the union with respect to such $K'$, we conclude that $(B')^K=B_K$.
 From this equality, we can deduce $(A')^K=A\cap (B')^K=A\cap B_K=A_K$ and
 $(\widetilde{A}')^K=\widetilde{A}\cap (B')^K=\widetilde{A}\cap B_K=\widetilde{A}_K$.
\end{prf}

\begin{lem}\label{lem:tilde-integral}
 \begin{enumerate}
  \item The ring $\widetilde{A}'$ is integral over its subring $A'$.
  \item For an open normal subgroup $K$ of $K_0$, $\widetilde{A}_K$ is integral over its subring $A_K$.
  \item For an open normal subgroup $K$ of $K_0$, $\widetilde{A}_K$ coincides with $B_K^\circ$.
 \end{enumerate}
\end{lem}

\begin{prf}
 We prove (i). Take $x\in \widetilde{A}'\subset \widetilde{A}$. Since $\widetilde{A}$ is integral over $A$,
 there exist $n\ge 1$ and $a_1,\ldots,a_n\in A$ such that $x^n+a_1x^{n-1}+\cdots+a_n=0$.
 As $A\subset \widetilde{A}\subset A[1/\varpi]$, we can find an integer $m\ge 0$ such that
 $\varpi^m x^i\in A$ for every $0\le i\le n-1$. Since $B'$ is dense in $A[1/\varpi]$, 
 we have $(a_i+\varpi^mA)\cap B'\neq \varnothing$; in other words, for each $i$ there exists
 $b_i\in B'$ such that $b_i-a_i\in\varpi^mA$. Note that $b_i\in A\cap B'=A'$.

 For such $b_i$'s, we have
 \[
  x^n+b_1x^{n-1}+\cdots+b_n=(b_1-a_1)x^{n-1}+\cdots+(b_n-a_n).
 \]
 The left hand side lies in $B'$, while the right hand side lies in $A$.
 If we put $c=(b_1-a_1)x^{n-1}+\cdots+(b_n-a_n)$, it is an element of $A\cap B'=A'$ and
 $x$ is a root of the monic polynomial
 $T^n+b_1T^{n-1}+\cdots+b_n-c\in A'[T]$. Hence $x$ is integral over $A'$.

 Next we prove (ii). By (i), it suffices to show that $A'$ is integral over $A_K$.
 Take an arbitrary element $a\in A'$ and an open normal subgroup $K'$ of $K_0$ contained in $K$
 such that $a\in B_{K'}$. Then $a\in A'\cap B_{K'}=A_{K'}$. The element $a$ is a root of the monic polynomial
 $\prod_{g\in K/K'}(T-g(a))$, whose coefficients belong to $(A_{K'})^K=A_K$. Therefore $a$ is integral
 over $A_K$.

 Finally we prove (iii). Since we have a homomorphism $(B_K,B_K^\circ)\to (A[1/\varpi],\widetilde{A})$
 of affinoid $(k,k^\circ)$-algebras, we have $B_K^\circ \subset \widetilde{A}\cap B_K=\widetilde{A}_K$.
 Let us prove the reverse inclusion $\widetilde{A}_K\subset B_K^\circ$. Take $a\in\widetilde{A}_K$.
 By \cite[Lemma 3.3 (i)]{MR1207303}, we have
 \[
  B_K^\circ=\{x\in B_K\mid v(x)\le 1\ \text{for every $v\in \Spa(B_K,B_K^\circ)$}\}.
 \]
 Therefore, it suffices to observe that $v(a)\le 1$ for every $v\in \Spa(B_K,B_K^\circ)$.
 By Lemma \ref{lem:limit-surj}, there exists 
 $\widetilde{v}\in \Spa(A[1/\varpi],\widetilde{A})$ such that $\widetilde{v}\vert_{B_K}=v$.
 As $a\in \widetilde{A}_K\subset \widetilde{A}$, we have $v(a)=\widetilde{v}(a)\le 1$, as desired.
\end{prf}

\begin{cor}\label{cor:integral-extension}
 For an open normal subgroup $K$ of $K_0$, $B_K^\circ$ contains $A_K$ over which $B_K^\circ$ is integral.
\end{cor}

\begin{prf}
 Clear from Lemma \ref{lem:tilde-integral} (ii), (iii).
\end{prf}

Recall that a $k^\circ$-algebra $C$ is said to be topologically finitely generated if
it is isomorphic to a quotient of $k^\circ\langle T_1,\ldots,T_n\rangle$ for some $n\ge 0$
(see \cite[Definition 8.4.1 and Proposition 8.4.4]{Fujiwara-Kato}).
Such $C$ is $\varpi$-adically complete. 
We say that $c_1,\ldots,c_{n'}\in C$ are topological generators of $C$ if 
the $k^\circ$-homomorphism $k^\circ\langle T_1,\ldots,T_{n'}\rangle\to C$
given by $T_i\mapsto c_i$ is surjective.

We would like to write the $k^\circ$-algebra $A'$ as an inductive limit of
topologically finitely generated $k^\circ$-algebras.
First we focus on each $B_K$. In the following proposition, we slightly change the notation for simplicity.

\begin{prop}\label{prop:open-limit}
 Let $B$ be a complete f-adic $k$-algebra of topologically finite type.
 Let $A$ be an open subring of $B^\circ$ such that $B^\circ$ is integral over $A$.
 We write $\mathcal{D}_A$ for the set consisting of topologically finitely generated $k^\circ$-subalgebras
 of $A$ which are rings of definition of $A$ (or equivalently, $B$).
 Then, the following hold:
 \begin{enumerate}
  \item The set $\mathcal{D}_A$ is a filtered ordered set with respect to inclusions.
  \item We have $A=\bigcup_{A_0\in\mathcal{D}_A}A_0$.
  \item For $A_0\in \mathcal{D}_A$, $A$ and $B^\circ$ are integral over $A_0$.
 \end{enumerate} 
\end{prop}

\begin{prf}
 First we prove that $\mathcal{D}_A$ is non-empty. Since $B$ is of topologically finite type,
 we can take a surjective continuous $k$-homomorphism $\phi\colon k\langle T_1,\ldots,T_n\rangle\to B$. 
 By the open mapping theorem, $\phi$ is an open map. Moreover, by \cite[Proposition 1.10]{MR1207303},
 $\phi$ is an adic map. Hence the image $B_0$ of $k^\circ\langle T_1,\ldots,T_n\rangle$ under $\phi$ is 
 a topologically finitely generated $k^\circ$-subalgebra of $B$ which is a ring of definition of $B$.
 Therefore $\varpi B_0$ is an ideal of definition of $B_0$ (see \cite[Proposition 1.5 (ii)]{MR1207303}),
 and $B_0$ is $\varpi$-adically complete.
 Further, by \cite[proof of Lemma 4.4]{MR1306024}, $B^\circ$ is the integral closure of $B_0$ inside $B$.

 Put $b_i=\phi(T_i)\in B_0\subset B^\circ$. Since $B^\circ$ is integral over $A$,
 for each $1\le i\le n$ there exists a monic polynomial $f_i\in A[T]$ such that $f_i(b_i)=0$.
 Since $A$ is open, there exists an integer $m\ge 1$ such that $\varpi^mB_0\subset A$.
 Let $A_0$ be the $k^\circ$-subalgebra of $A$ generated by the coefficients of $f_1,\ldots,f_n$ and 
 the elements of $\varpi^mB_0$. Since $A_0$ contains $\varpi^mB_0$, it is an open subring of $B$.
 As the coefficients of $f_1,\ldots,f_n$ are power-bounded in $B$, $A_0$ is bounded in $B$.
 Hence $A_0$ is a ring of definition of $B$. In particular, $A_0$ is $\varpi$-adically complete.
 We will show that $A_0\in \mathcal{D}_A$.

 Let us prove that $B_0\subset A_0[b_1,\ldots,b_n]\subset B^\circ$. For every $x\in B_0$, there exists
 $f\in k^\circ[T_1,\ldots,T_n]$ such that $x-f(b_1,\ldots,b_n)\in \varpi^mB_0$.
 Thus $x\in f(b_1,\ldots,b_n)+\varpi^mB_0\subset A_0[b_1,\ldots,b_n]$.
 The inclusion $A_0[b_1,\ldots,b_n]\subset B^\circ$ is clear.
 We put $B'_0=A_0[b_1,\ldots,b_n]$. As $b_1,\ldots,b_n$ are power-bounded in $B$, $B_0'$ is also
 a ring of definition of $B$. 
 Since $b_i$ is integral over $A_0$ by construction, $B'_0$ is a finite $A_0$-algebra.
 Note also that $\varpi^mB_0\subset A_0$ implies that $\varpi^mB_0'\subset A_0$.

 Now we prove $A_0\in\mathcal{D}_A$. The $k^\circ$-algebra $A_0/\varpi^mB_0$ is generated by the coefficients
 of $f_1,\ldots,f_n$. Therefore, its quotient $A_0/\varpi^mB_0'$ is also finitely generated over $k^\circ$.
 Let us observe that $A_0/\varpi^{m+1}B_0'$ is finitely generated over $k^\circ$.
 For simplicity, we set $R=A_0/\varpi^{m+1}B_0'$ and $I=\varpi^mB_0'/\varpi^{m+1}B_0'$.
 Since $B_0'$ is finitely generated as an $A_0$-module, $I$ is a finitely generated ideal of $R$.
 As $m\ge 1$, we have $I^2=0$.
 Take $a_1,\ldots,a_r\in R$ so that their images generate $R/I$ as a $k^\circ$-algebra, and generators
 $x_1,\ldots,x_s\in I$ as an $R$-module. Then it is immediate to see that $R$ is generated by 
 $a_1,\ldots,a_r,x_1,\ldots,x_s$ as a $k^\circ$-algebra.
 Since $\varpi^{m+1}B_0'\subset \varpi A_0$, we conclude that $A_0/\varpi A_0$ is a finitely generated
 $k^\circ$-algebra.
 By \cite[Chapter 0, Proposition 8.4.2]{Fujiwara-Kato}, this means that $A_0$ is topologically finitely generated
 over $k^\circ$, as desired.

 Next we prove (iii). For $A_0\in\mathcal{D}_A$, we can take a surjection
 $\phi\colon k\langle T_1,\ldots,T_n\rangle\to B$ above so that 
 $\phi(k^\circ\langle T_1,\ldots,T_n\rangle)=A_0$. As mentioned above, $B^\circ$ is integral over $A_0$.
 Hence $A$ is also integral over $A_0$.

 We prove (i). Take $A_0,A_0'\in \mathcal{D}_A$ and topological generators $a'_1,\ldots,a'_n$ of $A_0'$
 over $k^\circ$. Then, by the same argument as in the third paragraph of this proof, we can show that
 $A_0'\subset A_0[a'_1,\ldots,a'_n]$. Clearly $A_0[a'_1,\ldots,a'_n]$ is open and bounded in $B$,
 hence a ring of definition of $B$.
 In particular, $A_0[a'_1,\ldots,a'_n]$ is $\varpi$-adically complete.
 On the other hand, by (iii), $a'_i\in A$ is integral over $A_0$.
 Therefore $A_0[a'_1,\ldots,a'_n]$ is a finite $A_0$-algebra, hence topologically finitely generated
 $k^\circ$-algebra. 
 Thus $A_0[a'_1,\ldots,a'_n]$ is an element of $\mathcal{D}_A$ containing $A_0$ and $A_0'$.

 Finally we prove (ii). Fix $A_0\in\mathcal{D}_A$. For $a\in A$, consider the subring $A_0[a]$ of $A$.
 Since $a\in A\subset B^\circ$, $A_0[a]$ is a ring of definition of $B$. 
 In particular, $A_0[a]$ is $\varpi$-adically complete.
 By (iii), $a$ is integral over $A_0$. Therefore $A_0[a]$ is finite over $A_0$, and thus
 topologically finitely generated over $k^\circ$. 
 Hence $A_0[a]$ is an element of $\mathcal{D}_A$ containing $a$.
 This completes the proof.
\end{prf}

We return to the original setting. For an open normal subgroup $K$ of $K_0$,
we simply write $\mathcal{D}_K$ for $\mathcal{D}_{A_K}$.

\begin{lem}\label{lem:smaller-K}
 For open normal subgroups $K$, $K'$ of $K_0$ with $K'\subset K$ and $A_0\in \mathcal{D}_K$,
 there exists an element of $\mathcal{D}_{K'}$ containing $A_0$.
\end{lem}

\begin{prf}
 Take an arbitrary element $A'_0$ of $\mathcal{D}_{K'}$ and topological generators $a_1,\ldots,a_n$
 of $A_0$ over $k^\circ$. 
 Clearly $A_0'[a_1,\ldots,a_n]$ is a subring of $A_{K'}$.
 Since $a_1,\ldots,a_n$ are power-bounded in $B_K$ and $B_K\to B_{K'}$ is adic, they are power-bounded in $B_{K'}$.
 Therefore $A_0'[a_1,\ldots,a_n]$ is a ring of definition of $B_{K'}$.
 By Proposition \ref{prop:open-limit} (iii), $a_i\in A_K\subset A_{K'}$ is integral over $A_0'$.
 Therefore $A_0'[a_1,\ldots,a_n]$ is a finite $A_0'$-algebra. In particular, it is topologically finitely
 generated over $k^\circ$. Thus $A_0'[a_1,\ldots,a_n]$ belongs to $\mathcal{D}_{K'}$.
 On the other hand, since $A_0'[a_1,\ldots,a_n]\cap B_K$ is open in $B_K$, there exists $m\ge 1$
 such that $\varpi^mA_0\subset A_0'[a_1,\ldots,a_n]$. 
 Therefore, by the same way as in the third paragraph of the proof of Proposition \ref{prop:open-limit},
 we can check that $A_0\subset A_0'[a_1,\ldots,a_n]$. This concludes the proof.
\end{prf}

Put $\mathcal{D}=\coprod_{K\lhd K_0}\mathcal{D}_K$. We define a partial order on $\mathcal{D}$
as follows: for $A_0\in \mathcal{D}_K$ and $A'_0\in \mathcal{D}_{K'}$, $A_0\le A'_0$ if
$K'\subset K$ and $A_0\subset A_0'$.

\begin{cor}\label{cor:D-filtered}
 This makes $\mathcal{D}$ a filtered ordered set. 
\end{cor}

\begin{prf}
 Clear from Proposition \ref{prop:open-limit} (i) and Lemma \ref{lem:smaller-K}.
\end{prf}

\begin{cor}\label{cor:A'-limit}
 We have 
 \[
  A'=\bigcup_{A_0\in \mathcal{D}}A_0\cong \varinjlim_{A_0\in \mathcal{D}}A_0,\quad
 A/\mathfrak{m}A\cong \varinjlim_{A_0\in \mathcal{D}}A_0/\mathfrak{m}A_0.
 \]
\end{cor}

\begin{prf}
 The first claim is clear from Proposition \ref{prop:open-limit} (ii) and Corollary \ref{cor:D-filtered}.
 As tensor product commutes with inductive limit, we have 
 $A'/\mathfrak{m}A'\cong \varinjlim_{A_0\in \mathcal{D}}A_0/\mathfrak{m}A_0$.
 Therefore the second follows from Lemma \ref{lem:A'-special-fiber}.
\end{prf}

\begin{lem}\label{lem:transition-finite}
 Let $A_0$, $A_0'$ be elements of $\mathcal{D}$ such that $A_0\le A_0'$. 
 Then, the inclusion $A_0\hookrightarrow A_0'$ is finite and continuous.
\end{lem}

\begin{prf}
 Let $K$ and $K'$ be open normal subgroups of $K_0$ such that $A_0\in \mathcal{D}_K$ and $A_0'\in\mathcal{D}_{K'}$.
 Since $(A_{K'})^{K/K'}=A_K$, $A_{K'}$ is integral over $A_K$. By Proposition \ref{prop:open-limit} (iii),
 $A_K$ is integral over $A_0$. As $A_0\subset A_0'\subset A_{K'}$, we conclude that $A_0'$ is integral over $A_0$.
 In particular, $A_0'/\varpi A_0'$ is integral over $A_0/\varpi A_0$.
 On the other hand, as $A_0'$ is a topologically finitely generated $k^\circ$-algebra,
 $A_0'/\varpi A_0'$ is finitely generated over $k^\circ$. In particular, 
 $A_0'/\varpi A_0'$ is finitely generated over $A_0/\varpi A_0$.
 Thus $A_0'/\varpi A_0'$ is finite over $A_0/\varpi A_0$.
 By \cite[Chapter 0, Proposition 7.2.4]{Fujiwara-Kato}, we conclude that $A_0'$ is finite over $A_0$.

 Since the topology of $A_0$ and $A'_0$ are $\varpi$-adic (\cite[Proposition 1.5 (ii)]{MR1207303}),
 the $k^\circ$-homomorphism $A_0\to A'_0$ is clearly continuous.
\end{prf}

Recall that we have fixed a $\kappa$-isomorphism $\varinjlim_{x\mapsto x^p}R\cong A/\mathfrak{m}A$
in Assumption \ref{assump:fin-type-perfection}.
In particular, we are given a $\kappa$-homomorphism $R\to A/\mathfrak{m}A$.

\begin{defn}\label{defn:C}
 Let $\mathcal{C}$ be the set consisting of pairs $(A_0,\pi)$, where
 \begin{itemize}
  \item $A_0$ is an element of $\mathcal{D}$, and
  \item $\pi\colon R\to A_0/\mathfrak{m}A_0$ is a finite $\kappa$-homomorphism such that
	the composite of $\pi$ and the homomorphism $A_0/\mathfrak{m}A_0\to A/\mathfrak{m}A$
	induced from the inclusion $A_0\hookrightarrow A'\hookrightarrow A$ is equal to
	the fixed homomorphism $R\to A/\mathfrak{m}A$.
 \end{itemize}
 It is equipped with the partial order induced from that of $\mathcal{D}$.
 Namely, $(A_0,\pi)\le (A_0',\pi')$ if and only if $A_0\le A_0'$ and the composite
 $R\xrightarrow{\pi} A_0/\mathfrak{m}A_0\to A'_0/\mathfrak{m}A'_0$ equals $\pi'$.
 By Corollaries \ref{cor:D-filtered}, \ref{cor:A'-limit} and Lemma \ref{lem:transition-finite},
 this order is filtered
 (recall that $R$ is finitely generated over $\kappa$).
\end{defn}

\begin{lem}\label{lem:property-C}
 \begin{enumerate}
  \item The set $\mathcal{C}$ is non-empty.
  \item For $(A_0,\pi)\in\mathcal{C}$, the morphism $\Spec A_0/\mathfrak{m}A_0\to \Spec R$ induced by $\pi$
	is a homeomorphism on the underlying topological spaces.
 \end{enumerate}
\end{lem}

\begin{prf}
 We prove (i). Since $A/\mathfrak{m}A=\varinjlim_{A_0\in\mathcal{D}}A_0/\mathfrak{m}A_0$
 and $R$ is finitely presented over $\kappa$, we can find $A_1\in \mathcal{D}$ and a $\kappa$-homomorphism
 $R\to A_1/\mathfrak{m}A_1$ such that the composite $R\to A_1/\mathfrak{m}A_1\to A/\mathfrak{m}A$ is
 the fixed one. On the other hand, Since $A/\mathfrak{m}A=\varinjlim_{x\to x^p}R$, there exist
 a factorization $A_1/\mathfrak{m}A_1\to R\to A/\mathfrak{m}A$ of the homomorphism 
 $A_1/\mathfrak{m}A_1\to A/\mathfrak{m}A$ and an integer $m\ge 1$ such that the following diagram commutes:
 \[
 \xymatrix{%
 R\ar[rr]^-{x\mapsto x^{p^m}}\ar[rd]&& R\ar[r]^-{(*)}& A/\mathfrak{m}A\lefteqn{.}\\
 & A_1/\mathfrak{m}A_1\ar[ru]
 }
 \]
 Apply the same argument to the homomorphism $(*)$, we can find $A_0\in\mathcal{D}$ with
 $A_1\le A_0$ and a homomorphism $R\to A_0/\mathfrak{m}A_0$ such that the right triangle
 of the following diagram commutes:
 \[
 \xymatrix{%
 R\ar[rr]^-{x\mapsto x^{p^m}}\ar[rd]&& R\ar[rr]\ar[rd]_-{(1)}&& A/\mathfrak{m}A\lefteqn{.}\\
 & A_1/\mathfrak{m}A_1\ar[ur]\ar[rr]_-{(2)}&& A_0/\mathfrak{m}A_0\ar[ru]
 }
 \]
 Replacing $A_0$ by larger one, we can also make the lower triangle commute.

 By Lemma \ref{lem:transition-finite}, the homomorphism (2) is finite. Therefore, the homomorphism (1) is
 also finite. If we write $\pi$ for the composite 
 $R\xrightarrow{x\mapsto x^{p^m}}R\xrightarrow{(1)}A_0/\mathfrak{m}A_0$, it is a finite $\kappa$-homomorphism,
 and $(A_0,\pi)$ belongs to $\mathcal{C}$.

 Next we prove (ii). Since $\Spec A_0/\mathfrak{m}A_0\to \Spec R$ is finite and therefore closed,
 it suffices to show that it is a bijection.
 By the same argument as above, we can construct the commutative diagram
 \[
 \xymatrix{%
 R\ar[rr]^-{x\mapsto x^{p^{m'}}}\ar[rd]&& R\ar[rr]\ar[rd]&& A/\mathfrak{m}A\lefteqn{,}\\
 & A_0/\mathfrak{m}A_0\ar[ur]\ar[rr]&& A_2/\mathfrak{m}A_2\ar[ru]
 }
 \]
 where $m'\ge 1$ is an integer and $A_2\in\mathcal{D}$ with $A_0\le A_2$.
 We will observe that $\Spec A_2/\mathfrak{m}A_2\to \Spec A_0/\mathfrak{m}A_0$ is surjective.
 Take open normal subgroups $K$, $K'$ of $K_0$ such that $A_0\in\mathcal{D}_K$ and $A_2\in\mathcal{D}_{K'}$.
 Put $\mathcal{X}=\Spf A_0$ and $\mathcal{X}'=\Spf A_2$. 
 By Lemma \ref{lem:transition-finite}, a morphism $\mathcal{X}'\to \mathcal{X}$ is induced.
 By Proposition \ref{prop:open-limit} (iii),
 we have $d(\mathcal{X})=\Spa(B_K,B_K^\circ)=X_K$ and $d(\mathcal{X}')=\Spa(B_{K'},B_{K'}^\circ)=X_{K'}$.
 Further, we have the following commutative diagram of topological spaces
 (see \cite[Proposition 1.9.1]{MR1734903}):
 \[
  \xymatrix{%
 \mathcal{X}'_{\mathrm{red}}\ar[d]& d(\mathcal{X}')\ar[l]_-{\lambda_{\mathcal{X}'}}\ar[d]\ar@{=}[r]&X_{K'}\ar[d]\\
 \mathcal{X}_{\mathrm{red}}& d(\mathcal{X})\ar[l]_-{\lambda_{\mathcal{X}}}\ar@{=}[r]&X_K\lefteqn{.}
 }
 \]
 By assumption the map $X_{K'}\to X_K$ is surjective. On the other hand, as $\mathcal{X}$ is flat over $k^\circ$,
 the map $\lambda_{\mathcal{X}}$ is surjective 
 (see \cite[Chapter 2, Proposition 3.1.5, Theorem A.4.7]{Fujiwara-Kato}).
 Hence the map $\mathcal{X}'_{\mathrm{red}}\to \mathcal{X}_{\mathrm{red}}$ is also surjective.
 As topological spaces, $\mathcal{X}_{\mathrm{red}}$ (resp.\ $\mathcal{X}'_{\mathrm{red}}$) is identified with
 $\Spec A_0/\mathfrak{m}A_0$ (resp.\ $\Spec A_2/\mathfrak{m}A_2$). 
 Therefore the map $\Spec A_2/\mathfrak{m}A_2\to \Spec A_0/\mathfrak{m}A_0$ is surjective, as desired.

 By this surjectivity, the map $\Spec R\to \Spec A_0/\mathfrak{m}A_0$ is also surjective.
 On the other hand, as the composite $\Spec R\to \Spec A_0/\mathfrak{m}A_0\to \Spec R$ is the identity
 on the underlying topological space, $\Spec R\to \Spec A_0/\mathfrak{m}A_0$ is injective.
 Hence $\Spec R\to \Spec A_0/\mathfrak{m}A_0$ is bijective, and so is $\Spec A_0/\mathfrak{m}A_0\to \Spec R$,
 as desired.
\end{prf}

Now we can construct the specialization map
$\spp^*\colon H^i_c(Y,\overline{\Q}_\ell)\to \varinjlim_K H^i_c(X_K,\overline{\Q}_\ell)$.

\begin{defn}\label{defn:specialization-map-infinite-level}
 Take $(A_0,\pi)\in\mathcal{C}$ and put $\mathcal{X}=\Spf A_0$.
 Then, by Proposition \ref{prop:open-limit} (iii), we have $d(\mathcal{X})=\Spa(B_K,B_K^\circ)=X_K$,
 where $K$ is the open normal subgroup of $K_0$ such that $A_0\in\mathcal{D}_K$.
 We define $\spp^*\colon H^i_c(Y,\overline{\Q}_\ell)\to \varinjlim_K H^i_c(X_K,\overline{\Q}_\ell)$
 by the composite of
 \[
  H^i_c(Y,\overline{\Q}_\ell)\xrightarrow{\pi^*}H^i_c(\mathcal{X}_s,\overline{\Q}_\ell)\xrightarrow{\spp^*}H^i_c(d(\mathcal{X}),\overline{\Q}_\ell)=H^i_c(X_K,\overline{\Q}_\ell)\to \varinjlim_K H^i_c(X_K,\overline{\Q}_\ell),
 \]
 where we also write $\pi$ for the morphism $\mathcal{X}_s\to Y$ induced by $\pi\colon R\to A_0/\mathfrak{m}A_0$.
\end{defn}

\begin{lem}\label{lem:sp-well-def}
 The map $\spp^*$ in Definition \ref{defn:specialization-map-infinite-level} is independent of
 the choice of $(A_0,\pi)\in\mathcal{C}$.
\end{lem}

\begin{prf}
 Since $\mathcal{C}$ is filtered, it suffices to compare $\spp^*$ for $(A_0,\pi)$ and that for
 $(A_0',\pi')$ under the condition $(A_0,\pi)\le (A'_0,\pi')$. Put $\mathcal{X}'=\Spf A_0'$.
 By Lemma \ref{lem:transition-finite}, we have a finite morphism $\mathcal{X}'\to \mathcal{X}$,
 which we denote by $f$.
 We have only to prove that the following diagram is commutative for $\Lambda=\Z/\ell^m\Z$:
 \[
  \xymatrix{%
 H^i_c(\mathcal{X}_s,\Lambda)\ar[r]^-{\spp^*}\ar[d]^-{f_s^*}& H^i_c(d(\mathcal{X}),\Lambda)\ar[d]^-{d(f)^*}\\
 H^i_c(\mathcal{X}'_s,\Lambda)\ar[r]^-{\spp^*}&H^i_c(d(\mathcal{X}'),\Lambda)\lefteqn{.}
 }
 \]
 Note that we have the following commutative diagram of sites
 (\cite[(3.5.4)]{MR1734903}):
  \[
  \xymatrix{%
  d(\mathcal{X}')_{\mathrm{\et}}\ar[r]^-{\lambda_{\mathcal{X}'}}\ar[d]^-{d(f)}& \mathcal{X}'_{s,\et}\ar[d]^-{f_s}\\
  d(\mathcal{X})_{\mathrm{\et}}\ar[r]^-{\lambda_{\mathcal{X}}}& \mathcal{X}_{s,\et}\lefteqn{.}
 }
  \]
 This gives rise to the following commutative diagram in $D^+(\mathcal{X}_{s,\et},\Lambda)$:
 \[
 \xymatrix{%
 \Lambda\ar[r]^-{\adj}\ar[dd]^-{\adj}& R\lambda_{\mathcal{X}*}\Lambda\ar[d]^-{\adj}\\
 & R\lambda_{\mathcal{X}*}Rd(f)_*\Lambda\ar@{=}[d]\\
 Rf_{s*}\Lambda\ar[r]^-{\adj}& Rf_{s*}R\lambda_{\mathcal{X}'*}\Lambda\lefteqn{.}
 }
 \]
 By taking $H^i_c$ of this diagram, we get the desired commutativity.
\end{prf}

Next we discuss the functoriality of $\spp^*$. Here we only consider automorphisms of $A$
which come from those on the tower $\{X_K\}$.

\begin{defn}\label{defn:finite-level-aut}
 Let $\sigma$ be an automorphism of $A$ which preserves $k^\circ$ (i.e., $\sigma(k^\circ)=k^\circ$).
 We say that $\sigma$ is of finite level if it satisfies $\sigma(A')=A'$.
\end{defn}

\begin{lem}\label{lem:aut-Hecke-tower}
 Let $\sigma$ be an automorphism of finite level of $A$.
 Then the following hold.
 \begin{enumerate}
  \item For each open normal subgroup $K$ of $K_0$, there exists an open normal subgroup $K'$ of $K_0$
	such that $\sigma(B_K)\subset B_{K'}$. 
  \item The automorphism $\sigma$ naturally induces an automorphism of
	$\varinjlim_K H^i_c(X_K,\overline{\Q}_\ell)$.
  \item For any $A_0\in\mathcal{D}$, there exists $A_1\in \mathcal{D}$ such that $\sigma(A_0)\subset A_1$.
	For such $A_1$, the homomorphism $\sigma\colon A_0\to A_1$ is finite and continuous.
 \end{enumerate}
\end{lem}

\begin{prf}
 First note that $\sigma(\varpi)$ belongs to $\mathfrak{m}$, hence there exists an integer $r\ge 1$ such that
 $\sigma(\varpi)^r\in \varpi k^\circ$. In particular $\sigma\colon k^\circ\to k^\circ$ is continuous.
 The same holds for $\sigma^{-1}$.

 We prove (i). For an open normal subgroup $K$ of $K_0$, take $A_0\in\mathcal{D}_K$
 and its topological generators $a_1,\ldots,a_n\in A_0$ over $k^\circ$.
 Since $\sigma$ preserves $A'$, we have $\sigma(a_1),\ldots,\sigma(a_n)\in A'$.
 Therefore, we can find an open normal subgroup $K'$ of $K_0$ such that
 $\sigma(a_1),\ldots,\sigma(a_n)$ are fixed by $K'$.
 We will prove that $\sigma(A_0)$ is fixed by $K'$.
 Take $g\in K'$ and $x\in A_0$. Then, for every integer $m\ge 1$,
 there exists $f\in k^\circ[T_1,\ldots,T_n]$ such that
 $x-f(a_1,\ldots,a_n)\in\varpi^{rm} A_0\subset \varpi^{rm} A$.
 Then we have $\sigma(x)-(\sigma f)(\sigma(a_1),\ldots,\sigma(a_n))\in \sigma(\varpi)^{rm}A\subset \varpi^mA$.
 Since the action of $g\in K'$ is an automorphism over $k^\circ$ and $\sigma(a_1),\ldots,\sigma(a_n)$
 are fixed by $K'$, we obtain
 $g(\sigma(x))-(\sigma f)(\sigma(a_1),\ldots,\sigma(a_n))\in\varpi^m A$,
 and thus $\sigma(x)-g(\sigma(x))\in \varpi^mA$ for every $m\ge 1$.
 As $A$ is $\varpi$-adically complete, this means that $\sigma(x)$ is fixed by $g\in K'$.
 Hence we have $\sigma(A_0)\subset (A')^{K'}=A_{K'}$. Now we conclude that
 \[
  \sigma(B_K)=\sigma(A_0[1/\varpi])=\sigma(A_0)[1/\varpi]\subset A_{K'}[1/\varpi]=B_{K'},
 \]
 hence (i).

 If we twist the $k$-algebra structure of $B_{K'}$ by $\sigma\colon k\to k$,
 $B_{K'}$ is again a complete f-adic $k$-algebra of topologically finite type, and
 the map $\sigma\colon B_K\to B_{K'}$ becomes a homomorphism over $k$.
 By \cite[Theorem 6.1.3/1]{MR746961}, $\sigma$ is automatically continuous.
 The composite of $B_K\xrightarrow{\sigma}B_{K'}\hookrightarrow B'$ is equal to
 that of $B_K\hookrightarrow B'\xrightarrow{\sigma}B'$, which is integral
 (see the proof of Lemma \ref{lem:tilde-integral} (ii)). 
 Hence $\sigma\colon B_K\to B_{K'}$ is also integral.
 By \cite[Theorem 6.3.5/1]{MR746961}, this means that $\sigma\colon B_K\to B_{K'}$ is finite.
 Therefore a finite morphism $X_{K'}=\Spa(B_{K'},B_{K'}^\circ)\to \Spa(B_K,B_K^\circ)=X_K$ and
 a map $H^i_c(X_K,\overline{\Q}_\ell)\to H^i_c(X_{K'},\overline{\Q}_\ell)$ are induced.
 Passing to the inductive limit, we obtain
 $\varinjlim_K H^i_c(X_K,\overline{\Q}_\ell)\to \varinjlim_K H^i_c(X_K,\overline{\Q}_\ell)$,
 which we denote by $\sigma_*$.
 Applying the same argument to $\sigma^{-1}$, we obtain 
 $(\sigma^{-1})_*\colon \varinjlim_K H^i_c(X_K,\overline{\Q}_\ell)\to \varinjlim_K H^i_c(X_K,\overline{\Q}_\ell)$.
 It is immediate to show that $\sigma_*$ and $(\sigma^{-1})_*$ are inverse to each other.
 Therefore $\sigma_*$ is an isomorphism. This completes the proof of (ii).

 Let us prove (iii). By the argument above, for $A_0\in\mathcal{D}_K$, we can find an open normal subgroup $K'$
 of $K_0$ such that $\sigma(A_0)\subset A_{K'}$, and the map $\sigma\colon A_0\to A_{K'}$ is continuous.
 Take a system of topological generators $a_1,\ldots,a_n$ of $A_0$ over $k^\circ$ and 
 an element $A_0'\in\mathcal{D}_{K'}$.
 Put $A_1=A_0'[\sigma(a_1),\ldots,\sigma(a_n)]$.
 Since $\sigma(a_1),\ldots,\sigma(a_n)\in A_{K'}$, we can observe that
 $A_1\in\mathcal{D}_{K'}$ in the same way as in the proof of
 Lemma \ref{lem:smaller-K}.
 On the other hand, by the continuity of $\sigma$, there exists an integer $m\ge 1$ such that
 $\sigma(\varpi^mA_0)\subset A_0'$. Therefore, by the same way as in the third paragraph of
 the proof of Proposition \ref{prop:open-limit},
 we can check that $\sigma(A_0)\subset A_1$. This concludes the existence of $A_1$.
 
 Next, we take an arbitrary $A_1\in\mathcal{D}$ such that $\sigma(A_0)\subset A_1$.
 The composite of $A_0\xrightarrow{\sigma}A_1\hookrightarrow A'$ is equal to
 that of $A_0\hookrightarrow A'\xrightarrow{\sigma}A'$, which is integral
 by Proposition \ref{prop:open-limit} (iii) and the proof of Lemma \ref{lem:tilde-integral} (ii). 
 Hence $\sigma\colon A_0\to A_1$ is also integral.
 This implies that $\sigma\colon A_0/\varpi A_0\to A_1/\sigma(\varpi)A_1$ is integral.
 On the other hand, as $A_1/\sigma(\varpi)A_1$ is finitely generated over $k^\circ$
 and $\sigma$ preserves $k^\circ$,
 the map $\sigma\colon A_0/\varpi A_0\to A_1/\sigma(\varpi)A_1$ is of finite type.
 Hence $\sigma\colon A_0/\varpi A_0\to A_1/\sigma(\varpi)A_1$ is finite, and so is $\sigma\colon A_0\to A_1$
 by \cite[Chapter 0, Proposition 7.2.4]{Fujiwara-Kato}. 
 By \cite[Proposition 1.5 (ii)]{MR1207303}, the topology of $A_0$ and $A_1$ are $\varpi$-adic.
 Since $\sigma(\varpi^{rm} A_0)\subset \varpi^mA_1$, $\sigma\colon A_0\to A_1$ is continuous.
 This completes the proof of (iii).
\end{prf}

\begin{prop}\label{prop:specialization-map-functoriality}
 Let $\sigma$ be an automorphism of finite level of $A$ and $\overline{\sigma}$ an automorphism of $R$.
 Assume that the fixed isomorphism $A/\mathfrak{m}A\cong \varinjlim_{x\mapsto x^p}R$ is compatible with
 $\sigma$ and $\overline{\sigma}$. Then, the specialization map
 $\spp^*\colon H^i_c(Y,\overline{\Q}_\ell)\to \varinjlim_K H^i_c(X_K,\overline{\Q}_\ell)$ is
 compatible with the actions of $\overline{\sigma}$ and $\sigma$.
\end{prop}

\begin{prf}
 We may assume that $R\neq 0$. Then $A/\mathfrak{m}A\neq 0$.
 Since $\kappa\subset A/\mathfrak{m}A$ is preserved by $\sigma$,
 we conclude that $\kappa\subset R$ is preserved by $\overline{\sigma}$ and 
 $\overline{\sigma}\vert_{\kappa}=\sigma\vert_{k^\circ}\bmod \mathfrak{m}$.

 Take $(A_0,\pi)\in \mathcal{C}$. By Lemma \ref{lem:aut-Hecke-tower} (iii), we can find $A_1\in\mathcal{D}$
 containing $\sigma(A_0)$. We write $\pi'$ for the composite
 $R\xrightarrow{\overline{\sigma}^{-1}} R\xrightarrow{\pi}A_0/\mathfrak{m}A_0\xrightarrow{\sigma}A_1/\mathfrak{m}A_1$, which is a $\kappa$-homomorphism.
 By the commutative diagram
 \[
  \xymatrix{%
 R\ar[r]^-{\pi}\ar[d]^-{\overline{\sigma}}& A_0/\mathfrak{m}A_0\ar[r]\ar[d]^-{\sigma}& A/\mathfrak{m}A\ar[d]^-{\sigma}\\
 R\ar[r]^-{\pi'}& A_1/\mathfrak{m}A_1\ar[r]& A/\mathfrak{m}A
 }
 \]
 and Lemma \ref{lem:aut-Hecke-tower} (iii), we can observe that $(A_1,\pi')\in \mathcal{C}$.
 Put $\mathcal{X}=\Spf A_0$ and $\mathcal{X}'=\Spf A_1$. 
 By Lemma \ref{lem:aut-Hecke-tower} (iii), the homomorphism $\sigma\colon A_0\to A_1$
 induces a finite morphism $\sigma\colon \mathcal{X}'\to \mathcal{X}$ of formal schemes.
 As in the proof of Lemma \ref{lem:sp-well-def}, the commutative diagram of sites
  \[
  \xymatrix{%
  d(\mathcal{X}')_{\mathrm{\et}}\ar[r]^-{\lambda_{\mathcal{X}'}}\ar[d]^-{\sigma}& \mathcal{X}'_{s,\et}\ar[d]^-{\sigma}\ar[r]^-{\pi'}& Y_{\et}\ar[d]^-{\overline{\sigma}}\\
  d(\mathcal{X})_{\mathrm{\et}}\ar[r]^-{\lambda_{\mathcal{X}}}& \mathcal{X}_{s,\et}\ar[r]^-{\pi}& Y_\et
 }
  \]
 gives rise to the following commutative diagram for $\Lambda=\Z/\ell^m\Z$:
 \[
  \xymatrix{%
 H^i_c(Y,\Lambda)\ar[r]^-{\pi}\ar[d]^-{\overline{\sigma}}& H^i_c(\mathcal{X}_s,\Lambda)\ar[r]^-{\spp^*}\ar[d]^-{\sigma}& H^i_c(d(\mathcal{X}),\Lambda)\ar[d]^-{\sigma}\ar@{=}[r]&H^i_c(X_K,\Lambda)\ar[d]^-{\sigma}\\
 H^i_c(Y,\Lambda)\ar[r]^-{\pi'}& H^i_c(\mathcal{X}'_s,\Lambda)\ar[r]^-{\spp^*}&H^i_c(d(\mathcal{X}'),\Lambda)\ar@{=}[r]&H^i_c(X_{K'},\Lambda)\lefteqn{.}
 }
 \]
 Here $K$ (resp.\ $K'$) is the open normal subgroup of $K_0$ such that $A_0\in\mathcal{D}_K$
 (resp.\ $A_1\in\mathcal{D}_{K'}$). By the construction of the specialization map,
 this implies the commutativity of the diagram
 \[
 \xymatrix{%
 H^i_c(Y,\overline{\Q}_\ell)\ar[r]^-{\spp^*}\ar[d]^-{\overline{\sigma}}& \varinjlim_K H^i_c(X_K,\overline{\Q}_\ell)\ar[d]^-{\sigma}\\
  H^i_c(Y,\overline{\Q}_\ell)\ar[r]^-{\spp^*}& \varinjlim_K H^i_c(X_K,\overline{\Q}_\ell)\lefteqn{,}
 }
 \]
 which concludes the proof.
\end{prf}

Here is the infinite level version of Theorem \ref{thm:injectivity-criterion-finite-level}.

\begin{thm}\label{thm:injectivity-criterion-infinite-level}
 Let $\{Z_K\}_{K\lhd K_0}$ be a projective system of adic spaces locally of finite type, separated
 and taut over $\Spa(k,k^\circ)$ with proper transition maps.
 Suppose that we are given a system of open immersions $\{X_K\hookrightarrow Z_K\}_{K\lhd K_0}$
 over $\Spa(k,k^\circ)$ which makes the cartesian diagram
 \[
  \xymatrix{%
 X_{K'}\ar[r]\ar[d]& Z_{K'}\ar[d]\\ X_K\ar[r]& Z_K
 }
 \]
 for every open normal subgroups $K$, $K'$ of $K_0$ with $K'\subset K$.

 Assume that $Y$ is pure-dimensional and smooth over $\kappa$.
 Let $V$ be a subspace of $H^i_c(Y,\overline{\Q}_\ell)$ such that the composite
 $V\hookrightarrow H^i_c(Y,\overline{\Q}_\ell)\to H^i(Y,\overline{\Q}_\ell)$ is
 an injection.
 Then, the composite
 \[
  V\hookrightarrow H^i_c(Y,\overline{\Q}_\ell)\xrightarrow{\spp^*} \varinjlim_K H^i_c(X_K,\overline{\Q}_\ell)\to \varinjlim_K H^i_c(Z_K,\overline{\Q}_\ell)
 \]
 is also an injection.
\end{thm}

\begin{prf}
 Fix $(A_0,\pi)\in\mathcal{C}$ with $A_0\in\mathcal{D}_K$.
 Let $K'$ be an open normal subgroup of $K_0$ contained in $K$.
 Then, by Lemmas \ref{lem:smaller-K} and \ref{lem:transition-finite},
 we can find $(A_0',\pi')\in\mathcal{C}$ such that $(A_0,\pi)\le (A_0',\pi')$ and $A_0'\in\mathcal{D}_{K'}$.
 Put $\mathcal{X}'=\Spf A_0'$.
 It is of topologically finite type and flat over $\Spf k^\circ$, 
 and pseudo-compactifiable (see \cite[Example 4.25 (i)]{formalnearby}).
 By Proposition \ref{prop:open-limit} (iii), we have $d(\mathcal{X}')=X_{K'}$,
 which is smooth over $\Spa(k,k^\circ)$.
 By Lemma \ref{lem:property-C} (ii), the morphism 
 $\pi'\colon \mathcal{X}'_s\to \Spec R=Y$ is finite and surjective,
 and the dimension of $\mathcal{X}'_s$ is equal to $\dim Y$. 
 Therefore, Theorem \ref{thm:injectivity-criterion-finite-level} tells us that the composite
 \[
 V\hookrightarrow H^i_c(Y,\overline{\Q}_\ell)\xrightarrow{\spp^*\circ\pi'} H^i_c(d(\mathcal{X}'),\overline{\Q}_\ell)
 =H^i_c(X_{K'},\overline{\Q}_\ell)\to H^i_c(Z_{K'},\overline{\Q}_\ell)
 \]
 is injective. Obviously, the inductive limit of this map with respect to $K'$ coincides with
 \[
 V\hookrightarrow H^i_c(Y,\overline{\Q}_\ell)\xrightarrow{\spp^*} \varinjlim_{K'} H^i_c(X_{K'},\overline{\Q}_\ell)\to \varinjlim_{K'} H^i_c(Z_{K'},\overline{\Q}_\ell)
 \]
 considered in the theorem. Thus we obtain the desired injectivity.
\end{prf}

We will end this section by giving a simple but helpful example.

\begin{exa}\label{exa:non-reduced}
 Here we assume that $p\neq 2$. Take $\varpi'\in \mathfrak{m}$ such that $\varpi'^2=\varpi$,
 and fix a sequence $(\varpi'^{1/p^m})_{m\ge 0}$ of $p^m$th roots of $\varpi'$.
 Put $\varpi^{1/p^m}=(\varpi'^{1/p^m})^2$.
 Let $I$ be the ideal of $k^\circ\langle T^{p^{-\infty}}\rangle$ (see Example \ref{exa:perfection})
 generated by $T^{2/p^m}-\varpi^{1/p^m}$ for every $m\ge 0$, and $\overline{I}$ the closure of $I$
 in $k^\circ\langle T^{p^{-\infty}}\rangle$.
 Let us consider $A=k^\circ\langle T^{p^{-\infty}}\rangle/\overline{I}$.
 \begin{enumerate}
  \item The $k^\circ$-algebra $A$ is $\varpi$-adically complete and flat.
  \item Define continuous $k$-homomorphisms 
	\[
	 \phi\colon k[S]/(S^2-1)\to A[1/\varpi],\quad \psi\colon A[1/\varpi]\to k\times k
	\]
	by $\phi(S)=\varpi'^{-1}T$ and
	$\psi(T^{1/p^m})=(\varpi'^{1/p^m},-\varpi'^{1/p^m})$.
	Then, $\phi$ and $\psi$ are isomorphisms.
	In particular, $A[1/\varpi]$ is a complete f-adic $k$-algebra of topologically finite type.
  \item We have $\psi(A)=\{(x,y)\in k^\circ\times k^\circ\mid x-y\in\mathfrak{m}\}$. In particular,
	$A$ is a ring of definition of $A[1/\varpi]$ which is not topologically finitely generated over $k^\circ$.
  \item We have $A/\mathfrak{m}A\cong \kappa$.
  \item Let $\widetilde{A}$ be the integral closure of $A$ in $A[1/\varpi]$.
	We have $\psi(\widetilde{A})=k^\circ\times k^\circ$. Namely, $\widetilde{A}$ coincides with 
	$A[1/\varpi]^\circ$.
  \item Put $A_m=k^\circ\langle T_m\rangle/(T_m^2-\varpi^{1/p^m})=k^\circ[T_m]/(T_m^2-\varpi^{1/p^m})$.
	We have a $k^\circ$-homomorphism $A_m\to A$; $T_m\mapsto T^{1/p^m}$, which turns out to be injective.
	We write $\psi_m$ for the composite $A_m\to A\xrightarrow{\psi}k\times k$.
	Then, we have $\psi_m(A_m)=\{(x,y)\in k^\circ\times k^\circ\mid x-y\in\varpi'^{1/p^m}k^\circ\}$.
	In particular, $A\cong \varinjlim_m A_m$.
 \end{enumerate}
 By (i), (ii), (iv), (v), this $A$ with the obvious action of $K_0=1$ satisfy
 Assumptions \ref{assump:Hecke}, \ref{assump:fin-type-perfection}.
 By (iv), $Y=\Spec \kappa$ is a point, while $\Spa(A[1/\varpi],A[1/\varpi]^\circ)$ consists of two points
 by (ii). In particular, the specialization map $\spp^*$ on $H^0_c$ in
 Definition \ref{defn:specialization-map-infinite-level} is not surjective.
\end{exa}

\begin{prf}
 The assertion (i) follows from Corollary \ref{cor:flatness-complete}, as the image of $T^2-\varpi$
 in $\kappa[T^{p^{-\infty}}]$ is non-zero. 

 We prove (ii). It is easy to see that $\phi$ and $\psi$ are
 well-defined and $\psi\circ\phi$ is an isomorphism. Therefore $\phi$ is an injection.
 Put $C=k^\circ[S]/(S^2-1)$. Since $\phi(\varpi'C)\subset A$, $\phi$ is continuous.
 We shall prove that $\phi$ is surjective. 
 First we prove that $A\subset \phi(C)+\varpi A$. It suffices to show that $T^{1/p^m}\in \phi(C)$.
 Note that $(T^{1/p^m}/\varpi'^{1/p^m})^2=T^{2/p^m}/\varpi^{1/p^m}=1$ in $A[1/\varpi]$.
 Since we are assuming $p\neq 2$, 
 we have 
 \[
  \phi(S)=\frac{T}{\varpi'}=\Bigl(\frac{T^{1/p^m}}{\varpi'^{1/p^m}}\Bigr)^{p^m}=\frac{T^{1/p^m}}{\varpi'^{1/p^m}}.
 \]
 Hence $T^{1/p^m}=\phi(\varpi'^{1/p^m}S)\in \phi(C)$, as desired. Now we prove the surjectivity of $\phi$.
 Since $k[S]/(S^2-1)$ is $\varpi$-divisible, it suffices to see that every $a_0\in A$ lies in $\phi(C)$.
 By $A\subset \phi(C)+\varpi A$, we can find $x_0\in C$ and $a_1\in A$ such that $a_0=\phi(x_0)+\varpi a_1$.
 Continuing this process, we get $(x_m)_{m\ge 0}$ and $(a_m)_{m\ge 1}$ such that $x_m\in C$, $a_m\in A$
 and $a_m=\phi(x_m)+\varpi a_{m+1}$. Put $x=\sum_{m=0}^\infty \varpi^mx_m\in C$.
 Since $\phi$ is continuous and $A$ is $\varpi$-adically separated, we have
 \[
  \phi(x)=\lim_{n\to\infty}\phi\Bigl(\sum_{m=0}^n\varpi^mx_m\Bigr)=\lim_{n\to\infty}(a_0-\varpi^{n+1}a_{n+1})=a_0.
 \]
 Hence $a_0$ lies in $\phi(C)$. Therefore $\phi$ is a bijection. 
 As $A\subset \phi(C)$, $\phi^{-1}$ is continuous. Thus $\phi$ is an isomorphism.
 Since $\psi\circ\phi$ is an isomorphism, so is $\psi$.

 Let us prove (iii). We can easily check that 
 $\psi(A)\subset \{(x,y)\in k^\circ\times k^\circ\mid x-y\in\mathfrak{m}\}$.
 We prove the reverse inclusion. Let $x,y\in k^\circ$ be elements such that $x-y\in\mathfrak{m}$.
 Since $(x,y)=(\psi\circ\phi)(\frac{x+y}{2}+\frac{x-y}{2}S)$, we have 
 \[
  \psi^{-1}(x,y)=\phi\Bigl(\frac{x+y}{2}+\frac{x-y}{2}S\Bigr)=\frac{x+y}{2}+\frac{x-y}{2}\cdot \frac{T}{\varpi'}
 =\frac{x+y}{2}+\frac{x-y}{2}\cdot \frac{T^{1/p^m}}{\varpi'^{1/p^m}}.
 \]
 Since $x-y\in \mathfrak{m}$, for a sufficiently large $m$, $(x-y)/\varpi'^{1/p^m}$ belongs to $k^\circ$.
 Hence we have $\psi^{-1}(x,y)\in A$, that is, $(x,y)\in \psi(A)$.
 Suppose that $\psi(A)$ is topologically finitely generated over $k^\circ$, and take topological generators
 $(x_1,y_1),\ldots,(x_n,y_n)\in \psi(A)$. Then, there exists an integer $m\ge 0$ such that
 $x_i-y_i\in \varpi^{1/p^m}k^\circ$ for every $1\le i\le n$.
 Since $\psi(A)/\varpi^{1/p^m}(k^\circ\times k^\circ)$ is generated by the images of $(x_1,y_1),\ldots,(x_n,y_n)$,
 it is included in $\{(x,x)\mid x\in k^\circ/\varpi^{1/p^m}k^\circ\}$. It is absurd, as
 $\psi(A)/\varpi^{1/p^m}(k^\circ\times k^\circ)$ contains an element $(\varpi^{1/p^{m+1}},0)$.
 Thus $\psi(A)$ is not topologically finitely generated over $k^\circ$.

 For (iv), it suffices to prove that the natural map $\kappa=k^\circ/\mathfrak{m}\to A/\mathfrak{m}A$ is
 surjective, or equivalently, $A=k^\circ+\mathfrak{m}A$.
 By (ii), we may replace $A$ by $\psi(A)$. Take $(x,y)\in \psi(A)$. By (iii), we have $x,y\in k^\circ$ and
 $x-y\in\mathfrak{m}$.
 Therefore we have $(x,y)=(y,y)+(x-y,0)\in k^\circ+\mathfrak{m}\psi(A)$.
 Now we conclude that $\psi(A)=k^\circ+\mathfrak{m}\psi(A)$.

 We prove (v). Since $k^\circ\times k^\circ$ is integrally closed in $k\times k$, it suffices to show that
 $k^\circ\times k^\circ$ is integral over $\psi(A)$. Take any element $x=(a,b)\in k^\circ\times k^\circ$.
 Then $(a,a)$ and $(b,b)$ belong to $\psi(A)$, and $x$ is a root of the polynomial
 $(T-(a,a))(T-(b,b))$ in $\psi(A)[T]$. This means that $x$ is integral over $\psi(A)$.

 Finally consider (vi). It is immediate to see that the homomorphism $A_m[1/\varpi]\to k\times k$
 induced from $\psi_m$ is an isomorphism. Therefore $A_m[1/\varpi]\to A[1/\varpi]$ is also an isomorphism.
 Since $A_m$ is $\varpi$-torsion free, we conclude that $A_m\to A$ is an injection.
 Let us prove that $\psi_m(A_m)=\{(x,y)\in k^\circ\times k^\circ\mid x-y\in \varpi'^{1/p^m}k^\circ\}$.
 As $\psi_m(T_m)=(\varpi'^{1/p^m},-\varpi'^{1/p^m})$, $\psi_m(A_m)$ is contained in the right hand side.
 Since the $k^\circ$-algebra 
 $\{(x,y)\in k^\circ\times k^\circ\mid x-y\in \varpi'^{1/p^m}k^\circ\}$ is generated by
 $(\varpi'^{1/p^m},-\varpi'^{1/p^m})$, the desired equality holds.
\end{prf}

\section{The Lubin-Tate tower}\label{sec:LT-tower}
Let $F$ be a non-archimedean local field. We denote by $\mathcal{O}_F$ the ring of integers of $F$
and $\mathfrak{p}_F$ the maximal ideal of $\mathcal{O}_F$.
Write $q$ for the characteristic of the residue field $\mathcal{O}_F/\mathfrak{p}_F$, and $p$ for
the characteristic of $\F_q=\mathcal{O}_F/\mathfrak{p}_F$. 
The normalized valuation of $F$ is denoted by $v_F$.
Fix a uniformizer $\varpi\in \mathcal{O}_F$.

We fix an algebraic closure $\overline{F}$ of $F$. Let $F^{\mathrm{ur}}$ (resp.\ $F^{\mathrm{ab}}$)
be the maximal unramified (resp.\ abelian) extension of $F$ inside $\overline{F}$.
Denote by $\breve{F}$, $\widehat{F}^{\mathrm{ab}}$ and $C$ for the completion of 
$F^{\mathrm{ur}}$, $F^{\mathrm{ab}}$ and $\overline{F}$, respectively.
The residue field of $\breve{F}$ is denoted by $\overline{\F}_q$, which is
an algebraic closure of $\F_q$.
We write $W_F$ for the Weil group of $F$,
and $\Art_F\colon F^\times\xrightarrow{\cong}W_F^{\mathrm{ab}}$ for the isomorphism of the local
class field theory, which is normalized so that $\Art_F(\varpi)$ is a lift of
the geometric Frobenius automorphism on $\overline{\F}_q$.

We denote by $\mathbf{Nilp}$ the category of schemes over $\mathcal{O}_{\breve{F}}$ on which
$\varpi$ is locally nilpotent, and by $\mathbf{Set}$ the category of sets.
For $S\in\mathbf{Nilp}$, we set $\overline{S}=S\otimes_{\mathcal{O}_{\breve{F}}}\overline{\F}_q$.

Let $n\ge 1$ be an integer.
We fix a one-dimensional formal $\mathcal{O}_F$-module $\mathbb{X}$ of height $n$ over $\overline{\F}_q$,
which is unique up to isomorphism. Put $D=\End_{\mathcal{O}_F}(\mathbb{X})\otimes_{\mathcal{O}_F}F$.
It is known to be a central division algebra over $F$ with invariant $1/n$.

For an integer $m\ge 0$, 
let $\mathcal{M}_m\colon \mathbf{Nilp}\to \mathbf{Set}$ be the functor
that sends $S$ to the set of isomorphism classes of triples $(X,\rho,\eta)$,
where 
\begin{itemize}
 \item $X$ is a formal $\mathcal{O}_F$-module on $S$,
 \item $\rho\colon \mathbb{X}\otimes_{\overline{\F}_q}\overline{S}\to X\times_S\overline{S}$ is
       an $\mathcal{O}_F$-quasi-isogeny, and
 \item $\eta$ is a Drinfeld $m$-level structure on $X$.
\end{itemize}
For an integer $\delta$, let $\mathcal{M}_m^{(\delta)}$ be the subfunctor of $\mathcal{M}_m$
corresponding to the triple $(X,\rho,\eta)$ with $\height_{\mathcal{O}_F}\rho=\delta$.
Then, $\mathcal{M}_m^{(0)}$ can be identified with the deformation functor
of $\mathbb{X}$ with Drinfeld $m$-level structures.
In \cite[\S 4]{MR0384707}, it is proved that $\mathcal{M}_m^{(0)}$ is represented by $\Spf A_m$, 
where $A_m$ is an $n$-dimensional Noetherian regular complete local
$\mathcal{O}_{\breve{F}}$-algebra.
Since $\mathcal{M}_m^{(\delta)}$ is (non-canonically) isomorphic to $\mathcal{M}_m^{(0)}$,
we conclude that $\mathcal{M}_m=\coprod_{\delta\in\Z}\mathcal{M}_m^{(\delta)}$ is represented
by a locally Noetherian formal scheme over $\mathcal{O}_{\breve{F}}$.
Moreover, $\{\mathcal{M}_m\}_{m\ge 0}$ form a projective system of formal schemes with finite flat transition 
morphisms. We call it the Lubin-Tate tower.

The Lubin-Tate tower is equipped with several actions.
First, the group $D^\times=\mathbf{QIsog}_{\mathcal{O}_F}(\mathbb{X})$ of self-quasi-isogenies of $\mathbb{X}$
acts on each $\mathcal{M}_m$ on the right. The formal scheme $\mathcal{M}_m$ also has a natural Weil descent datum
in the sense of \cite[Definition 3.45]{MR1393439}.
These are compatible with the transition maps of $\{\mathcal{M}_m\}_{m\ge 0}$.
Further, the group $\GL_n(F)$ acts on the tower $\{\mathcal{M}_m\}_{m\ge 0}$ on the right as a pro-object.
This action may change the level $m$. However, the subgroup $\GL_n(\mathcal{O}_F)$ preserves each $\mathcal{M}_m$,
and the $m$th principal congruence subgroup $K_m=\Ker(\GL_n(\mathcal{O}_F)\to \GL_n(\mathcal{O}_F/\mathfrak{p}_F^m))$
acts trivially on $\mathcal{M}_m$.

We denote by $M_m$ the rigid generic fiber of $\mathcal{M}_m$, which is known to be an $n-1$-dimensional
smooth rigid space over $\breve{F}$. From $\{\mathcal{M}_m\}_{m\ge 0}$ we obtain
a projective system $\{M_m\}_{m\ge 0}$ of rigid spaces over $\breve{F}$, whose transition
maps are finite, \'etale and Galois. For $m'\ge m\ge 0$, the action of $K_m$ on $M_{m'}$ gives
identification between $K_m/K_{m'}$ and the Galois group of $M_{m'}\to M_m$.
For each compact open subgroup $K$ of $\GL_n(\mathcal{O}_F)$, 
we can define the Lubin-Tate space $M_K$ of level $K$ as follows:
take an integer $m\ge 0$ such that $K_m\subset K$ and put $M_K=M_m/(K/K_m)$,
which is in fact independent of the choice of $m$.
Now we obtain a projective system $\{M_K\}_{K\subset\GL_n(\mathcal{O}_F)}$,
whose transition maps are finite and \'etale.
It is also equipped with actions of $D^\times$ and $\GL_n(F)$, and a Weil descent datum.

For our purpose, it is convenient to consider the quotient of $M_K$ by the discrete subgroup
$\varpi^\Z\subset F^\times\subset D^\times$, which is denoted by $M_{K,\varpi^\Z}$.
Note that the decomposition $\mathcal{M}_m=\coprod_{\delta\in \Z}\mathcal{M}_m^{(\delta)}$ induces 
$M_K=\coprod_{\delta\in \Z}M_K^{(\delta)}$, and $\varpi\in D^\times$ maps $M_K^{(\delta)}$ to $M_K^{(\delta+n)}$.
Therefore, the quotient $M_{K,\varpi^\Z}$ is isomorphic to $\coprod_{0\le \delta<n}M_K^{(\delta)}$.

Take a prime number $\ell\neq p$. Fix an isomorphism $\overline{\Q}_\ell\cong \C$ and identify them.
We consider the compactly supported $\ell$-adic cohomology
of the tower $\{M_{K,\varpi^\Z}\}_{K\subset\GL_n(\mathcal{O}_F)}$.

\begin{defn}
 We put $H_{\mathrm{LT}}=\varinjlim_K H^{n-1}_c(M_{K,\varpi^\Z}\otimes_{\breve{F}}C,\overline{\Q}_\ell)$.
\end{defn}

 The group $\GL_n(F)\times D^\times$ naturally acts on $H_{\mathrm{LT}}$. This action is known to be smooth
 (see \cite[Lemma 2.5.1]{MR2383890}).
 By \cite[Lemma 5.36]{MR1393439},
 the subgroup $\varpi^\Z\subset F^\times\subset \GL_n(F)$ acts trivially on $H_{\mathrm{LT}}$.
 On the other hand, by using the Weil descent datum on $M_{K,\varpi^\Z}$, one can define an action of $W_F$
 on $H_{\mathrm{LT}}$. As a result, $H_{\mathrm{LT}}$ becomes a representation of 
 $(\GL_n(F)/\varpi^\Z)\times (D^\times/\varpi^\Z)\times W_F$.
 The $\GL_n(F)$-supercuspidal part of this representation can be described
 by using the local Langlands correspondence and the local Jacquet-Langlands correspondence:

\begin{thm}[Non-abelian Lubin-Tate theory, \cite{MR1876802}, \cite{MR1719811}]\label{thm:NALT}
 Let $\pi$ be an irreducible supercuspidal representation of $\GL_n(F)$ whose central character
 is trivial on $\varpi^\Z$. We write $\rec_F(\pi)$ (resp.\ $\JL(\pi)$) for
 the irreducible smooth representation of $W_F$ (resp.\ $D^\times$) 
 corresponding to $\pi$ under the local Langlands correspondence
 (resp.\ the local Jacquet-Langlands correspondence).
 Then we have
 \[
  H_{\mathrm{LT},\pi^\vee}=\pi^\vee\boxtimes \JL(\pi)\boxtimes\rec_F(\pi)(\tfrac{1-n}{2})
 \]
 as representations of $\GL_n(F)\times D^\times\times W_F$, where $H_{\mathrm{LT},\pi^\vee}$ denotes
 the $\pi^\vee$-isotypic part of $H_{\mathrm{LT}}$.
\end{thm}

To simplify our argument, we also use the following result, which is in some sense stronger than
the theorem above:

\begin{thm}\label{thm:W_F-smooth}
 The action of $W_F$ on $H_{\mathrm{LT}}$ is smooth.
\end{thm}

\begin{prf}
 Since $D^\times/\varpi^\Z$ is compact, we have the isotypic decomposition
 \[
  H_{\mathrm{LT}}=\bigoplus_{\rho}H_{\mathrm{LT},\rho}, 
 \]
 where $\rho$ runs through irreducible smooth representations of $D^\times$ whose central characters are
 trivial on $\varpi^\Z$, and $H_{\mathrm{LT},\rho}$ denotes the $\rho$-isotypic part of $H_{\mathrm{LT}}$.
 We fix $\rho$ and denote by $\pi$ the discrete series representation of $\GL_n(F)$ satisfying $\JL(\pi)=\rho$.
 Such $\pi$ can be written as $\mathrm{St}_s(\pi')$, where $s$ is an integer dividing $n$ and 
 $\pi'$ is an irreducible supercuspidal representation of $\GL_{n/s}(F)$.
 By \cite{MR2511742} and \cite{MR2569585}, we have
 \[
  H_{\mathrm{LT},\rho}=\pi^\vee\boxtimes\rho\boxtimes \rec_F(\pi')(\tfrac{m}{2})
 \]
 for some $m\in\Z$. Since $\rec_F(\pi')$ is a smooth representation of $W_F$,
 we conclude that the action of $W_F$ on $H_{\mathrm{LT},\rho}$ is smooth.
\end{prf}

Let us recall the determinant map. 
We write $\{\mathcal{N}_m\}_m$ for the Lubin-Tate tower with respect to $\wedge^n\mathbb{X}$,
where $\wedge^n$ is the exterior product defined in \cite{MR3168925}.
Then, as in \cite[\S 2.5]{Weinstein-stabred}, we can construct the determinant map 
$\mathcal{M}^{(0)}_m\to \mathcal{N}^{(0)}_m$. 
On the other hand, by the classical Lubin-Tate theory, we have $\mathcal{N}^{(0)}_m=\Spf\mathcal{O}_{\breve{F}_m}$
where $\breve{F}_m$ denotes the completion of the $m$th Lubin-Tate extension of $F^{\mathrm{ur}}$.
Choose an $\breve{F}$-embedding $\varinjlim_{m}\breve{F}_m\hookrightarrow \widehat{F}^{\mathrm{ab}}$
and consider the base change 
$\mathcal{M}_{m,\mathcal{O}_C}^{(0)}=\mathcal{M}_m^{(0)}\widehat{\otimes}_{\mathcal{O}_{\breve{F}_m}}\mathcal{O}_C$
of $\mathcal{M}_m^{(0)}$
under $\mathcal{O}_{\breve{F}_m}\hookrightarrow \mathcal{O}_{\widehat{F}^{\mathrm{ab}}}\hookrightarrow\mathcal{O}_C$.
We put 
\[
 G^1=\{(g,d,\sigma)\in \GL_n(F)\times D^\times\times W_F\mid \det(g)^{-1}\Nrd(d)\Art_F^{-1}(\sigma)=1\},
\]
where $\Nrd\colon D^\times\to F^\times$ denotes the reduced norm.
By using the action of $\GL_n(F)\times D^\times$ and the Weil descent datum on $\{\mathcal{M}_m\}_{m}$,
we can define an action of $G^1$ on the tower $\{\mathcal{M}_{m,\mathcal{O}_C}^{(0)}\}_{m}$.
The action of $(g,d,\sigma)\in G^1$ makes the following diagram commute:
\[
 \xymatrix{%
 \{\mathcal{M}^{(0)}_{m,\mathcal{O}_C}\}_m\ar[r]^-{(g,d,\sigma)}\ar[d]&\{\mathcal{M}^{(0)}_{m,\mathcal{O}_C}\}_m\ar[d]\\
 \Spf\mathcal{O}_C\ar[r]^-{\sigma^*}& \Spf\mathcal{O}_C\lefteqn{.}
 }
\]
Passing to the rigid generic fiber, we obtain the rigid space $M^{(0)}_{m,C}=M^{(0)}_m\otimes_{\breve{F}_m}C$
over $C$.
For an integer $m\ge 0$, put $K'_m=\SL_n(\mathcal{O}_F)\cap K_m$.
It is easy to observe that the transition map $M^{(0)}_{m',C}\to M^{(0)}_{m,C}$ for $m'\ge m\ge 0$
is finite \'etale Galois with Galois group $K'_m/K'_{m'}$.
Therefore, by the same method as in the case of $\{M_K\}_{K\subset\GL_n(\mathcal{O}_F)}$,
we can extend $\{M^{(0)}_{m,C}\}_m$ to the projective system $\{M^{(0)}_{K,C}\}_{K\subset \SL_n(\mathcal{O}_F)}$, where
$K$ runs through open subgroups of $\SL_n(\mathcal{O}_F)$.
The group $G^1$ acts on this tower.

We put $H'_{\mathrm{LT}}=\varinjlim_K H^{n-1}_c(M^{(0)}_{K,C},\overline{\Q}_\ell)$.
The following proposition gives a connection between $H_{\mathrm{LT}}$ and $H'_{\mathrm{LT}}$.

\begin{prop}\label{prop:NALT-cInd}
 Put $G=\GL_n(F)\times (D^\times/\varpi^\Z)\times W_F$.
 The group $G^1$ can be regarded as a cocompact closed normal subgroup of $G$.
 Moreover, the action of $G^1$ on $H'_{\mathrm{LT}}$ is smooth,
 and we have a $G$-equivariant isomorphism 
 \[
  \Ind_{G^1}^G H'_{\mathrm{LT}}\cong H_{\mathrm{LT}}.
 \]
\end{prop}

\begin{prf}
 First we will prove that $G^1$ is a cocompact closed normal subgroup of $G$.
 The natural homomorphism $G^1\to G$ is clearly injective.
 Further, its image is equal to the kernel of the homomorphism $\nu\colon G\to F^\times/\varpi^{n\Z}$;
 $(g,\overline{d},\sigma)\mapsto \overline{\det(g)^{-1}\Nrd(d)\Art_F^{-1}(\sigma)}$.
 Hence the image of $G^1$ is closed normal and cocompact in $G$.

 Next we construct a $G^1$-equivariant surjection $H_{\mathrm{LT}}\to H'_{\mathrm{LT}}$.
 Take an integer $m\ge 0$. We have $M_m^{(0)}\otimes_{\breve{F}}C=\coprod_{\sigma\in\Gal(\breve{F}_m/\breve{F})}M_m^{(0)}\otimes_{\breve{F}_m,\sigma}C$.
 Therefore, the morphism $M^{(0)}_{m,C}=M_m^{(0)}\otimes_{\breve{F}_m}C\hookrightarrow M_m^{(0)}\otimes_{\breve{F}}C\hookrightarrow M_{m,\varpi^\Z}\otimes_{\breve{F}}C$ induces a surjective map
 \[
  H^{n-1}_c(M_{m,\varpi^\Z}\otimes_{\breve{F}}C,\overline{\Q}_\ell)\to H^{n-1}_c(M_m^{(0)}\otimes_{\breve{F}_m}C,\overline{\Q}_\ell).
 \]
 By taking inductive limit, we obtain a surjective linear map $\phi\colon H_{\mathrm{LT}}\to H'_{\mathrm{LT}}$.
 It is easy to see that this is $G^1$-equivariant.
 By Theorem \ref{thm:W_F-smooth}, the action of $G^1$ on $H_{\mathrm{LT}}$ is smooth.
 Therefore $H'_{\mathrm{LT}}$ is also a smooth representation of $G^1$.

 By the Frobenius reciprocity, we obtain a desired $G$-equivariant homomorphism 
 $\widetilde{\phi}\colon H_{\mathrm{LT}}\to\Ind_{G^1}^G H'_{\mathrm{LT}}$.
 We shall prove that this map is an isomorphism. 
 It suffices to show that $\widetilde{\phi}\colon H_{\mathrm{LT}}^{K_m}\to (\Ind_{G^1}^G H'_{\mathrm{LT}})^{K_m}$
 is an isomorphism for each $m\ge 1$. 
 We have 
 \[
   H_{\mathrm{LT}}^{K_m}=H^{n-1}_c(M_{m,\varpi^\Z}\otimes_{\breve{F}}C,\overline{\Q}_\ell)=\bigoplus_{0\le \delta<n}H^{n-1}_c(M_m^{(\delta)}\otimes_{\breve{F}}C,\overline{\Q}_\ell).
 \]
 On the other hand, note that $\nu\colon G\to F^\times/\varpi^{n\Z}$ induces a bijection
 $G^1\backslash G/K_m\xrightarrow{\cong}F^\times/\varpi^{n\Z}(1+\mathfrak{p}_F^m)$.
 In particular, it is a finite set with cardinality $r=n(\mathcal{O}_F^\times:1+\mathfrak{p}_F^m)$.
 Since $\Nrd\colon D^\times/\varpi^\Z\to F^\times/\varpi^{n\Z}$ is surjective, we can take
 a system of representatives $\{(1,d_i,1)\}_{1\le i\le r}$ of $G^1\backslash G/K_m$
 such that $d_1=1$.
 We have an injection
 \[
  (\Ind_{G^1}^G H'_{\mathrm{LT}})^{K_m}\hookrightarrow \bigoplus_{1\le i\le r}(H'_{\mathrm{LT}})^{K'_m}=\bigoplus_{1\le i\le r}H^{n-1}_c(M^{(0)}_{m,C},\overline{\Q}_\ell);\quad f\mapsto (f(1,d_i,1))_{1\le i\le r}.
 \]
 Take an element $x\in H^{n-1}_c(M^{(0)}_{m,C},\overline{\Q}_\ell)\subset H^{n-1}_c(M^{(0)}_m\otimes_{\breve{F}}C,\overline{\Q}_\ell)\subset H_{\mathrm{LT}}^{K_m}$ and consider $\widetilde{\phi}(x)\in (\Ind_{G^1}^G H'_{\mathrm{LT}})^{K_m}$.
 Let $d$ be an element of $D^\times$ and write $\overline{d}$ for the image of $d$ in $D^\times/\varpi^\Z$.
 If $\Nrd(d)\notin \varpi^{n\Z}\mathcal{O}_F^\times$, $(1,\overline{d},1)$ maps $M_{m,C}^{(0)}$ into
 $M_m^{(\delta)}\otimes_{\breve{F}}C$ with $\delta\neq 0$.
 If $\Nrd(d)\in \varpi^{n\Z}\mathcal{O}_F^\times\setminus \varpi^{n\Z}(1+\mathfrak{p}_F^m)$,
 $(1,\overline{d},1)$ maps $M_{m,C}^{(0)}$ to $M_m^{(0)}\otimes_{\breve{F}_m,\sigma}C$ with $\sigma\neq 1$.
 Therefore, unless $d\in \varpi^{n\Z}(1+\mathfrak{p}_F^m)$, we have $\phi((1,\overline{d},1)x)=0$.
 Hence we conclude that
 \[
  (\widetilde{\phi}(x))(1,d_i,1)=\phi((1,d_i,1)x)=\begin{cases}
				 x& i=1,\\ 0& i\neq 1.
				\end{cases}
 \]
 If we put $x_j=(1,d_j,1)^{-1}x\in H_{\mathrm{LT}}^{K_m}$, we have
 \[
  (\widetilde{\phi}(x_j))(1,d_i,1)=\phi((1,d_id_j^{-1},1)x)=\begin{cases}
				 x& i=j,\\ 0& i\neq j.
				\end{cases}
 \]
 This means that the composite $H_{\mathrm{LT}}^{K_m}\xrightarrow{\widetilde{\phi}}(\Ind_{G^1}^G H'_{\mathrm{LT}})^{K_m}\hookrightarrow \bigoplus_{1\le i\le r}H^{n-1}_c(M^{(0)}_{m,C},\overline{\Q}_\ell)$ is surjective.
 Therefore we conclude that $(\Ind_{G^1}^G H'_{\mathrm{LT}})^{K_m}\cong \bigoplus_{1\le i\le r}H^{n-1}_c(M^{(0)}_{m,C},\overline{\Q}_\ell)$, and $\widetilde{\phi}$ is surjective.

 Now we have only to compare the dimensions of $H_{\mathrm{LT}}^{K_m}$ and
 $(\Ind_{G^1}^G H'_{\mathrm{LT}})^{K_m}$.
 As $M_m^{(\delta)}\cong M_m^{(0)}$ (an isomorphism is given by an element $d\in D^\times$
 with $v_F(\Nrd(d))=-\delta$) and $M_m^{(0)}\otimes_{\breve{F}_m,\sigma}C\cong M_{m,C}^{(0)}$
 for $\sigma\in \Gal(\breve{F}_m/\breve{F})$, we obtain
 \[
  \dim H_{\mathrm{LT}}^{K_m}=n\dim H^{n-1}_c(M_m^{(0)}\otimes_{\breve{F}}C,\overline{\Q}_\ell)=n[\breve{F}_m:\breve{F}]\dim H^{n-1}_c(M_{m,C}^{(0)},\overline{\Q}_\ell).
 \]
 On the other hand, we have
 $\dim (\Ind_{G^1}^G H'_{\mathrm{LT}})^{K_m}=r\dim H^{n-1}_c(M_{m,C}^{(0)},\overline{\Q}_\ell)$.
 Since $r=n(\mathcal{O}_F^\times:1+\mathfrak{p}_F^m)=n[\breve{F}_m:\breve{F}]$, we conclude that
 $\dim H_{\mathrm{LT}}^{K_m}=\dim (\Ind_{G^1}^G H'_{\mathrm{LT}})^{K_m}$.
\end{prf}

Next we introduce the infinite level setting.
Let $A_\infty$ be the $\mathfrak{m}_0$-adic completion of $\varinjlim_m A_m$, where
$\mathfrak{m}_0$ denotes the maximal ideal of $A_0$.
Put $A_{m,\mathcal{O}_C}=A_m\mathbin{\widehat{\otimes}}_{\mathcal{O}_{\breve{F}_m}}\mathcal{O}_C$ and
$A_{\infty,\mathcal{O}_C}=A_\infty\mathbin{\widehat{\otimes}}_{\mathcal{O}_{\widehat{F}^{\mathrm{ab}}}}\mathcal{O}_C$.
We can observe that $A_{\infty,\mathcal{O}_C}$ is the $\mathfrak{m}_0$-adic
completion of $\varinjlim_m A_{m,\mathcal{O}_C}$.
The group $G^1$ acts on $A_{\infty,\mathcal{O}_C}$.
The action of $(g,d,\sigma)\in G^1$ preserves $\mathcal{O}_C\subset A_{\infty,\mathcal{O}_C}$
and induces $\sigma$ on $\mathcal{O}_C$.
It is known that the ring $A_{\infty,\mathcal{O}_C}$ with the $G^1$-action has very simple description;
see \cite[\S 2]{Weinstein-stabred} and \cite[\S 1]{Imai-Tsushima-tame}.

We put $\mathcal{M}^{(0)}_{\infty,\mathcal{O}_C}=\Spf A_{\infty,\mathcal{O}_C}$,
$t(\mathcal{M}^{(0)}_{\infty,\mathcal{O}_C})=\Spa (A_{\infty,\mathcal{O}_C},A_{\infty,\mathcal{O}_C})$,
and write $M^{(0)}_{\infty,C}$ for the generic fiber (i.e., the open subspace given by the condition $\lvert \varpi(x)\rvert\neq 0$)
of $t(\mathcal{M}^{(0)}_{\infty,\mathcal{O}_C})$.
By \cite{MR1207303} and \cite{MR1306024}, $t(\mathcal{M}^{(0)}_{\infty,\mathcal{O}_C})$ and
$M^{(0)}_{\infty,C}$ are topological spaces equipped with structure presheaves.
In fact $M^{(0)}_{\infty,C}$ is known to be a perfectoid space (\cite[Lemma 2.10.1]{Weinstein-stabred})
and its structure presheaf is a sheaf, but we do not need this result.
Since $\mathcal{M}^{(0)}_{m,\mathcal{O}_C}=\Spf A_{m,\mathcal{O}_C}$
and $M^{(0)}_{m,C}$ is the generic fiber of $t(\mathcal{M}^{(0)}_{m,\mathcal{O}_C})=\Spa(A_{m,\mathcal{O}_C},A_{m,\mathcal{O}_C})$,
we have a morphism $M^{(0)}_{\infty,C}\to M^{(0)}_{m,C}$
compatible with the transition maps of $\{M^{(0)}_{m,C}\}_m$.
Hence we obtain a morphism $M^{(0)}_{\infty,C}\to \{M^{(0)}_{K,C}\}_{K\subset \SL_n(\mathcal{O}_F)}$,
which is clearly $G^1$-equivariant.

The following proposition enables us to apply the results in Section \ref{sec:infinite-level}
to an affinoid of $M^{(0)}_{\infty,C}$.

\begin{prop}\label{prop:LT-tower-affinoid}
 Let $U$ be a rational subset of $t(\mathcal{M}^{(0)}_{\infty,\mathcal{O}_C})$ contained in $M^{(0)}_{\infty,C}$,
 and $A$ a $\varpi$-adically complete flat $\mathcal{O}_C$-algebra 
 such that $\mathcal{X}=\Spf A$ is a formal model of $U$;
 in other words, $(A[1/\varpi],\widetilde{A})$ is isomorphic to
 $(\mathcal{O}(U),\mathcal{O}^+(U))$ over $(C,\mathcal{O}_C)$ (recall that $\widetilde{A}$ denotes
 the integral closure of $A$ in $A[1/\varpi]$).
 Assume that there exists an open subgroup $K_U$ of $\SL_n(\mathcal{O}_F)$ which stabilizes $U$ and
 whose induced action on $U$ extends to $\mathcal{X}$; namely, the induced action of
 $K_U$ on $\mathcal{O}(U)=A[1/\varpi]$ stabilizes $A$.
 \begin{enumerate}
  \item There exist an integer $m\ge 0$ such that $K'_m\subset K_U$ and
	an open affinoid $U_m$ in $M^{(0)}_{m,C}$ whose inverse image
	under $M^{(0)}_{\infty,C}\to M^{(0)}_{m,C}$ is equal to $U$.
	We put $K_0=K'_m$ and $U_{K_0}=U_m$.
  \item For an open normal subgroup $K$ of $K_0$, let $U_K=\Spa(B_K,B_K^\circ)$ denote the affinoid
	obtained as the inverse image of $U_{K_0}$ under $M^{(0)}_{K,C}\to M^{(0)}_{K_0,C}$.
	Then, $(A,\{B_K\}_{K\lhd K_0})$ satisfies Assumption \ref{assump:Hecke}.
  \item Assume that $g\in G^1$ satisfies $Ug=U$ and extends to an automorphism of $\mathcal{X}$.
	Then, the induced automorphism of $A$ is of finite level in the sense of
	Definition \ref{defn:finite-level-aut}.
 \end{enumerate}
\end{prop}

\begin{prf}
 We write $U=R(\frac{t_1,\ldots,t_r}{s})$, where $t_1,\ldots,t_r,s$ are elements of $A_{\infty,\mathcal{O}_C}$
 such that $t_1,\ldots,t_r$ generate an open ideal of $A_{\infty,\mathcal{O}_C}$.
 As $A_{\infty,\mathcal{O}_C}$ is the completion of $\varinjlim_m A_{m,\mathcal{O}_C}$,
 \cite[Lemma 3.10]{MR1207303} tells us that we may assume that $t_1,\ldots,t_r,s\in A_{m,\mathcal{O}_C}$
 for some $m\ge 0$. We take such $m$ so that $K_m'\subset K_U$.
 Since $\lvert\varpi(x)\rvert\neq 0$ for $x\in U$, by \cite[Lemma 3.11]{MR1207303} there exists
 an integer $j\ge 0$ such that $\lvert a(x)\rvert<\lvert \varpi(x)\rvert$ for every
 $a\in \mathfrak{m}_0^j A_{\infty,\mathcal{O}_C}$ and $x\in U$.
 Take generators $b_1,\ldots,b_l\in \mathfrak{m}_0^j$ and 
 consider the rational subset $U_m=R(\frac{t_1,\ldots,t_r}{s},\frac{b_1,\ldots,b_l}{\varpi})$
 of $t(\mathcal{M}^{(0)}_{m,\mathcal{O}_C})$.
 Then, it is contained in the generic fiber $M^{(0)}_{m,C}$ and its inverse image under
 $M^{(0)}_{\infty,C}\to M^{(0)}_{m,C}$ is equal to $U$.
 This proves (i).
 
 Next we prove (ii). Since $K_0=K_m'\subset K_U$, $A$ is equipped with an action of $K_0$.
 The conditions (a) and (b) in Assumption \ref{assump:Hecke} follow from
 the corresponding properties for $\{M^{(0)}_{K,C}\}_{K\lhd K_0}$.
 Consider the condition (c).
 For $m'\ge m$, we simply write $U_{m'}$ and $B_{m'}$ for $U_{K'_{m'}}$ and $B_{K'_{m'}}$, respectively.
 As the map $t(\mathcal{M}^{(0)}_{\infty,\mathcal{O}_C})\to t(\mathcal{M}^{(0)}_{m',\mathcal{O}_C})$
 carries $U$ to $U_{m'}$, we have a continuous homomorphism of affinoid rings
 $(B_{m'},B_{m'}^\circ)\to (\mathcal{O}(U),\mathcal{O}^+(U))\cong(A[1/\varpi],\widetilde{A})$
 (see \cite[Proposition 1.3 (i)]{MR1306024}).
 Therefore we obtain a homomorphism $\{(B_{m'},B_{m'}^\circ)\}_{m'\ge m}\to (A[1/\varpi],\widetilde{A})$,
 which is clearly $K_0$-equivariant.
 This extends to a $K_0$-equivariant homomorphism
 $\{(B_K,B_K^\circ)\}_{K\lhd K_0}\to (A[1/\varpi],\widetilde{A})$.
 By \cite[Lemma 1.5 (ii)]{MR1306024}, the induced map
 \[
  \Spa(A[1/\varpi],\widetilde{A})\to \varprojlim_{K\lhd K_0}\Spa(B_K,B_K^\circ)=\varprojlim_{m'\ge m}\Spa(B_{m'},B_{m'}^\circ)
 \]
 is identified with $U\to \varprojlim_{K\lhd K_0} U_K=\varprojlim_{m'\ge m} U_{m'}$.
 On the other hand, by \cite[Proposition 2.4.2]{Scholze-Weinstein},
 the continuous map between underlying topological spaces
 \[
  t(\mathcal{M}^{(0)}_{\infty,\mathcal{O}_C})\to\varprojlim_{m'} t(\mathcal{M}^{(0)}_{m',\mathcal{O}_C})
 \]
 is a homeomorphism.
 Therefore $U\to \varprojlim_{m'\ge m} U_{m'}$ is also a homeomorphism.
 The density of the image of $\varinjlim_K B_K=\varinjlim_m B_m\to A[1/\varpi]$
 follows from the density of the image of $\varinjlim_m A_{m,\mathcal{O}_C}\to A_{\infty,\mathcal{O}_C}$;
 see the proof of \cite[Remark 2.4.3 ii)]{MR1734903}. This completes the proof of (ii).

 We prove (iii). Clearly $g$ preserves $\mathcal{O}_C\subset A$.
 By the definition of the action of $G^1$ on $\{M^{(0)}_{K,C}\}_{K\subset\SL_n(\mathcal{O}_F)}$,
 for an open normal subgroup $K$ of $K_0$ and $g\in G^1$, there exists an open normal subgroup $K_g$ of $K_0$
 such that $g$ gives a morphism $M^{(0)}_{K_g,C}\to M^{(0)}_{K,C}$.
 As the maps $U\to U_K$ and $U\to U_{K_g}$ are surjective (see Lemma \ref{lem:limit-surj}),
 $U_{K_g}$ is mapped to $U_K$ by $g$.
 Therefore a continuous homomorphism $B_K\to B_{K_g}$ is induced. Passing to the inductive limit,
 we conclude that the action of $g$ on $A[1/\varpi]$ preserves the image of
 $\varinjlim_{K\lhd K_0}B_K\to A[1/\varpi]$. Since we are assuming that $g$ preserves $A\subset A[1/\varpi]$,
 this means that the action of $g$ is of finite level.
\end{prf}

\begin{cor}\label{cor:LT-specialization}
 Let $U$ and $\mathcal{X}=\Spf A$ be as in Proposition \ref{prop:LT-tower-affinoid}.
 Let $J$ be a subgroup of $G^1$ whose action on $M^{(0)}_{\infty,C}$ stabilizes $U$ and
 extends to an action on $\mathcal{X}$.
 Assume that there exists an affine scheme $Y$ of finite type over $\overline{\F}_q$ endowed with a $J$-action
 such that $\mathcal{X}_s=\Spec (A\otimes_{\mathcal{O}_C}\overline{\F}_q)$ is isomorphic to the perfection of $Y$
 as schemes over $\overline{\F}_q$ with $J$-actions.
 \begin{enumerate}
  \item  We have a $J$-equivariant homomorphism
	 \[
	 \spp^*\colon H^{n-1}_c(Y,\overline{\Q}_\ell)\to \varinjlim_K H^{n-1}_c(M^{(0)}_{K,C},\overline{\Q}_\ell)=H'_{\mathrm{LT}}.
	 \]
  \item Assume moreover that $Y$ is pure-dimensional and smooth over $\overline{\F}_q$.
	Let $V$ be a subspace of $H^{n-1}_c(Y,\overline{\Q}_\ell)$ such that the composite
	$V\hookrightarrow H^{n-1}_c(Y,\overline{\Q}_\ell)\to H^{n-1}(Y,\overline{\Q}_\ell)$ is injective.
	Then, the composite
	$V\hookrightarrow H^{n-1}_c(Y,\overline{\Q}_\ell)\xrightarrow{\spp^*} H'_{\mathrm{LT}}$	
	is also injective.
  \item In the situation of (ii), suppose that $V$ is stable under $J$.
	Let $J'$ be a subgroup of $G^1$ which contains $J$ as a finite index subgroup.
	Assume that $U$ and $Ug$ are disjoint for every $g\in J'\setminus J$.
	Then, the composite
	\[
	\Ind_J^{J'}V\hookrightarrow \Ind_J^{J'}H^{n-1}_c(Y,\overline{\Q}_\ell)\to H'_{\mathrm{LT}}
	\]
	is injective, where the second arrow denotes the $J'$-homomorphism induced from $\spp^*$
	by the Frobenius reciprocity.
 \end{enumerate}
\end{cor}

\begin{prf}
 Let $K_0$, $\{U_K\}_{K_0\lhd K}$ and $\{B_K\}_{K_0\lhd K}$ be as in Proposition \ref{prop:LT-tower-affinoid}.
 Write $Y=\Spec R$.
 By Proposition \ref{prop:LT-tower-affinoid} (ii), $(A,\{B_K\}_{K\lhd K_0},R)$ satisfies 
 Assumptions \ref{assump:Hecke} and \ref{assump:fin-type-perfection}.
 Therefore, by Definition \ref{defn:specialization-map-infinite-level},
 we have a specialization map
 \[
  \spp^*\colon H^{n-1}_c(Y,\overline{\Q}_\ell)\to \varinjlim_K H^{n-1}_c(U_K,\overline{\Q}_\ell).
 \]
 By Propositions \ref{prop:specialization-map-functoriality} and \ref{prop:LT-tower-affinoid} (iii),
 this map is $J$-equivariant.
 Composing with the natural map 
 $\varinjlim_K H^{n-1}_c(U_K,\overline{\Q}_\ell)\to \varinjlim_K H^{n-1}_c(M^{(0)}_{K,C},\overline{\Q}_\ell)$,
 we obtain a desired $J$-equivariant map.
 The second assertion is a direct consequence of Theorem \ref{thm:injectivity-criterion-infinite-level}.

 Let us prove (iii). Take a system of representatives $g_1,\ldots,g_r$ of $J\backslash J'$ such that
 $g_1=1$. Put $U_i=Ug_i$. Then, $\mathcal{X}=\Spf A$ gives a formal model of $U_i$
 by the isomorphism $(A[1/\varpi],\widetilde{A})\cong (\mathcal{O}(U),\mathcal{O}^+(U))\xrightarrow[\cong]{g_i^{-1}}(\mathcal{O}(U_i),\mathcal{O}^+(U_i))$. We write $\mathcal{X}_i$ for this formal model.
 Then, the pair $(U_i,\mathcal{X}_i)$ satisfies the condition in Proposition \ref{prop:LT-tower-affinoid};
 indeed, we can take $K_{U_i}$ as $g_i^{-1}K_Ug_i$.
 Let $m\ge 0$ be an integer which satisfies the condition in Proposition \ref{prop:LT-tower-affinoid} (i)
 for every $1\le i\le r$. Then, by Proposition \ref{prop:LT-tower-affinoid} (ii),
 we obtain towers $\{U_{i,K}\}_{K\lhd K_0}$. Put $\mathcal{X}'=\coprod_{1\le i\le r}\mathcal{X}_i$ and
 $U'_K=\coprod_{1\le i\le r}U'_{i,K}$.
 The reduction $\mathcal{X}'_s$ is identified with the perfection of $Y'=\coprod_{1\le i\le r}Y_i$, 
 where $Y_i$ is a copy of $Y$.
 By assumption, the natural morphism $U'_K\to M^{(0)}_{K,C}$ is an open immersion.
 Therefore, by Definition \ref{defn:specialization-map-infinite-level} we have homomorphisms
 \begin{align*}
 \bigoplus_{1\le i\le r}H^{n-1}_c(Y_i,\overline{\Q}_\ell)=H^{n-1}_c(Y',\overline{\Q}_\ell)&\xrightarrow{\spp^*}\varinjlim_K H^{n-1}_c(U'_K,\overline{\Q}_\ell)\\
  &\to \varinjlim_K H^{n-1}_c(M^{(0)}_{K,C},\overline{\Q}_\ell)=H'_{\mathrm{LT}},
 \end{align*}
 whose composite is denoted by $\spp'^*$.
 The restriction of $\spp'^*$ to $H^{n-1}_c(Y_1,\overline{\Q}_\ell)$ coincides with $\spp^*$ in (i).
 Furthermore, by Theorem \ref{thm:injectivity-criterion-infinite-level}, the composite
 \[
 \bigoplus_{1\le i\le r}V\hookrightarrow \bigoplus_{1\le i\le r}H^{n-1}_c(Y_i,\overline{\Q}_\ell)\xrightarrow{\spp'^*} H'_{\mathrm{LT}}
 \]
 is injective.

 On the other hand, we can define an action of $J'$ on $\mathcal{X}'$ as follows.
 Take $h\in J'$. For $1\le i\le r$, let $j_i$ be the unique integer such that $g_ih\in Jg_{j_i}$.
 The map $i\mapsto j_i$ is a permutation of $\{1,\ldots,r\}$. 
 Define $h\colon \mathcal{X}_i\to \mathcal{X}_{j_i}$ as
 \[
  \mathcal{X}_i=\mathcal{X}\xrightarrow{g_ihg_{j_i}^{-1}}\mathcal{X}=\mathcal{X}_{j_i},
 \]
 and $h\colon \mathcal{X}'\to \mathcal{X}'$ as the disjoint union of them.
 Similarly, we can define an action of $J'$ on $Y'$.
 Since $\bigcup_{1\le i\le r}U_i$ is stable under $J'$,
 the action of $J'$ on $\{M^{(0)}_{K,C}\}_K$ induces that on $\{U'_K\}_K$, and it is compatible
 with the action of $J'$ on $\mathcal{X}'$ (see the proof of Proposition \ref{prop:LT-tower-affinoid} (iii)).
 Therefore, the action of $J'$ on $\mathcal{X}'$ is of finite level and
 the homomorphism $\spp'^*\colon H^{n-1}_c(Y',\overline{\Q}_\ell)\to H'_{\mathrm{LT}}$ is $J'$-equivariant.

 It is immediate to see that the $J'$-equivariant map
 \[
  \Ind_J^{J'}H^{n-1}_c(Y,\overline{\Q}_\ell)\to H^{n-1}_c(Y',\overline{\Q}_\ell)
 \]
 induced from $H^{n-1}_c(Y_1,\overline{\Q}_\ell)\hookrightarrow H^{n-1}_c(Y',\overline{\Q}_\ell)$
 is an isomorphism. Under this isomorphism,
 $\Ind_J^{J'}V\subset \Ind_J^{J'}H^{n-1}_c(Y,\overline{\Q}_\ell)$ is mapped to
 the subspace $\bigoplus_{1\le i\le r}V\subset \bigoplus_{1\le i\le r}H^{n-1}_c(Y_i,\overline{\Q}_\ell)=H^{n-1}_c(Y',\overline{\Q}_\ell)$.
 The composite
 \[
  \Ind_J^{J'}H^{n-1}_c(Y,\overline{\Q}_\ell)\cong H^{n-1}_c(Y',\overline{\Q}_\ell)\xrightarrow{\spp'^*}H'_{\mathrm{LT}}
 \]
 is clearly the $J'$-homomorphism induced from $\spp^*$ by the Frobenius reciprocity.
 Therefore, we conclude that the composite
 \[
 \Ind_J^{J'}V\hookrightarrow \Ind_J^{J'}H^{n-1}_c(Y,\overline{\Q}_\ell)\to H'_{\mathrm{LT}}
 \]
 is injective.
\end{prf}

We end this section by giving representation-theoretic lemmas
used in the next two sections.

\begin{lem}\label{lem:induction-easy}
 Let $G$ be a totally disconnected locally compact group and $J$ an open subgroup of $G$.
 \begin{enumerate}
  \item For a smooth representation $\pi$ of $J$, there exists a $J$-equivariant injection
	$\pi\hookrightarrow \Ind_J^G\pi$.
  \item Let $N$ be a closed subgroup of $G$ satisfying $NJ=JN$,
	and $\rho$ a smooth representation of $N$.
	Then, we have a $J$-equivariant isomorphism 
	$(\Ind_N^{NJ}\rho)\vert_J\cong \Ind_{N\cap J}^J (\rho\vert_{N\cap J})$
	(note that $NJ$ is an open subgroup of $G$).
 \end{enumerate}
\end{lem}

\begin{prf}
 We write $V$ for the representation space of $\pi$.
 For $v\in V$, define the map $\phi_v\colon G\to V$ as follows:
 \[
  \phi_v(g)=\begin{cases}
	     gv& g\in J,\\ 0& g\notin J.
	    \end{cases}
 \]
 It is immediate to see that $\phi_v\in \Ind_J^GV$ and $v\mapsto \phi_v$ gives a $J$-equivariant injection
 $\pi\hookrightarrow \Ind_J^G\pi$. This concludes the proof of (i).

 Next consider (ii). 
 We denote by $\algInd$ the algebraic induction functor.
 Then, we can easily check that 
 \[
  (\algInd_N^{NJ}\rho)\vert_J\cong \algInd_{N\cap J}^J(\rho\vert_{N\cap J}).
 \]
 Indeed, the map $\algInd_N^{NJ}\rho\to \algInd_{N\cap J}^J\rho$ is given by $\phi\mapsto \phi\vert_J$,
 and the map $\algInd_{N\cap J}^J\rho\to \algInd_N^{NJ}\rho$ by $\psi\mapsto \psi'$, where
 $\psi'(ng)=n\psi(g)$ for $n\in N$ and $g\in J$. These maps are inverse to each other.
 By taking the $J$-smooth part $(-)^{\text{$J$-sm}}$, we have
 \[
  \bigl((\algInd_N^{NJ}\rho)\vert_J\bigr)^{\text{$J$-sm}}\cong \Ind_{N\cap J}^J(\rho\vert_{N\cap J}).
 \]
 Since $J$ is an open subgroup of $NJ$, for any representation $W$ of $NJ$,
 we have $W^{\text{$J$-sm}}=W^{\text{$NJ$-sm}}$.
 Therefore, we have $(\algInd_N^{NJ}\rho)^{\text{$J$-sm}}=\Ind_N^{NJ}\rho$. This concludes the proof of (ii).
\end{prf}

\begin{lem}\label{lem:representation-theory}
 Let $G$ and $H$ be totally disconnected locally compact groups and $\Delta$ a closed subgroup of $G\times H$.
 Assume that the composite $\Delta\hookrightarrow G\times H\twoheadrightarrow G$ is surjective.
 Let $(\pi,V)$ (resp.\ $(\rho,W)$) be a smooth representation of $\Delta$ (resp.\ $H$).
 \begin{enumerate}
  \item The group $G$ acts on the space $\Hom_{\Delta\cap H}(\rho,\pi)$ as follows:
	for $g\in G$ and $\phi\in\Hom_{\Delta\cap H}(\rho,\pi)$,
	\[
	(g\phi)(w)=(g,h_g)\phi(h_g^{-1}w)\quad (w\in W).
	\]
	Here $h_g$ denotes an element of $H$ such that $(g,h_g)\in \Delta$, whose existence follows from
	the assumption on $\Delta$. Note that the right hand side is independent of the choice of $h_g$.
  \item If $\rho$ is finitely generated as a representation of $H$, we have a natural isomorphism
	$\Hom_{\Delta\cap H}(\rho,\pi)^{\text{$G$-}\mathrm{sm}}\cong \Hom_H(\rho,\Ind_\Delta^{G\times H}\pi)$
	as representations of $G$,
	where $(-)^{\text{$G$-}\mathrm{sm}}$ denotes the $G$-smooth part.
  \item Assume that $\rho$ is finitely generated as a representation of $H$.
	For a closed subgroup $\Delta'$ of $G\times H$ containing $\Delta$, we have
	\[
	\Hom_{\Delta\cap H}(\rho,\pi)^{\text{$G$-}\mathrm{sm}}\cong \Hom_{\Delta'\cap H}(\rho,\Ind_\Delta^{\Delta'}\pi)^{\text{$G$-}\mathrm{sm}}
	\]
	as representations of $G$.
 \end{enumerate}
\end{lem}

\begin{prf}
 As in the proof of Lemma \ref{lem:induction-easy} (ii), we denote by $\algInd$ the algebraic induction functor.
 First we prove the following:
 \begin{equation*}
  \Hom_{\Delta\cap H}(\rho,\pi)\cong \Hom_H(\rho,\algInd_\Delta^{G\times H}\pi).\tag{$*$}
 \end{equation*}
 By the surjectivity of $\Delta\hookrightarrow G\times H\twoheadrightarrow G$, we have $\Delta H=H\Delta=G\times H$.
 Therefore, as in the proof of Lemma \ref{lem:induction-easy} (ii), we can prove that
 \[
  (\algInd_\Delta^{G\times H}\pi)\vert_H\cong \algInd_{\Delta\cap H}^H(\pi\vert_{\Delta\cap H}).
 \]
 Therefore, by the Frobenius reciprocity we have
 \[
  \Hom_{\Delta\cap H}(\rho,\pi)\cong \Hom_H(\rho,\algInd_{\Delta\cap H}^H\pi)\cong \Hom_H(\rho,\algInd_\Delta^{G\times H}\pi).
 \]
 This concludes the proof of $(*)$.
 Under the isomorphism $(*)$, $\phi\in \Hom_{\Delta\cap H}(\rho,\pi)$ is mapped to
 $\Phi\colon W\to \algInd_\Delta^{G\times H}V$ as follows:
 \[
  \Phi(w)(g,h)=(g,h_g)\phi(h_g^{-1}hw)\quad (w\in W).
 \]
 The inverse map $\Phi\mapsto \phi$ is given by $\phi(w)=\Phi(w)(1,1)$ $(w\in W)$.
 
 Note that $G$ acts on $\Hom_H(\rho,\algInd_\Delta^{G\times H}\pi)$. 
 Therefore, by using the isomorphism $(*)$, we can define an action of $G$ on $\Hom_{\Delta\cap H}(\rho,\pi)$.
 Let us describe it explicitly.
 For $g\in G$, $w\in W$ and $\phi\in \Hom_{\Delta\cap H}(\rho,\pi)$ corresponding to
 $\Phi\in \Hom_H(\rho,\algInd_\Delta^{G\times H}\pi)$, we have
 \begin{align*}
  (g\phi)(w)&=\bigl((g\Phi)(w)\bigr)(1,1)=\bigl((g,1)\Phi(w)\bigr)(1,1)=\bigl((g,h_g)(1,h_g^{-1})\Phi(w)\bigr)(1,1)\\
  &=\bigl((g,h_g)\Phi(h_g^{-1}w)\bigr)(1,1)=\Phi(h_g^{-1}w)(g,h_g)=(g,h_g)\bigl(\Phi(h_g^{-1}w)(1,1)\bigr)\\
  &=(g,h_g)\phi(h_g^{-1}w).
 \end{align*}
 This action coincides with that in (i). This concludes the proof of (i).

 Next we prove (ii). By $(*)$ and the smoothness of $\rho$, we have an isomorphism
 \[
 \Hom_{\Delta\cap H}(\rho,\pi)\cong \Hom_H\bigl(\rho,(\algInd_\Delta^{G\times H}\pi)^{\text{$H$-sm}}\bigr),
 \]
 which is $G$-equivariant by the proof of (i).
 By taking the $G$-smooth parts, we have
 \[
 \Hom_{\Delta\cap H}(\rho,\pi)^{\text{$G$-sm}}\cong \Hom_H\bigl(\rho,(\algInd_\Delta^{G\times H}\pi)^{\text{$H$-sm}}\bigr)^{\text{$G$-sm}}.
 \]
 Since $\Ind_{\Delta}^{G\times H}\pi=(\algInd_{\Delta}^{G\times H}\pi)^{\text{$(G\times H)$-sm}}=((\algInd_{\Delta}^{G\times H}\pi)^{\text{$H$-sm}})^{\text{$G$-sm}}$, it suffices to show the following:
 \begin{quote}
  for a $G\times H$-representation $U$,
  we have $\Hom_H(\rho,U)^{\text{$G$-sm}}=\Hom_H(\rho,U^{\text{$G$-sm}})$.
 \end{quote}
 Clearly we have $\Hom_H(\rho,U)^{\text{$G$-sm}}\subset \Hom_H(\rho,U^{\text{$G$-sm}})$.
 Take $\phi\in \Hom_H(\rho,U^{\text{$G$-sm}})$. Since $\rho$ is finitely generated,
 we can find a finite system of generators $x_1,\ldots,x_n\in W$. 
 There exists an open subgroup $K$ of $G$ which stabilizes $\phi(x_1),\ldots,\phi(x_n)\in U^{\text{$G$-sm}}$.
 Then, for every $x\in W$, $\phi(x)$ is fixed by $K$. 
 In other words, $\phi$ lies in $\Hom_H(\rho,U)^K$. 
 This proves that $\Hom_H(\rho,U^{\text{$G$-sm}})\subset \Hom_H(\rho,U)^{\text{$G$-sm}}$.
 Now the proof of (ii) is complete.

 For (iii), just note that
 \begin{align*}
  \Hom_{\Delta\cap H}(\rho,\pi)^{\text{$G$-sm}}&\cong\Hom_H(\rho,\Ind^{G\times H}_\Delta\pi)
 \cong \Hom_H(\rho,\Ind^{G\times H}_{\Delta'}\Ind^{\Delta'}_\Delta\pi)\\
  &\cong \Hom_{\Delta'\cap H}(\rho,\Ind^{\Delta'}_\Delta\pi)^{\text{$G$-sm}}.
 \end{align*}
\end{prf}

\section{Example I: depth 0 supercuspidal representations}\label{sec:depth-0}
In this section, we determine $\rec_F(\pi)$ and $\JL(\pi)$ for a depth $0$ supercuspidal representation
$\pi$ of $\GL_n(F)$. First, we will introduce some notation.
Let $\theta\colon \F_{q^n}^\times\to \C^\times$ be a character. We assume that $\theta$ is regular, namely,
it does not factor through the norm map $\Nr_{\F_{q^n}/\F_{q^m}}\colon \F_{q^n}^\times\to \F_{q^m}^\times$
to any subfield $\F_{q^m}\subsetneq \F_{q^n}$.
It is well-known that we can attach to $\theta$ an irreducible cuspidal representation $R_\theta$
of $\GL_n(\F_q)$. Here we recall a construction of $R_\theta$ by the Deligne-Lusztig theory \cite{MR0393266}.
Let $\mathrm{DL}_n$ be the affine algebraic variety over $\overline{\F}_q$ defined by the equation
\[
 \bigl(\det(x_i^{q^{j-1}})_{1\le i,j\le n}\bigr)^{q-1}=(-1)^{n-1}.
\]
The group $\GL_n(\F_q)$ naturally acts on $\mathrm{DL}_n$ on the right.
On the other hand, the group $\F_{q^n}^\times$ acts on $\mathrm{DL}_n$ (on the right) by
$x_i\mapsto \zeta x_i$ $(\zeta\in \F_{q^n}^\times)$.
The representation $R_\theta$ is given as the $\theta$-isotypic part
$H^{n-1}_c(\mathrm{DL}_n,\overline{\Q}_\ell)_\theta$ of the $(n-1)$th compactly supported
cohomology $H^{n-1}_c(\mathrm{DL}_n,\overline{\Q}_\ell)$ of $\mathrm{DL}_n$.
The central character of $R_\theta$ equals $\theta\vert_{\F_q^\times}$
(\cite[Corollary 1.22]{MR0393266}).

We denote by $F_n$ the degree $n$ unramified extension of $F$.
Let $\chi\colon F_n^\times\to \C^\times$ be a tame character, namely, a character
which is trivial on $1+\varpi\mathcal{O}_{F_n}$. Its restriction to $\mathcal{O}_{F_n}^\times$
factors as $\mathcal{O}_{F_n}^\times\twoheadrightarrow \F_{q^n}^\times\xrightarrow{\overline{\chi}}\C^\times$.
We say that $\chi$ is regular if $\overline{\chi}$ is regular.

\begin{defn}\label{defn:GL-depth-0}
 For a regular tame character $\chi\colon F_n^\times\to\C^\times$,
 let $\overline{\pi}_\chi$ be the representation of $F^\times\GL_n(\mathcal{O}_F)$ such that
 \begin{itemize}
  \item $\overline{\pi}_\chi\vert_{\GL_n(\mathcal{O}_F)}$ is the inflation of $R_{\overline{\chi}}$
	by $\GL_n(\mathcal{O}_F)\twoheadrightarrow \GL_n(\F_q)$, and
  \item $F^\times$ acts on $\overline{\pi}_\chi$ by $\chi\vert_{F^\times}$.
 \end{itemize}
 We put $\pi_\chi=\cInd_{F^\times\GL_n(\mathcal{O}_F)}^{\GL_n(F)}\overline{\pi}_\chi$,
 which is known to be an irreducible supercuspidal representation of $\GL_n(F)$.
\end{defn}

From a regular tame character $\chi$, we can also construct a representation of $D^\times$.
Recall that the quotient of $\mathcal{O}_D$ by its Jacobson radical is equal to $\F_{q^n}$.
Therefore we have a surjection $\mathcal{O}_D^\times\twoheadrightarrow \F_{q^n}^\times$.

\begin{defn}\label{defn:D-depth-0}
 For a regular tame character $\chi\colon F_n^\times\to\C^\times$, let $\overline{\rho}_\chi$
 be the character of $F^\times\mathcal{O}_D^\times$ such that
 \begin{itemize}
  \item $\overline{\rho}_\chi\vert_{\mathcal{O}_D^\times}$ is the composite
	$\mathcal{O}_D^\times\twoheadrightarrow \F_{q^n}^\times\xrightarrow{\overline{\chi}}\C^\times$,	and
  \item $\overline{\rho}_\chi\vert_{F^\times}=\chi\vert_{F^\times}$.
 \end{itemize}
 We put $\rho_\chi=\cInd_{F^\times\mathcal{O}_D^\times}^{D^\times}\overline{\rho}_\chi$,
 which is known to be an irreducible smooth representation of $D^\times$.
\end{defn}

The goal of this section is the following:

\begin{thm}\label{thm:LLC-JLJC-depth-0}
 For a regular tame character $\chi\colon F_n^\times\to \C^\times$, we have
 \[
  \rec_F(\pi_\chi)=\Ind_{W_{F_n}}^{W_F}\chi\delta^{n-1},\quad \JL(\pi_\chi)=\rho_\chi,
 \]
 where $\delta\colon F_n^\times\to \C^\times$; $a\mapsto (-1)^{v_{F_n}(a)}$
 is the unramified quadratic character. As usual, we regard characters of $F_n^\times$ as those of $W_{F_n}$
 by the isomorphism $\Art_{F_n}\colon F_n^\times\xrightarrow{\cong}W_{F_n}^{\mathrm{ab}}$.
\end{thm}

\begin{rem}\label{rem:history-depth-0}
 \begin{enumerate}
  \item The identity $\rec_F(\pi_\chi)=\Ind_{W_{F_n}}^{W_F}\chi\delta^{n-1}$ was proved in
	\cite{MR1235293}. The proof in \cite{MR1235293} requires case-by-case argument when
	$(n,q)$ is either of $(2,2)$, $(2,3)$, $(4,2)$, $(6,2)$.
	Our proof works for any case uniformly.
  \item The identity $\rec_F(\pi_\chi)=\rho_\chi$ was obtained in \cite[Theorem 3]{MR2164626};
	see also \cite[Theorem 1]{MR2736858}.
 \end{enumerate}
\end{rem}

\begin{rem}\label{rem:unram-character-twist}
 Let $\xi\colon F^\times\to \C^\times$ be an unramified character.
 For a tame character $\chi\colon F_n^\times\to \C^\times$, put $\chi_\xi=\chi(\xi\circ\Nr_{F_n/F})$.
 Then, $\chi_\xi$ is also tame and $\overline{\chi_\xi}=\overline{\chi}$.
 Therefore, $\chi$ is regular if and only if $\chi_\xi$ is regular.
 Moreover, we have $\overline{\pi}_{\chi_\xi}=\overline{\pi}_{\chi}\otimes(\xi\circ\det)$ and
 \[
  \pi_\chi\otimes(\xi\circ\det)=\cInd_{F^\times\GL_n(\mathcal{O}_F)}^{\GL_n(F)}(\overline{\pi}_\chi\otimes (\xi\circ\det))=\pi_{\chi_\xi}.
 \]
 Similarly we have $\rho_\chi\otimes(\xi\circ\Nrd)=\rho_{\chi_\xi}$ and 
 $(\Ind_{W_{F_n}}^{W_F}\chi)\otimes \xi=\Ind_{W_{F_n}}^{W_F}\chi_\xi$.
 As $\rec_F$ and $\JL$ are compatible with character twist, Theorem \ref{thm:LLC-JLJC-depth-0} for $\chi$
 is equivalent to that for $\chi_\xi$.

 If we take $\xi$ such that $\xi(\varpi)^n=\chi(\varpi)^{-1}$, we have $\chi_\xi(\varpi)=1$.
 Therefore, to prove Theorem \ref{thm:LLC-JLJC-depth-0} we may assume that $\chi(\varpi)=1$.
 In this case the central characters of $\pi_\chi$ and $\rho_\chi$ are trivial on $\varpi^\Z\subset F^\times$.

 Later, for a regular character $\theta\colon \F_{q^n}^\times\to \C^\times$, we also write $\theta$
 for the regular tame character $F_n^\times=\varpi^\Z\times \mathcal{O}_{F_n}^\times\twoheadrightarrow \mathcal{O}_{F_n}^\times\twoheadrightarrow\F_{q^n}^\times\xrightarrow{\theta}\C^\times$.
 It suffices to show that $\rec_F(\pi_\theta)=\Ind_{W_{F_n}}^{W_F}\theta\delta^{n-1}$
 and $\JL(\pi_\theta)=\rho_\theta$.
\end{rem}

For $b\in\overline{\F}_q$ with $b^{q-1}=(-1)^{n-1}$, let $Y_b$ be
the affine algebraic variety over $\overline{\F}_q$ defined by the equation
$\det(x_i^{q^{j-1}})_{1\le i,j\le n}=b$. Clearly we have $\mathrm{DL}_n=\coprod_{b^{q-1}=(-1)^{n-1}}Y_b$.
In the following, we construct an affinoid in $M^{(0)}_{\infty,C}$ and its formal model
whose reduction is isomorphic to $Y_b$ for some $b\in\overline{\F}_q$ with $b^{q-1}=(-1)^{n-1}$.
Since it is well-known to specialists, we omit the detail.

Let $\widetilde{\mathbb{X}}$ be the formal $\mathcal{O}_F$-module over $\mathcal{O}_F$ whose logarithm
is $\sum_{i=0}^\infty\frac{T^{q^{in}}}{\varpi^i}$ (it is called the standard formal $\mathcal{O}_F$-module
in \cite[2.3]{MR3402698}). 
In the following, take $\mathbb{X}$ as $\widetilde{\mathbb{X}}\otimes_{\mathcal{O}_F}\overline{\F}_q$.
Let $y=(\xi^{q^{-m}})_{m\ge 0}\in T\widetilde{\mathbb{X}}(\mathcal{O}_C)$,
$\alpha_1,\ldots,\alpha_n\in\mathcal{O}_{F_n}$ and $t=(\tau^{q^{-m}})_{m\ge 0}$
be as in \cite[3.7]{MR3402698}.
Since $\widetilde{\mathbb{X}}\otimes_{\mathcal{O}_F}\mathcal{O}_C$ has CM by $\mathcal{O}_{F_n}$,
we have embeddings $F_n\hookrightarrow M_n(F)$ and $F_n\hookrightarrow D$,
by which we regard $F_n$ as subfields of $M_n(F)$ and $D$ (see \cite[3.1]{MR3402698}).

By \cite[2.9.2]{Weinstein-stabred}, we have
\[
 A_{\infty,\mathcal{O}_C}\cong \mathcal{O}_C[[X_1^{q^{-\infty}},\ldots,X_n^{q^{-\infty}}]]/(\delta(X_1,\ldots,X_n)^{q^{-m}}-\tau^{q^{-m}}\mid m\ge 0)^-,
\]
where $(-)^-$ denotes the closure.
For the definition of $\delta^{q^{-m}}$, see \cite[\S 1.1]{Imai-Tsushima-tame}.

Let $U$ be the subset of $M^{(0)}_{\infty,C}$
defined by $\lvert X_i\rvert\le \lvert\varpi\rvert^{\frac{1}{q^n-1}}$ for every $1\le i\le n$.
It is clearly a rational subset of $t(\mathcal{M}^{(0)}_{\infty,\mathcal{O}_C})$.
Recall that $\lvert\xi\rvert=\lvert\varpi\rvert^{\frac{1}{q^n-1}}$ (\cite[3.7]{MR3402698}),
and put $x^{q^{-m}}_i=X^{q^{-m}}_i/\xi^{q^{-m}}$.
Then, 
\[
 \delta'(\boldsymbol{x}_1,\ldots,\boldsymbol{x}_n)^{q^{-m}}=(\xi^{q^{-m}})^{-(1+q+\cdots+q^{n-1})}\delta(\boldsymbol{\xi}\boldsymbol{x}_1,\ldots,\boldsymbol{\xi}\boldsymbol{x}_n)^{q^{-m}}
\]
lies in $\mathcal{O}_C\langle x^{q^{-\infty}}_1,\ldots,x^{q^{-\infty}}_n\rangle$,
where $\boldsymbol{x}_i$ (resp.\ $\boldsymbol{\xi}$) denotes the system $(x_i^{q^{-m}})$ (resp.\ $(\xi^{q^{-m}})$).
On the other hand, since $\lvert \tau\rvert=\lvert\varpi\rvert^{\frac{1}{q-1}}=\lvert\xi\rvert^{1+q+\cdots+q^{n-1}}$,
\[
 \tau'^{q^{-m}}=(\xi^{q^{-m}})^{-(1+q+\cdots+q^{n-1})}\tau^{q^{-m}}
\]
lies in $\mathcal{O}_C^\times$. Put 
\[
 \mathcal{X}=\Spf \mathcal{O}_C\langle x^{q^{-\infty}}_1,\ldots,x^{q^{-\infty}}_n\rangle/(\delta'(x_1,\ldots,x_n)^{q^{-m}}-\tau'^{q^{-m}}\mid m\ge 0)^-.
\]
As in \cite[Theorem 2.5]{Imai-Tsushima-tame}, we can check that $\mathcal{X}$ gives a formal model of $U$.

\begin{prop}\label{prop:depth-0-formal-model}
 The formal scheme $\mathcal{X}$ is flat over $\mathcal{O}_C$ and its reduction $\mathcal{X}_s$ is
 isomorphic to the perfection of $Y_b$ for some $b\in \overline{\F}_q$ with $b^{q-1}=(-1)^{n-1}$.
\end{prop}

\begin{prf}
 First we compute the reduction $\mathcal{X}_s$.
 By \cite[Lemma 2.10.4]{MR3402698}, the image of $\delta'^{q^{-m}}$ in 
 $\overline{\F}_q[x^{q^{-\infty}}_1,\ldots,x^{q^{-\infty}}_n]$ is equal to 
 $(\det(x_i^{q^{j-1}})_{1\le i,j\le n})^{q^{-m}}$.
 On the other hand, \cite[Lemma 3.7.1]{MR3402698} tells us that the image of $\tau'^{q^{-m}}$
 in $\overline{\F}_q$ is equal to
 $(\det (\overline{\alpha}_i^{q^{j-1}}))^{q^{-m}}$.
 Put $b=\det (\overline{\alpha}_i^{q^{j-1}})$. Since $\alpha_i\in\mathcal{O}_{F_n}$,
 we have $\overline{\alpha}_i^{q^n}=\overline{\alpha}_i$. Therefore, we have
 \[
  b^q=\det (\overline{\alpha}_i^{q^j})=(-1)^{n-1}b.
 \]
 By definition, $\mathcal{X}_s$ is identified with
 \[
  \Spec \overline{\F}_q[x^{q^{-\infty}}_1,\ldots,x^{q^{-\infty}}_n]/(\det(x_i^{q^{j-1}})^{q^{-m}}-b^{q^{-m}}\mid m\ge 0),
 \]
 which is the perfection of $Y_b$.
 
 In particular, we have $\dim \mathcal{X}_s=\dim Y_b=n-1$. Hence Corollary \ref{cor:flatness-dimension} tells us
 that $\mathcal{X}$ is flat over $\mathcal{O}_C$. This completes the proof.
\end{prf}

Next we consider the group action on $U$ and $\mathcal{X}$.

\begin{defn}\label{defn:depth-0-J}
 We put $J=F^\times\GL_n(\mathcal{O}_F)\times(\mathcal{O}_D^\times/\varpi^\Z)\times W_{F_n}\subset G$
 and $J^1=G^1\cap J$.
 For $\sigma\in W_{F_n}$, we put $n_\sigma=v_{F_n}(\Art_{F_n}^{-1}(\sigma))\in \Z$
 and $u_\sigma=\varpi^{-n_\sigma}\Art_{F_n}^{-1}(\sigma)\in\mathcal{O}_{F_n}^\times$.
 We define a homomorphism $\Theta\colon J\to \GL_n(\F_q)\times\F_{q^n}^\times\times \Z$ by
 \begin{align*}
  (\varpi^mg,d,\sigma)\mapsto (\overline{g},\overline{d^{-1}u^{-1}_\sigma},n_\sigma),
 \end{align*}
 where $m\in\Z$, $g\in\GL_n(\mathcal{O}_F)$, $d\in\mathcal{O}_D^\times$, and $\sigma\in W_{F_n}$.

 As explained before, the group $\GL_n(\F_q)\times\F_{q^n}^\times$ acts on $\mathrm{DL}_n$.
 We can also define an action of $\Z$ on $\mathrm{DL}_n$ over $\F_{q^n}$ (not $\overline{\F}_q$)
 by $m\colon x_i\mapsto x_i,$ $a\mapsto a^{q^{-mn}}$
 $(a\in \overline{\F}_q)$,
 which commutes with the action of $\GL_n(\F_q)\times\F_{q^n}^\times$. 
 Therefore, $J$ acts on $\mathrm{DL}_n$ through the homomorphism $\Theta$.

 Recall that we have a decomposition $\mathrm{DL}_n=\coprod_{b^{q-1}=(-1)^{n-1}}Y_b$.
 It is easy to see that an element $(g,a)\in \GL_n(\F_q)\times\F_{q^n}^\times$ maps $Y_b$ to $Y_{b'}$,
 where $b'=\det(g)\Nr_{\F_{q^n}/\F_q}(a)b$. In particular, the group
 \[
  (\GL_n(\F_q)\times \F_{q^n}^\times)^1=\{(g,a)\in \GL_n(\F_q)\times \F_{q^n}^\times\mid \det(g)\Nr_{\F_{q^n}/\F_q}(a)=1\}
 \]
 acts on $Y_b$. On the other hand, if $b^q=(-1)^{n-1}b$, we have
 \[
  b^{q^n}=(-1)^{(n-1)(1+q+\cdots+q^{n-1})}b=b,
 \]
 hence $b\in\F_{q^n}$. Therefore the action of $\Z$ on $\mathrm{DL}_n$ preserves $Y_b$.

 Since the image of $J^1$ under $\Theta$ is contained in
 $(\GL_n(\F_q)\times \F_{q^n})^1\times\Z$, the action of $J^1$ on $\mathrm{DL}_n$ preserves $Y_b$.
\end{defn}

\begin{lem}\label{lem:depth-0-J-generator}
 The group $J^1$ is generated by the following elements:
 \begin{itemize}
  \item $(g,1,1)$ for $g\in \GL_n(\mathcal{O}_F)$ with $\det g=1$,
  \item $(1,d,1)$ for $d\in \mathcal{O}_D^\times$ with $\Nrd d=1$,
  \item $(a,a,1)$ for $a\in F_n^\times$,
  \item $(1,\Art_{F_n}^{-1}(\sigma)^{-1},\sigma)$ for $\sigma\in I_{F_n}$, and
  \item $(1,\varpi^{-1},\sigma)$ for $\sigma\in W_{F_n}$ with $\Art_{F_n}^{-1}(\sigma)=\varpi$.
 \end{itemize}
 Here $I_{F_n}$ denotes the inertia group of $F_n$.
\end{lem}

\begin{prf}
 This is immediate from the surjectivity of $\Nr_{F_n/F}\colon \mathcal{O}_{F_n}^\times\to \mathcal{O}_F^\times$.
\end{prf}

\begin{prop}\label{prop:depth-0-group-action}
 \begin{enumerate}
  \item The action of $J^1\subset G^1$ on $M^{(0)}_{\infty,C}$ stabilizes $U$ and extends to $\mathcal{X}$.
  \item The induced action of $J^1$ on $\mathcal{X}_s$ is compatible with the action of $J^1$
	on $Y_b$ under the isomorphism in Proposition \ref{prop:depth-0-formal-model}.
 \end{enumerate}
\end{prop}

\begin{prf}
 This follows easily from Lemma \ref{lem:depth-0-J-generator} and the description of the group action
 in \cite[\S 1.2]{Imai-Tsushima-tame}.
 For the element $(1,\Art^{-1}_{F_n}(\sigma)^{-1},\sigma)$ with $\sigma\in I_{F_n}$, 
 note that $\lvert\xi\rvert=\lvert\varpi\rvert^{\frac{1}{q^n-1}}$ implies that
 $\xi\in\varpi^{\frac{1}{q^n-1}}\mathcal{O}_C^\times$, hence
 $\sigma(\xi)/\xi\equiv \sigma(\varpi^{\frac{1}{q^n-1}})/\varpi^{\frac{1}{q^n-1}}\equiv\Art_{F_n}^{-1}(\sigma)$
 modulo the maximal ideal of $\mathcal{O}_C$.
 For the element $(1,\varpi^{-1},\sigma)$ with $\Art_{F_n}^{-1}(\sigma)=\varpi$, we use the fact
 $\sigma(\xi)=\xi$, which is a consequence of the classical Lubin-Tate theory.
\end{prf}

\begin{lem}\label{lem:coh-DL}
 Let $\theta\colon \F_{q^n}^\times\to\C^\times$ be a regular character.
 Then, $H^{n-1}_c(\mathrm{DL}_n,\overline{\Q}_\ell)_\theta$ is a $J^1$-stable
 subspace of $H^{n-1}_c(\mathrm{DL}_n,\overline{\Q}_\ell)$. Moreover, we have
 an isomorphism 
 \[
  \Hom_{J^1\cap (\mathcal{O}_D^\times/\varpi^\Z)}\bigl(\overline{\rho}_{\theta^{-1}},H^{n-1}_c(\mathrm{DL}_n,\overline{\Q}_\ell)_\theta\bigr) \cong \overline{\pi}_\theta\boxtimes (\theta\delta^{n-1})^{-1}(\tfrac{1-n}{2})
 \]
 of representations of $F^\times\GL_n(\mathcal{O}_F)\times W_{F_n}$.
 Note that the action of $F^\times\GL_n(\mathcal{O}_F)\times W_{F_n}$ on the left hand side is defined
 as in Lemma \ref{lem:representation-theory} (i), since the composite
 $J^1\hookrightarrow J\twoheadrightarrow F^\times\GL_n(\mathcal{O}_F)\times W_{F_n}$ is surjective.
\end{lem}

\begin{prf}
 Since $J^1$ obviously normalizes $(\mathcal{O}_D^\times/\varpi^\Z,\theta)$,
 $H^{n-1}_c(\mathrm{DL}_n,\overline{\Q}_\ell)_\theta$ is a $J^1$-stable
 subspace of $H^{n-1}_c(\mathrm{DL}_n,\overline{\Q}_\ell)$.
 As $J^1\cap (\mathcal{O}_D^\times/\varpi^\Z)$ acts on $H^{n-1}_c(\mathrm{DL}_n,\overline{\Q}_\ell)_\theta$
 by $\overline{\rho}_{\theta^{-1}}=\theta^{-1}$, as a vector space we have
 \[
  \Hom_{J^1\cap (\mathcal{O}_D^\times/\varpi^\Z)}\bigl(\overline{\rho}_{\theta^{-1}},H^{n-1}_c(\mathrm{DL}_n,\overline{\Q}_\ell)_\theta\bigr)=H^{n-1}_c(\mathrm{DL}_n,\overline{\Q}_\ell)_\theta.
 \]
 Under this identification, the action of $F^\times\GL_n(\mathcal{O}_F)\times W_{F_n}$
 on the left hand side is described as follows:
 \begin{itemize}
  \item $(g,1)$ for $g\in\GL_n(\mathcal{O}_F)$ with $\det g=1$ acts as $\overline{g}\in\GL_n(\F_q)$
	on $H^{n-1}_c(\mathrm{DL}_n,\overline{\Q}_\ell)_\theta$,
  \item $(\varpi,1)$ acts trivially,
  \item $(a,1)$ for $a\in\mathcal{O}_{F_n}^\times$ acts as $\overline{a}\in\GL_n(\F_q)$ on 
	$H^{n-1}_c(\mathrm{DL}_n,\overline{\Q}_\ell)_\theta$,
  \item $(1,\sigma)$ for $\sigma\in I_{F_n}$ acts by the scalar $\theta(\Art_{F_n}^{-1}(\sigma))^{-1}$,
  \item $(1,\sigma)$ for $\sigma\in W_{F_n}$ with $\Art_{F_n}(\sigma)=\varpi$
	acts as $\Frob_q^n$ on $H^{n-1}_c(\mathrm{DL}_n,\overline{\Q}_\ell)_\theta$.
 \end{itemize}
 By \cite{MR840835}, 
 $\Frob_q^n$ acts on $H^{n-1}_c(\mathrm{DL}_n,\overline{\Q}_\ell)_\theta$ by 
 $(-1)^{n-1}q^{\frac{n(n-1)}{2}}=(\delta^{1-n}(\tfrac{1-n}{2}))(\varpi)$
 (see also the proof of \cite[Th\'eor\`eme 3.1.12]{MR3278894}).
 Therefore we conclude that 
 \[
  \Hom_{J\cap (\mathcal{O}_D^\times/\varpi^\Z)}\bigl(\overline{\rho}_{\theta^{-1}},H^{n-1}_c(\mathrm{DL}_n,\overline{\Q}_\ell)_\theta\bigr) \cong \overline{\pi}_\theta\boxtimes (\theta\delta^{n-1})^{-1}(\tfrac{1-n}{2}).
 \]
\end{prf}

\begin{prop}\label{prop:sp-inj-depth-0}
 There exists a $J^1$-equivariant injection
 $H^{n-1}_c(\mathrm{DL}_n,\overline{\Q}_\ell)_\theta\to H'_{\mathrm{LT}}$. 
\end{prop}

\begin{prf}
 Let $b\in\overline{\F}_q$ be as in Proposition \ref{prop:depth-0-formal-model}.
 We put 
 \[
  (\F_{q^n}^\times)^1=\{a\in \F_{q^n}^\times\mid \Nr_{\F_{q^n}/\F_q}(a)=1\},
 \]
 and write $H^{n-1}_c(Y_b,\overline{\Q}_\ell)_\theta$ for the 
 $\theta\vert_{(\F_{q^n}^\times)^1}$-isotypic part of $H^{n-1}_c(Y_b,\overline{\Q}_\ell)$.
 By Propositions \ref{prop:depth-0-formal-model} and \ref{prop:depth-0-group-action},
 the pair $(U,\mathcal{X})$ satisfies the condition in Proposition \ref{prop:LT-tower-affinoid}
 (we may take $K_U=J^1\cap \SL_n(\mathcal{O}_F)$).
 By Proposition \ref{prop:depth-0-group-action} and Corollary \ref{cor:LT-specialization} (i),
 we have a $J^1$-equivariant map
 \[
  \spp^*\colon H^{n-1}_c(Y_b,\overline{\Q}_\ell)\to H'_{\mathrm{LT}}.
 \]
 Now note that $Y_b$ is purely $n-1$-dimensional and smooth over $\overline{\F}_q$.
 Since $Y_b$ is $(\F_{q^n}^\times)^1$-equivariantly isomorphic to the Deligne-Lusztig variety for $\SL_n(\F_q)$,
 \cite[Theorem 9.8]{MR0393266} tells us that the composite
 \[
  H^{n-1}_c(Y_b,\overline{\Q}_\ell)_\theta\hookrightarrow H^{n-1}_c(Y_b,\overline{\Q}_\ell)\to H^{n-1}(Y_b,\overline{\Q}_\ell)
 \]
 is an injection. Therefore, Corollary \ref{cor:LT-specialization} (ii) tells us that
 the composite
 \[
  H^{n-1}_c(Y_b,\overline{\Q}_\ell)_\theta\hookrightarrow H^{n-1}_c(Y_b,\overline{\Q}_\ell)\xrightarrow{\spp^*}H'_{\mathrm{LT}}
 \]
 is an injection.
 
 Thus we have only to prove that the pull-back map 
 $H^{n-1}_c(\mathrm{DL}_n,\overline{\Q}_\ell)\to H^{n-1}_c(Y_b,\overline{\Q}_\ell)$,
 which is $J^1$-equivariant, induces an isomorphism
 $H^{n-1}_c(\mathrm{DL}_n,\overline{\Q}_\ell)_\theta\xrightarrow{\cong} H^{n-1}_c(Y_b,\overline{\Q}_\ell)_\theta$.
 It is easy to see that $H^{n-1}_c(\mathrm{DL}_n,\overline{\Q}_\ell)\cong\Ind_{(\F_{q^n}^\times)^1}^{\F_{q^n}^\times}H^{n-1}_c(Y_b,\overline{\Q}_\ell)$ as representations of $\F_{q^n}^\times$.
 By the Frobenius reciprocity, we have
 \begin{align*}
  H^{n-1}_c(\mathrm{DL}_n,\overline{\Q}_\ell)_\theta&=\Hom_{\F_{q^n}^\times}(\theta,H^{n-1}_c(\mathrm{DL}_n,\overline{\Q}_\ell))\cong\Hom_{\F_{q^n}^\times}\bigl(\theta,\Ind_{(\F_{q^n}^\times)^1}^{\F_{q^n}^\times}H^{n-1}_c(Y_b,\overline{\Q}_\ell)\bigr)\\
  &=\Hom_{(\F_{q^n}^\times)^1}\bigl(\theta\vert_{(\F_{q^n}^\times)^1},H^{n-1}_c(Y_b,\overline{\Q}_\ell)\bigr)
  =H^{n-1}_c(Y_b,\overline{\Q}_\ell)_\theta.
 \end{align*}
 It is easy to see that the composite of the above identifications is equal to the pull-back map
 $H^{n-1}_c(\mathrm{DL}_n,\overline{\Q}_\ell)_\theta\to H^{n-1}_c(Y_b,\overline{\Q}_\ell)_\theta$.
 This concludes the proof.
\end{prf}

\begin{prf}[of Theorem \ref{thm:LLC-JLJC-depth-0}]
 By Proposition \ref{prop:sp-inj-depth-0}, we have an injection
 \[
  \Hom_{J^1\cap (\mathcal{O}_D^\times/\varpi^\Z)}\bigl(\overline{\rho}_{\theta^{-1}},H^{n-1}_c(\mathrm{DL}_n,\overline{\Q}_\ell)_\theta\bigr)\hookrightarrow \Hom_{J^1\cap (\mathcal{O}_D^\times/\varpi^\Z)}(\overline{\rho}_{\theta^{-1}},H'_{\mathrm{LT}}),
 \]
 which is $F^\times\GL_n(\mathcal{O}_F)\times W_{F_n}$-equivariant.
 Lemma \ref{lem:coh-DL} says that 
 \[
  \Hom_{J^1\cap (\mathcal{O}_D^\times/\varpi^\Z)}\bigl(\overline{\rho}_{\theta^{-1}},H^{n-1}_c(\mathrm{DL}_n,\overline{\Q}_\ell)_\theta\bigr)\cong \overline{\pi}_\theta\boxtimes (\theta\delta^{n-1})^{-1}(\tfrac{1-n}{2}),
 \]
 which is a smooth representation of $F^\times\GL_n(\mathcal{O}_F)\times W_{F_n}$.
 Therefore, we obtain an $F^\times\GL_n(\mathcal{O}_F)\times W_{F_n}$-equivariant injection
 \[
  \overline{\pi}_\theta\boxtimes (\theta\delta^{n-1})^{-1}(\tfrac{1-n}{2})
 \hookrightarrow \Hom_{J^1\cap (\mathcal{O}_D^\times/\varpi^\Z)}(\overline{\rho}_{\theta^{-1}},H'_{\mathrm{LT}})^{\mathrm{sm}},
 \]
 where $(-)^{\mathrm{sm}}$ denotes the $F^\times\GL_n(\mathcal{O}_F)\times W_{F_n}$-smooth part.
 On the other hand, Lemma \ref{lem:representation-theory} (ii) tells us that
 \begin{align*}
  \Hom_{J^1\cap (\mathcal{O}_D^\times/\varpi^\Z)}(\overline{\rho}_{\theta^{-1}},H'_{\mathrm{LT}})^{\mathrm{sm}}&\cong \Hom_{\mathcal{O}_D^\times/\varpi^\Z}(\overline{\rho}_{\theta^{-1}},\Ind_{J^1}^J H'_{\mathrm{LT}}).
 \end{align*}
 Since $J\subset G$ is open, $G^1\subset G$ is closed normal and $J^1=G^1\cap J$,
 by Lemma \ref{lem:induction-easy} (ii)
 we have a $J$-equivariant isomorphism $\Ind_{J^1}^J H'_{\mathrm{LT}}\cong \Ind_{G^1}^{G^1J} H'_{\mathrm{LT}}$.
 Further, as $G^1J$ is an open subgroup of $G$, 
 we have a $G^1J$-equivariant injection
 $\Ind_{G^1}^{G^1J} H'_{\mathrm{LT}}\hookrightarrow \Ind_{G^1J}^{G}(\Ind_{G^1}^{G^1J} H'_{\mathrm{LT}})=\Ind_{G^1}^{G} H'_{\mathrm{LT}}$ by Lemma \ref{lem:induction-easy} (i).
 Together with Proposition \ref{prop:NALT-cInd},
 we obtain an $F^\times\GL_n(\mathcal{O}_F)\times W_{F_n}$-equivariant
 injection
 \[
  \overline{\pi}_\theta\boxtimes (\theta\delta^{n-1})^{-1}(\tfrac{1-n}{2})\hookrightarrow 
 \Hom_{\mathcal{O}_D^\times/\varpi^\Z}(\overline{\rho}_{\theta^{-1}},H_{\mathrm{LT}}),
 \]
 which induces a non-zero $J$-equivariant map
 $\overline{\pi}_\theta\boxtimes\overline{\rho}_{\theta^{-1}}\boxtimes (\theta\delta^{n-1})^{-1}(\tfrac{1-n}{2})\to H_{\mathrm{LT}}$.
 By the Frobenius reciprocity, this corresponds to a non-zero $G$-equivariant map
 \[
  \pi_\theta\boxtimes \rho_{\theta^{-1}}\boxtimes \Ind_{W_{F_n}}^{W_F}(\theta\delta^{n-1})^{-1}(\tfrac{1-n}{2})
 \to H_{\mathrm{LT}}.
 \]
 Since $\pi_\theta$ is supercuspidal and its central character is trivial on $\varpi^\Z$,
 Theorem \ref{thm:NALT} tells us that there exists a non-zero map 
 $\rho_{\theta^{-1}}\boxtimes \Ind_{W_{F_n}}^{W_F}(\theta\delta^{n-1})^{-1}\to \JL(\pi_\theta)^\vee\boxtimes \rec_F(\pi_\theta)^\vee$.
 As $\rho_{\theta^{-1}}$ and $\JL(\pi_\theta)^\vee$ are irreducible, 
 we have $\rho_{\theta^{-1}}=\JL(\pi_\theta)^\vee$ and $\JL(\pi_\theta)=\rho_{\theta^{-1}}^\vee=\rho_\theta$.
 As $\rec_F(\pi_\theta)^\vee$ is irreducible and 
 $\dim \Ind_{W_{F_n}}^{W_F}(\theta\delta^{n-1})^{-1}=n=\dim \rec_F(\pi_\theta)^\vee$, we conclude that
 $\rec_F(\pi_\theta)^\vee=\Ind_{W_{F_n}}^{W_F}(\theta\delta^{n-1})^{-1}$,
 hence $\rec_F(\pi_\theta)=\Ind_{W_{F_n}}^{W_F}(\theta\delta^{n-1})$.
\end{prf}

\section{Example II: simple supercuspidal representations}\label{sec:ssc}
In this section, we determine $\rec_F(\pi)$ and $\JL(\pi)$ for a simple supercuspidal representation
$\pi$ of $\GL_n(F)$. In the following we fix a uniformizer $\varpi$ of $F$ and
a non-trivial additive character $\psi_0\colon \F_p\to\C^\times$. We write $\psi$ for the composite
$\F_q\xrightarrow{\Tr_{\F_q/\F_p}}\F_p\xrightarrow{\psi_0}\C^\times$.

We briefly recall the notion of simple supercuspidal representations of $\GL_n(F)$ and $D^\times$.
See \cite[\S 4.1]{parity-csd} for more detail.
Let $\Iw$ be the standard Iwahori subgroup of $\GL_n(F)$ and $\Iw_+$ the pro-$p$ radical of $\Iw$.
The character $\psi$ determines a character of $\Iw_+$:
\[
 \psi\colon \Iw_+\to \C^\times; (a_{ij})\mapsto \psi(\overline{a_{12}}+\overline{a_{23}}+\cdots+\overline{a_{n-1,n}}+\overline{\varpi^{-1}a_{n1}}).
\]
Here the image of $a\in \mathcal{O}_F$ under the map $\mathcal{O}_F\twoheadrightarrow\F_q$ is denoted by
$\overline{a}$.
Put
\[
 \varphi=
 \begin{pmatrix}
 0& 1& 0& \cdots& 0\\
 0& 0& 1& \cdots& 0\\
 \vdots& \vdots& \vdots& & \vdots\\
 0& 0& 0& \cdots& 1\\
 \varpi& 0& 0& \cdots& 0
 \end{pmatrix}\in \GL_n(F).
\]
It normalizes $(\Iw_+,\psi)$.

\begin{defn}\label{defn:ssc-GL}
 For a character $\chi\colon\F_q^\times\to \C^\times$ and $c\in \C^\times$, 
 let $\Lambda_{\chi,c}\colon \mathcal{O}_F^\times\varphi^\Z\Iw_+\to \C^\times$ be the character
 defined as follows:
 \[
 \Lambda_{\chi,c}(a)=\chi(\overline{a})\ (a\in\mathcal{O}_F^\times),\quad \Lambda_{\chi,c}(\varphi)=c,\quad
 \Lambda_{\chi,c}\vert_{\Iw_+}=\psi.
 \]
 We put $\pi_{\chi,c}=\cInd_{\mathcal{O}_F^\times\varphi^\Z\Iw_+}^{\GL_n(F)}\Lambda_{\chi,c}$,
 which is known to be an irreducible supercuspidal representation of $\GL_n(F)$.
\end{defn}

Next we consider the group $D^\times$. Recall that $F_n$ denotes the degree $n$ unramified extension of $F$.
We write $\tau\in\Gal(F_n/F)$ for the arithmetic Frobenius lift.
Since the invariant of $D$ is $1/n$, there exists an isomorphism
$D\cong F_n[\Pi]$, where $\Pi^n=\varpi$ and $\Pi a=\tau(a)\Pi$ for $a\in F_n$.
In the following, we fix this isomorphism and identify them.

\begin{defn}\label{defn:ssc-D}
 For a character $\chi\colon\F_q^\times\to \C^\times$ and $c\in \C^\times$, 
 let $\Lambda^D_{\chi,c}\colon \mathcal{O}_F^\times\Pi^\Z(1+\Pi\mathcal{O}_D)\to \C^\times$ be the character
 defined as follows:
 \[
 \Lambda^D_{\chi,c}(a)=\chi(\overline{a})\ (a\in\mathcal{O}_F^\times),\quad \Lambda^D_{\chi,c}(\Pi)=c,\quad
 \Lambda^D_{\chi,c}(1+\Pi d)=\psi(\Tr_{\F_{q^n}/\F_q}(\overline{d})).
 \]
 We put $\rho_{\chi,c}=\cInd_{\mathcal{O}_F^\times\Pi^\Z(1+\Pi\mathcal{O}_D)}^{D^\times}\Lambda^D_{\chi,c}$,
 which is known to be an irreducible smooth representation of $D^\times$.
 Note that the isomorphism class of $\rho_{\chi,c}$ is independent of 
 the choice of an isomorphism $D\cong F_n[\Pi]$, as every automorphism of $D$ is inner.
\end{defn}

\begin{rem}
 \begin{enumerate}
  \item As in \cite[\S 4.1]{parity-csd}, simple supercuspidal representations $\{\pi_{\zeta,\chi,c}\}$
	of $\GL_n(F)$ are in fact parameterized by
	$\zeta\in \F_q^\times$ and $(\chi,c)$ as above. The representation $\pi_{\chi,c}$
	in Definition \ref{defn:ssc-GL} corresponds to $\pi_{1,\chi,c}$.
	On the other hand, $\pi_{\zeta,\chi,c}$ equals $\pi_{\chi,c}$
	for the uniformizer $\widetilde{\zeta}\varpi$, where $\widetilde{\zeta}\in \mathcal{O}_F^\times$
	is the Teichm\"uller lift of $\zeta$.
	Therefore, if we vary the fixed uniformizer $\varpi$, the representations of the form $\pi_{\chi,c}$
	cover all simple supercuspidal representations.
	Similar remark holds for $\rho_{\chi,c}$.
  \item Let $\xi\colon F^\times\to \C^\times$ be an unramified character.
	Then, we can easily see that
	$\pi_{\chi,c}\otimes (\xi\circ\det)\cong\pi_{\chi,\xi(\varpi)c}$ and
	$\rho_{\chi,c}\otimes (\xi\circ\Nrd)\cong\rho_{\chi,\xi(\varpi)c}$.
	Since $\rec_F$ and $\JL$ are compatible with character twist,
	to determine $\rec_F(\pi_{\chi,c})$ and $\JL(\pi_{\chi,c})$, we may assume that $c=1$.
	In the following we put $\Lambda_\chi=\Lambda_{\chi,1}$, $\pi_\chi=\pi_{\chi,1}$,
	$\Lambda^D_\chi=\Lambda^D_{\chi,(-1)^{n-1}}$ and $\rho_\chi=\rho_{\chi,(-1)^{n-1}}$.
 \end{enumerate}
\end{rem}

We write $n=n'p^e$ with $p\nmid n'$.
Take an $n$th root $\varpi_L$ of $\varpi$ and put $L=F(\varpi_L)$,
$\varpi_E=\varpi_L^{p^e}$ and $E=F(\varpi_E)$.
For $x\in E$, we write
\[
 x=\varpi_E^{v_E(x)}\zeta_x(1+u_x\varpi_E),
\]
where $\zeta_x\in\mu_{q-1}(E)=\{a\in E\mid a^{q-1}=1\}\subset \mathcal{O}_F$ and $u_x\in\mathcal{O}_E$.
The goal of this section is as follows:

\begin{thm}\label{thm:LLC-JLJC-ssc}
 Let $\chi$ be a character of $\F_q^\times$.
 \begin{enumerate}
  \item We have $\JL(\pi_\chi)=\rho_\chi$.
  \item Assume that $p\nmid n$. 
	We define two characters $\xi_\chi,\mu\colon E^\times\to \C^\times$ by
	\begin{align*}
	 \xi_\chi(x)=\chi(\overline{\zeta}_x)\psi(\overline{u}_x)^n,\quad
	 \mu(x)=\lambda^{v_E(x)}\delta(\overline{\zeta}_x)^{n-1},
	\end{align*}
	where
	\[
	 \lambda=\begin{cases}
		  q^{-\frac{n-1}{2}}\sum_{y_1,\ldots,y_{n-1}\in \F_q}\psi\bigl(\sum_{1\le i\le j\le n-1}y_iy_j\bigr)& p\neq 2,\\
		  (\tfrac{q}{n})& p=2,
		 \end{cases}
	\]
	and $\delta\colon \F_q^\times\to \C^\times$ is the quadratic character.
	Then we have
	\[
	 \rec_F(\pi_\chi)=\Ind_{W_E}^{W_F} \mu^{-1}\xi_\chi.
	\]
  \item Assume that $p\mid n$. 
	We define two characters $\nu_\chi,\phi\colon E^\times\to \C$ by
	\[
	\nu_\chi(x)=\chi(\overline{\zeta}{}_x^{p^{-e}}),\quad \phi(x)=(-1)^{(n-1)v_E(x)}.
	\]
	Put $f=[\F_q:\F_p]$ and $m=\gcd(e,f)$.
	Let $\ITw$ be the affine algebraic variety over $\overline{\F}_q$ defined by the equation
	\[
	z^{p^m}-z=y^{p^e+1}-\frac{1}{n'}\sum_{1\le i\le j\le n-2}y_iy_j.
	\]
	As defined in \cite[(3.17), (3.18)]{Imai-Tsushima-wild}, the group $W_E$ acts on $\ITw$.
	On the other hand, the group $\F_{p^m}$ acts on $\ITw$ by
	$a\colon z\mapsto z+a$, $y\mapsto y$, $y_i\mapsto y_i$ $(a\in \F_{p^m})$. These two actions commute.
	Let $\psi'$ denote the additive character $\psi_0^{n'}\circ\Tr_{\F_{p^m}/\F_p}$ of $\F_{p^m}$.
	We write $H^{n-1}_c(\ITw,\overline{\Q}_\ell)_{\psi'}$ for the subspace of 
	$H^{n-1}_c(\ITw,\overline{\Q}_\ell)$ on which $\F_{p^m}$ acts by $\psi'$.
	It is a smooth representation of $W_E$.
	Put $\tau=H^{n-1}_c(\ITw,\overline{\Q}_\ell)_{\psi'}(\frac{n-1}{2})$.

	Then we have
	\[
	 \rec_F(\pi_\chi)=\Ind_{W_E}^{W_F}(\tau\otimes \nu_\chi\otimes \phi).
	\]
 \end{enumerate}
\end{thm}

\begin{rem}\label{rem:history-ssc}
 \begin{enumerate}
  \item Theorem \ref{thm:LLC-JLJC-ssc} (i) was proved in \cite{Imai-Tsushima-JL}
	by explicit computation of characters.
  \item Theorem \ref{thm:LLC-JLJC-ssc} (ii), (iii) were obtained in \cite{MR2148193}, \cite{2015arXiv150902960I},
	respectively (for (ii), see also \cite[Theorem 2.1 (1)]{MR2148193} and \cite[\S 5]{Imai-Tsushima-tame}).
 \end{enumerate}
\end{rem}

The proof of Theorem \ref{thm:LLC-JLJC-ssc} is divided into the cases $p\nmid n$ and $p\mid n$.
In this paper we only treat the case $p\mid n$; the other is similar but easier. 

\begin{defn}\label{defn:ssc-J}
 We put 
 \[
  J=\mathcal{O}_F^\times\varphi^\Z\Iw_+\times (\mathcal{O}_F^\times\Pi^\Z(1+\Pi\mathcal{O}_D)/\varpi^\Z)\times W_E
 \subset G,\quad J^1=G^1\cap J.
 \]
 Let $J^1_{\mathrm{s}}$ be the subgroup of $J^1$ generated by the following elements:
 \begin{itemize}
  \item $(g,d,1)$ for $g\in \Iw_+$ and $d\in 1+\Pi\mathcal{O}_D$ with $\det g=\Nrd d$,
  \item $(a,a,1)$ for $a\in \mathcal{O}_F^\times$,
  \item $(\varphi,\Pi,1)$, and
  \item $(g_\sigma,\Pi^{-n_\sigma},\sigma)$ for $\sigma\in W_E$,
	where $g_\sigma\in \mathcal{O}_F^\times\Iw_+$ is the element in \cite[(3.14)]{Imai-Tsushima-wild},
	and $n_\sigma=v_E(\Art^{-1}_E(\sigma))$.
 \end{itemize}
\end{defn} 

\begin{lem}\label{lem:J-property}
 \begin{enumerate}
  \item The group $J^1_{\mathrm{s}}$ is an open normal subgroup of $J^1$.
  \item Put $n'_1=\gcd(n',q-1)$ and $\mu=\mu_{n'_1}(\mathcal{O}_F)$. We regard $\mu$ as a subgroup of 
	$D^\times/\varpi^\Z$. Then, we have
	$J^1=\bigcup_{\zeta\in \mu}\zeta J^1_{\mathrm{s}}$ (later we will also prove that the union
	in the right hand side is disjoint).
  \item The composite $J^1_{\mathrm{s}}\hookrightarrow J\twoheadrightarrow \mathcal{O}_F^\times\varphi^\Z\Iw_+\times W_E$ is surjective.
 \end{enumerate}
\end{lem}

\begin{prf}
 Let $I$ be the open subgroup of $W_E$ consisting of $\sigma$ such that $\Art_E^{-1}(\sigma)\in 1+\mathfrak{p}_E$.
 If $\sigma\in I$, then $g_\sigma$ belongs to $\Iw_+$ and $n_\sigma=0$.
 Therefore, we can easily check that $J^1_{\mathrm{s}}$ contains 
 $J^1\cap (\Iw_+\times (1+\Pi\mathcal{O}_D)\times I)$. Hence $J^1_{\mathrm{s}}$ is open in $J^1$.

 We prove (ii). Take $(g,d,\sigma)\in J^1$. If we put $g'=gg_\sigma^{-1}$ and $d'=d\Pi^{n_\sigma}$,
 then $(g,d,\sigma)\in (g',d',1)J^1_{\mathrm{s}}$. Since $\det g'=\Nrd d'$, there exists an integer $m$
 such that $\varphi^mg'\in\mathcal{O}_F^\times\Iw_+$ and $\Pi^md'\in \mathcal{O}_F^\times(1+\Pi\mathcal{O}_D)$.
 Write $\varphi^mg'=a g''$ and $\Pi^md'=a'd''$, where $a,a'\in \mu_{q-1}(\mathcal{O}_F)$,
 $g''\in \Iw_+$ and $d''\in 1+\Pi\mathcal{O}_D$.
 Since $\det(\varphi^mg')=\Nrd(\Pi^md')$ and $\det g'',\Nrd d''\in 1+\mathfrak{p}_F$,
 we have $a^n=a'^n$ and $\det g''=\Nrd d''$. Put $\zeta=a'/a$, which is an element of $\mu$.
 Now we have 
 \[
  (g',d',1)=(1,\zeta,1)(a,a,1)(\varphi,\Pi,1)^{-m}(g'',d'',1)\in \zeta J^1_{\mathrm{s}}
 \]
 and $(g,d,\sigma)\in \zeta J^1_{\mathrm{s}}$. This concludes (ii).
 As $\mu$ is contained in the center of $J^1$, $J^1_{\mathrm{s}}$ is normal in $J^1$, hence (i).

 Finally consider (iii). Take an element $(g,\sigma)$ of $\mathcal{O}^\times_F\varphi^\Z\Iw_+\times W_E$.
 As $gg_\sigma^{-1}\in \mathcal{O}^\times_F\varphi^\Z\Iw_+$, there exist
 $a\in \mathcal{O}^\times_F$, $m\in \Z$ and $g'\in \Iw_+$ such that
 $gg_\sigma^{-1}=a\varphi^mg'$.
 Since $\Nrd\colon 1+\Pi\mathcal{O}_D\to 1+\mathfrak{p}_F$ is surjective (see \cite[Lemma 5]{MR0262250}),
 we can take $d\in 1+\Pi\mathcal{O}_D$ such that $\Nrd d=\det g'$.
 Then, the element
 \[
  (a,a,1)(\varphi,\Pi,1)^m(g',d,1)(g_\sigma,\Pi^{-n_\sigma},\sigma)\in J^1_{\mathrm{s}}
 \]
 is mapped to $(g,\sigma)$ under $J^1_{\mathrm{s}}\hookrightarrow J\twoheadrightarrow \mathcal{O}_F^\times\varphi^\Z\Iw_+\times W_E$.
\end{prf}

Now we use results in \cite{Imai-Tsushima-wild}. We may assume that $\varphi_{D,1}$ in 
\cite[\S 2.1]{Imai-Tsushima-wild} is equal to $\Pi$.
Let $U\subset M^{(0)}_{\infty,C}$ be the affinoid and $\mathcal{X}$ the formal model of $U$
constructed in \cite[\S 2]{Imai-Tsushima-wild} (in \cite{Imai-Tsushima-wild}, $U$ and $\mathcal{X}$
are denoted by $\mathcal{X}_1$ and $\mathfrak{X}_1$, respectively).
By construction, $U$ is a rational subset of $t(\mathcal{M}^{(0)}_{\infty,\mathcal{O}_C})$.

\begin{prop}\label{prop:affinoid-ssc}
 \begin{enumerate}
  \item The action of $J^1_{\mathrm{s}}\subset G^1$ on $M^{(0)}_{\infty,C}$ stabilizes $U$ and extends to $\mathcal{X}$.
  \item The formal scheme $\mathcal{X}$ is flat over $k^\circ$ and its reduction $\mathcal{X}_s$ is
	isomorphic to the perfection of $\ITw$.
  \item The induced action of $J^1_{\mathrm{s}}$ on $\mathcal{X}_s$ comes from an action on $\ITw$,
	which is described as follows:
	\begin{itemize}
	 \item $(g,d,1)$ for $g=(g_{ij})\in \Iw_+$ and $d=1+\Pi d'\in 1+\Pi\mathcal{O}_D$ with $\det g=\Nrd d$
	       acts on $\ITw$ as the element
	       \[
	       \frac{1}{n'}\Tr_{\F_q/\F_{p^m}}\bigl(\Tr_{\F_{q^n}/\F_q}(\overline{d'})-(\overline{g_{12}}+\cdots+\overline{g_{n-1,n}}+\overline{\varpi^{-1}g_{n1}})\bigr)\in \F_{p^m},
	       \]
	 \item $(a,a,1)$ for $a\in \mathcal{O}_F^\times$ acts on $\ITw$ trivially,
	 \item $(\varphi,\Pi,1)$ acts on $\ITw$ as in \cite[(3.2)]{Imai-Tsushima-wild}, and
	 \item $(g_\sigma,\Pi^{-n_\sigma},\sigma)$ for $\sigma\in W_E$ acts on $\ITw$
	       as mentioned in Theorem \ref{thm:LLC-JLJC-ssc} (iii).
	\end{itemize}
  \item For $\zeta\in \mu\setminus \{1\}$, we have $U\zeta\cap U=\varnothing$.
  \item We have $J^1=\coprod_{\zeta\in \mu}\zeta J^1_{\mathrm{s}}$
	and $J^1_{\mathrm{s}}\cap (\mathcal{O}_F^\times\Pi^\Z(1+\Pi\mathcal{O}_D)/\varpi^\Z)=(1+\Pi\mathcal{O}_D)^1$, where $(1+\Pi\mathcal{O}_D)^1$ denotes the subgroup of $1+\Pi\mathcal{O}_D$ consisting of elements
	with reduced norm $1$.
 \end{enumerate}
\end{prop}

\begin{prf}
 The assertion (i) is obtained in \cite[\S 3]{Imai-Tsushima-wild}.
 Consider (ii).
 The reduction $\mathcal{X}_s$ is computed in \cite[Theorem 2.5]{Imai-Tsushima-wild}.
 In particular, $\dim \mathcal{X}_s=\dim \ITw=n-1$.
 By Corollary \ref{cor:flatness-dimension} and the presentation of $\mathcal{X}$ in \cite[\S 2]{Imai-Tsushima-wild}, 
 we conclude that $\mathcal{X}$ is flat over $k^\circ$.
 The claim (iii) is also included in \cite[\S 3]{Imai-Tsushima-wild}.

 We prove (iv). We use the notation in \cite{Imai-Tsushima-wild}. 
 Recall that $U$ is defined by the inequalities \cite[(2.5)]{Imai-Tsushima-wild}:
 \begin{align*}
  &v\Bigl(\frac{\boldsymbol{x}_i}{\boldsymbol{x}_{i+1}}-\Bigl(\frac{\boldsymbol{x}_{n-1}}{\boldsymbol{x}_n}\Bigr)^{q^{n-1-i}}\Bigr)\ge \frac{1}{2nq^i}\quad \text{for $1\le i\le n-2$,}\\
  &v(\boldsymbol{x}_i-1)\ge \frac{1}{nq^{n-1}(p^e+1)}\quad \text{for $n-1\le i\le n$.}
 \end{align*}
 By \cite[(1.7)]{Imai-Tsushima-wild}, $\zeta$ acts by $\boldsymbol{x}_i\mapsto \zeta^{-1}\boldsymbol{x}_i$.
 Since $n_1'$ is prime to $p$ and $\zeta\neq 1$, we have $\zeta^{-1}-1\in \mathcal{O}_F^\times$.
 Therefore, if $v(\boldsymbol{x}_n-1)>0$, we have $v(\boldsymbol{x}_n)=0$ and
 \[
  v(\zeta^{-1}\boldsymbol{x}_n-1)=v((\zeta^{-1}-1)\boldsymbol{x}_n+\boldsymbol{x}_n-1)=0.
 \]
 Hence we have $U\cap U\zeta=\varnothing$, as desired.

 Finally we prove (v). By (i) and (iv), we have $\mu\cap J^1_{\mathrm{s}}=\{1\}$.
 Hence the cosets $\{\zeta J^1_{\mathrm{s}}\}_{\zeta\in\mu}$ are disjoint to each other.
 Together with Lemma \ref{lem:J-property} (ii), we obtain $J^1=\coprod_{\zeta\in \mu}\zeta J^1_{\mathrm{s}}$.
 On the other hand, clearly we have
 \begin{align*}
  (1+\Pi\mathcal{O}_D)^1&\subset J^1_{\mathrm{s}}\cap (\mathcal{O}_F^\times\Pi^\Z(1+\Pi\mathcal{O}_D)/\varpi^\Z)\\
  &\subset J^1\cap (\mathcal{O}_F^\times\Pi^\Z(1+\Pi\mathcal{O}_D)/\varpi^\Z)=\mu(1+\Pi\mathcal{O}_D)^1.
 \end{align*}
 Since $\mu\cap J^1_{\mathrm{s}}=\{1\}$, we conclude that $(1+\Pi\mathcal{O}_D)^1=J^1_{\mathrm{s}}\cap (\mathcal{O}_F^\times\Pi^\Z(1+\Pi\mathcal{O}_D)/\varpi^\Z)$.
\end{prf}

\begin{cor}\label{cor:coh-ITw}
 The subspace $H^{n-1}_c(\ITw,\overline{\Q}_\ell)_{\psi'}$ of $H^{n-1}_c(\ITw,\overline{\Q}_\ell)$
 in Theorem \ref{thm:LLC-JLJC-ssc} (iii) is $J^1_{\mathrm{s}}$-stable. Moreover, we have
 an isomorphism 
 \[
  \Hom_{J^1_{\mathrm{s}}\cap (\mathcal{O}_F^\times\Pi^\Z(1+\Pi\mathcal{O}_D)/\varpi^\Z)}\bigl(\Lambda^D_\chi,H^{n-1}_c(\ITw,\overline{\Q}_\ell)_{\psi'}\bigr) \cong \Lambda_\chi^{-1}\boxtimes  (\tau\otimes \nu_\chi\otimes \phi)(\tfrac{1-n}{2})
 \]
 of representations of $\mathcal{O}_F^\times\varphi^\Z\Iw_+\times W_E$.
 Note that Lemma \ref{lem:J-property} (iii) enables us to define an action of
 $\mathcal{O}_F^\times\varphi^\Z\Iw_+\times W_E$ on the left hand side 
 as in Lemma \ref{lem:representation-theory} (i).
\end{cor}

\begin{prf}
 Since the action of $\F_{p^m}$ on $H^{n-1}_c(\ITw,\overline{\Q}_\ell)$ commutes with the action of 
 $J^1_{\mathrm{s}}$, the subspace $H^{n-1}_c(\ITw,\overline{\Q}_\ell)_{\psi'}\subset H^{n-1}_c(\ITw,\overline{\Q}_\ell)$ is $J^1_{\mathrm{s}}$-stable.
 
 Recall that in Proposition \ref{prop:affinoid-ssc} (v) we proved that 
 $J^1_{\mathrm{s}}\cap (\mathcal{O}_F^\times\Pi^\Z(1+\Pi\mathcal{O}_D)/\varpi^\Z)=(1+\Pi\mathcal{O}_D)^1$.
 For $1+\Pi d\in (1+\Pi\mathcal{O}_D)^1$, we have
 $\Lambda^D_\chi(1+\Pi d)=\psi(\Tr_{\F_{q^n}/\F_q}(\overline{d}))$.
 On the other hand, by Proposition \ref{prop:affinoid-ssc} (iii),
 $1+\Pi d$ acts on $H^{n-1}_c(\ITw,\overline{\Q}_\ell)_{\psi'}$ by
 $\psi'(n'^{-1}\Tr_{\F_{q^n}/\F_{p^m}}(\overline{d}))=\psi(\Tr_{\F_{q^n}/\F_q}(\overline{d}))$.
 Therefore, as a vector space we have
 \begin{equation*}
  \Hom_{J^1_{\mathrm{s}}\cap (\mathcal{O}_F^\times\Pi^\Z(1+\Pi\mathcal{O}_D)/\varpi^\Z)}\bigl(\Lambda^D_\chi,H^{n-1}_c(\ITw,\overline{\Q}_\ell)_{\psi'}\bigr)=H^{n-1}_c(\ITw,\overline{\Q}_\ell)_{\psi'}.\tag{$*$}
 \end{equation*}
 Under this identification, the action of $\mathcal{O}_F^\times\varphi^\Z\Iw_+\times W_E$
 on the left hand side is described as follows:
 \begin{enumerate}
  \item[(a)] $(a,1)$ for $a\in \mathcal{O}^\times_F$ acts by the scalar
	     $\chi(\overline{a})^{-1}=\Lambda_\chi^{-1}(a)$,
  \item[(b)] $(\varphi,1)$ acts trivially,
  \item[(c)] $(g,1)$ for $g\in \Iw_+$ acts by the scalar $\psi(g)^{-1}=\Lambda_\chi^{-1}(g)$ and
  \item[(d)] $(g_\sigma,\sigma)$ for $\sigma\in W_E$ acts on $H^{n-1}_c(\ITw,\overline{\Q}_\ell)_{\psi'}$
	as $(-1)^{(n-1)n_\sigma}\sigma$ on $H^{n-1}_c(\ITw,\overline{\Q}_\ell)_{\psi'}$.
 \end{enumerate}
 These are consequences of Proposition \ref{prop:affinoid-ssc} (iii). 
 For (b), note that $(\varphi,\Pi)$ acts on $H^{n-1}_c(\ITw,\overline{\Q}_\ell)$ by the scalar $(-1)^{n-1}$
 (see \cite[Proposition 3.2.2, Corollary 4.5]{Imai-Tsushima-wild} and \cite[Proposition 4.2.3]{Imai-Tsushima-tame}).

 For $\sigma\in W_E$, put $u_\sigma=\varpi_E^{-n_\sigma}\Art_E^{-1}(\sigma)$ and take 
 $b_\sigma\in \mu_{q-1}(\mathcal{O}_F)$ such that $\overline{b}{}_\sigma^{p^e}=\overline{u}_\sigma$.
 Recall that $g_\sigma=b_\sigma\mathrm{diag}(1,\ldots,1,c_\sigma)$ for some $c_\sigma\in 1+\mathfrak{p}_F$
 (see \cite[(3.14)]{Imai-Tsushima-wild}).
 Therefore, by (a) and (c), $(g_\sigma,1)$ acts by the scalar 
 $\chi(\overline{b}_\sigma)^{-1}=\nu_\chi(\Art_E^{-1}(\sigma))^{-1}$.
 Hence the action of $(1,\sigma)$ on the left hand side of $(*)$ is $\nu_\chi(\Art_E^{-1}(\sigma))\phi(\Art_E^{-1}(\sigma))\sigma$. Now we conclude that
 \[
    \Hom_{J^1_{\mathrm{s}}\cap (\mathcal{O}_F^\times\Pi^\Z(1+\Pi\mathcal{O}_D)/\varpi^\Z)}\bigl(\Lambda^D_\chi,H^{n-1}_c(\ITw,\overline{\Q}_\ell)_{\psi'}\bigr)\cong \Lambda^{-1}_\chi\boxtimes (\tau\otimes \nu_\chi\otimes \phi)(\tfrac{1-n}{2}).
 \]
\end{prf}

\begin{prop}\label{prop:Artin-Schreier}
 Let $\mathcal{L}_{\psi'}$ denote the rank $1$ sheaf over $\A^1_{\overline{\F}_q}=\Spec \overline{\F}_q[t]$
 defined by the $\F_{p^m}$-torsor $z^{p^m}-z=t$ and the additive character 
 $\psi'\colon \F_{p^m}\to \C^\times$,
 and $\mathcal{L}_{\psi'}^{\ITw}$ the pull-back of $\mathcal{L}_{\psi'}$ under
 \[
  \A^{n-1}_{\overline{\F}_q}\to \A^1_{\overline{\F}_q};\quad
 (y,y_1,\ldots,y_{n-2})\mapsto y^{p^e+1}-\frac{1}{n'}\sum_{1\le i\le j\le n-2}y_iy_j.
 \]
 Then, $\dim H^{n-1}_c(\A^{n-1}_{\overline{\F}_q},\mathcal{L}_{\psi'}^{\ITw})=p^e$ and
 the map 
 \[
  H^{n-1}_c(\A^{n-1}_{\overline{\F}_q},\mathcal{L}_{\psi'}^{\ITw})\to H^{n-1}(\A^{n-1}_{\overline{\F}_q},\mathcal{L}_{\psi'}^{\ITw})
 \]
 is an isomorphism.

 In particular, $\dim H^{n-1}_c(\ITw,\overline{\Q}_\ell)_{\psi'}=p^e$ and the composite
 \[
  H^{n-1}_c(\ITw,\overline{\Q}_\ell)_{\psi'}\hookrightarrow H^{n-1}_c(\ITw,\overline{\Q}_\ell)\to H^{n-1}(\ITw,\overline{\Q}_\ell)
 \]
 is injective.
\end{prop}

\begin{prf}
 We denote by $\mathcal{L}^1_{\psi'}$ (resp.\ $\mathcal{L}^2_{\psi'}$) the pull-back of $\mathcal{L}_{\psi'}$ under
 $\A^1_{\overline{\F}_q}\to \A^1_{\overline{\F}_q}$; $y\mapsto y^{p^e+1}$ (resp.\
 $\A^{n-2}_{\overline{\F}_q}\to \A^1_{\overline{\F}_q}$;
 $(y_1,\ldots,y_{n-2})\mapsto -n'^{-1}\sum_{1\le i\le j\le n-2}y_iy_j$).
 Then we have $\mathcal{L}^{\ITw}_{\psi'}=\mathcal{L}^1_{\psi'}\boxtimes \mathcal{L}^2_{\psi'}$.

 We use results in Appendix \ref{sec:Artin-Schreier}.
 By Remark \ref{rem:m=1} and Lemma \ref{lem:AS-A^1}, we have
 \begin{itemize}
  \item $\dim H^1_c(\A^1_{\overline{\F}_q},\mathcal{L}^1_{\psi'})=p^e$,
  \item $H^i_c(\A^1_{\overline{\F}_q},\mathcal{L}^1_{\psi'})=H^i(\A^1_{\overline{\F}_q},\mathcal{L}^1_{\psi'})=0$
	for $i\neq 1$, and
  \item $H^1_c(\A^1_{\overline{\F}_q},\mathcal{L}^1_{\psi'})\to H^1(\A^1_{\overline{\F}_q},\mathcal{L}^1_{\psi'})$
	is an isomorphism.
 \end{itemize}
 By Remark \ref{rem:m=1}, Lemma \ref{lem:AS-quadratic} and Example \ref{exa:non-degenerate}, we have
 \begin{itemize}
  \item $\dim H^{n-2}_c(\A^{n-2}_{\overline{\F}_q},\mathcal{L}^2_{\psi'})=1$ and
  \item $H^{n-2}_c(\A^{n-2}_{\overline{\F}_q},\mathcal{L}^2_{\psi'})\to H^{n-2}(\A^{n-2}_{\overline{\F}_q},\mathcal{L}^2_{\psi'})$ is an isomorphism.
 \end{itemize}
 Hence we conclude the proposition by the K\"unneth formula.
\end{prf}

\begin{cor}\label{cor:sp-inj-ssc}
 We have a $J^1$-equivariant injection
 $\Ind_{J^1_{\mathrm{s}}}^{J^1} H^{n-1}_c(\ITw,\overline{\Q}_\ell)_{\psi'}\to H'_{\mathrm{LT}}$.
\end{cor}

\begin{prf}
 By Proposition \ref{prop:affinoid-ssc} (i), (ii), the pair $(U,\mathcal{X})$ satisfies
 the condition in Proposition \ref{prop:LT-tower-affinoid}
 (we may take $K_U=J^1_{\mathrm{s}}\cap \SL_n(\mathcal{O}_F)$).
 It is easy to see that $\ITw$ is purely $n-1$-dimensional and smooth.
 By Proposition \ref{prop:affinoid-ssc} (i), (iv) and (v), for $g\in J^1\setminus J^1_{\mathrm{s}}$ we have
 $U\cap Ug=\varnothing$. 
 Hence, by Corollary \ref{cor:LT-specialization} (iii) and Proposition \ref{prop:Artin-Schreier},
 we have $J^1$-equivariant homomorphisms
 \[
  \Ind_{J^1_{\mathrm{s}}}^{J^1} H^{n-1}_c(\ITw,\overline{\Q}_\ell)_{\psi'}\hookrightarrow \Ind_{J^1_{\mathrm{s}}}^{J^1} H^{n-1}_c(\ITw,\overline{\Q}_\ell)\to H'_{\mathrm{LT}},
 \]
 whose composite is injective.
\end{prf}

\begin{prf}[of Theorem \ref{thm:LLC-JLJC-ssc}]
 By Corollary \ref{cor:sp-inj-ssc}, we have an injection
 \begin{align*}
  &\Hom_{J^1\cap (\mathcal{O}_F^\times\Pi^\Z(1+\Pi\mathcal{O}_D)/\varpi^\Z)}\bigl(\Lambda_\chi^D,\Ind_{J^1_{\mathrm{s}}}^{J^1}H^{n-1}_c(\ITw,\overline{\Q}_\ell)_{\psi'}\bigr)\\
  &\qquad\qquad\hookrightarrow 
 \Hom_{J^1\cap (\mathcal{O}_F^\times\Pi^\Z(1+\Pi\mathcal{O}_D)/\varpi^\Z)}(\Lambda_\chi^D,H'_{\mathrm{LT}}),
 \end{align*}
 which is $\mathcal{O}_F^\times\varphi^\Z\Iw_+\times W_E$-equivariant.
 Taking the $\mathcal{O}_F^\times\varphi^\Z\Iw_+\times W_E$-smooth part $(-)^{\mathrm{sm}}$,
 we obtain an $\mathcal{O}_F^\times\varphi^\Z\Iw_+\times W_E$-equivariant injection
 \begin{align*}
  &\Hom_{J^1\cap (\mathcal{O}_F^\times\Pi^\Z(1+\Pi\mathcal{O}_D)/\varpi^\Z)}\bigl(\Lambda_\chi^D,\Ind_{J^1_{\mathrm{s}}}^{J^1}H^{n-1}_c(\ITw,\overline{\Q}_\ell)_{\psi'}\bigr)^{\mathrm{sm}}\\
  &\qquad\qquad\hookrightarrow 
 \Hom_{J^1\cap (\mathcal{O}_F^\times\Pi^\Z(1+\Pi\mathcal{O}_D)/\varpi^\Z)}(\Lambda_\chi^D,H'_{\mathrm{LT}})^{\mathrm{sm}}.
 \end{align*}
 By Lemma \ref{lem:representation-theory} (iii), the left hand side is identified with
 \[
  \Hom_{J^1_{\mathrm{s}}\cap (\mathcal{O}_F^\times\Pi^\Z(1+\Pi\mathcal{O}_D)/\varpi^\Z)}\bigl(\Lambda_\chi^D,H^{n-1}_c(\ITw,\overline{\Q}_\ell)_{\psi'}\bigr)^{\mathrm{sm}}.
 \]
 On the other hand, in the same way as in the proof of Theorem \ref{thm:LLC-JLJC-depth-0},
 we have
 \begin{align*}
  \Hom_{J^1\cap (\mathcal{O}_F^\times\Pi^\Z(1+\Pi\mathcal{O}_D)/\varpi^\Z)}(\Lambda_\chi^D,H'_{\mathrm{LT}})^{\mathrm{sm}}
  &\cong \Hom_{\mathcal{O}_F^\times\Pi^\Z(1+\Pi\mathcal{O}_D)/\varpi^\Z}(\Lambda_\chi^D,\Ind_{J^1}^J H'_{\mathrm{LT}})\\
  &=\Hom_{\mathcal{O}_F^\times\Pi^\Z(1+\Pi\mathcal{O}_D)/\varpi^\Z}(\Lambda_\chi^D,\Ind_{G^1}^{G^1J} H'_{\mathrm{LT}})\\
  &\hookrightarrow \Hom_{\mathcal{O}_F^\times\Pi^\Z(1+\Pi\mathcal{O}_D)/\varpi^\Z}(\Lambda_\chi^D,\Ind_{G^1}^G H'_{\mathrm{LT}})\\
  &=\Hom_{\mathcal{O}_F^\times\Pi^\Z(1+\Pi\mathcal{O}_D)/\varpi^\Z}(\Lambda_\chi^D,H_{\mathrm{LT}}).
 \end{align*}
 Therefore, by Corollary \ref{cor:coh-ITw} we have an $\mathcal{O}_F^\times\varphi^\Z\Iw_+\times W_E$-equivariant
 injection
 \[
  \Lambda_\chi^{-1}\boxtimes  (\tau\otimes \nu_\chi\otimes \phi)(\tfrac{1-n}{2})\hookrightarrow \Hom_{\mathcal{O}_F^\times\Pi^\Z(1+\Pi\mathcal{O}_D)/\varpi^\Z}(\Lambda_\chi^D,H_{\mathrm{LT}}).
 \]
 This induces a non-zero $J$-equivariant map
 \[
  \Lambda_\chi^{-1}\boxtimes \Lambda_\chi^D\boxtimes (\tau\otimes \nu_\chi\otimes \phi)(\tfrac{1-n}{2})\to H_{\mathrm{LT}},
 \]
 which corresponds to a non-zero $G$-equivariant map
 \[
 \pi_\chi^\vee\boxtimes \rho_\chi\boxtimes \Ind_{W_E}^{W_F}(\tau\otimes \nu_\chi\otimes \phi)(\tfrac{1-n}{2})\to H_{\mathrm{LT}}
 \]
 by the Frobenius reciprocity.
 Since $\pi^\vee_\chi$ is supercuspidal and its central character is trivial on $\varpi^\Z$,
 Theorem \ref{thm:NALT} tells us that there exists a non-zero map 
 $\rho_{\chi}\boxtimes \Ind_{W_E}^{W_F}(\tau\otimes \nu_\chi\otimes \phi)\to \JL(\pi_\chi)\boxtimes \rec_F(\pi_\chi)$.
 As $\rho_\chi$ and $\JL(\pi_\chi)$ are irreducible, 
 we have $\rho_\chi=\JL(\pi_\chi)$.
 As $\rec_F(\pi_\chi)$ is irreducible and 
 $\dim \Ind_{W_E}^{W_F}(\tau\otimes \nu_\chi\otimes \phi)=n'\dim \tau=n'p^e=n=\dim \rec_F(\pi_\chi)$
 (see Proposition \ref{prop:Artin-Schreier}), we conclude that
 $\rec_F(\pi_\chi)=\Ind_{W_E}^{W_F}(\tau\otimes \nu_\chi\otimes \phi)$.
\end{prf}

\begin{rem}
 In the proof above, we do not need any information about the representation $\tau$ except its dimension.
 The irreducibility of $\Ind_{W_E}^{W_F}(\tau\otimes \nu_\chi\otimes \phi)$ is also a consequence
 of the proof.
\end{rem}

\appendix

\section{Some ring theory over $k^\circ$}\label{sec:ring-theory}
Let $k$ be a complete non-archimedean field. As in Section \ref{sec:finite-level},
$k^\circ$ denotes the valuation ring of $k$,
$\mathfrak{m}$ the maximal ideal of $k^\circ$, and $\kappa$ the residue field of $k^\circ$.
Choose a non-zero element $\varpi\in\mathfrak{m}$.
In this appendix, we collect some useful facts on algebras over $k^\circ$.

For a $k^\circ$-algebra $A$, we consider the following three properties:
\begin{enumerate}
 \item[(a)] $A$ is $\varpi$-torsion free (in other words, flat over $k^\circ$).
 \item[(b)] $A/\mathfrak{m}A$ is an integral domain.
 \item[(c)] For every non-zero element $f\in A$, there exist $c_f\in k^\circ$ and $f'\in A\setminus \mathfrak{m}A$
	    such that $f=c_ff'$.
\end{enumerate}

\begin{prop}\label{prop:int-dom}
 Assume that $A$ satisfies the properties (a), (b) and (c).
 Then the following hold.
 \begin{enumerate}
  \item The ring $A$ is an integral domain.
  \item For $f\in A\setminus \mathfrak{m}A$, $A/(f)$ is flat over $k^\circ$.
 \end{enumerate}
\end{prop}

\begin{prf}
 We prove (i). First, by the property (b), we have $A\neq 0$.
 Let $f$ and $g$ be non-zero elements in $A$. Take decompositions $f=c_ff'$ and $g=c_gg'$
 as in the property (c). Clearly $c_fc_g\neq 0$ in $k^\circ$.
 By (b), we have $f'g'\in A\setminus \mathfrak{m}A$, and then $f'g'\neq 0$.
 Hence, by (a) we conclude that $fg=(c_fc_g)(f'g')\neq 0$, as desired.

 Next consider (ii). Take $a\in A$ such that $\varpi a\in (f)$ and prove $a\in (f)$.
 Write $\varpi a=fg$ with $g\in A$. We may assume that $a\neq 0$, which implies that $g\neq 0$ by (a).
 Take a decomposition $g=c_gg'$ as in (c). Since $c_g\neq 0$, we have either $c_g/\varpi\in k^\circ$ or
 $\varpi/c_g\in\mathfrak{m}$. If the latter holds, then $fg'=(\varpi/c_g)a$ lies in $\mathfrak{m}A$.
 As $f,g'\in A\setminus \mathfrak{m}A$, this contradicts the property (b).
 Therefore $c_g/\varpi$ lies in $k^\circ$, and $a=(c_g/\varpi)fg'$ lies in $(f)$.
\end{prf}

\begin{cor}\label{cor:flatness}
 Assume that $A$ satisfies the properties (a), (b) and (c).
 Fix an integer $q>1$.
 Let $(f_m)_{m\ge 1}$ and $(g_m)_{m\ge 1}$ be sequences of elements in $A$ satisfying the following:
 \[
  f_{m+1}^q=f_m,\quad g_{m+1}^q=g_m, \quad f_1-g_1\in A\setminus\mathfrak{m}A.
 \]
 Then the $k^\circ$-algebra $A/(f_m-g_m\mid m\ge 1)$ is flat over $k^\circ$.
\end{cor}

\begin{prf}
 Put $I=(f_m-g_m\mid m\ge 1)$ and $I_m=(f_m-g_m)$.
 Since $f_m-g_m=f_{m+1}^q-g_{m+1}^q=(f_{m+1}-g_{m+1})(f_{m+1}^{q-1}+\cdots+g_{m+1}^{q-1})\in I_{m+1}$, 
 we have $I_1\subset I_2\subset \cdots$ and $I=\bigcup_{m\ge 1}I_m$.
 Therefore $A/I\cong \varinjlim_m A/I_m$, and it suffices to prove that $A/I_m$ is flat over $k^\circ$
 for each $m$.

 If $f_m-g_m$ belongs to $\mathfrak{m}A$, $I_m$ is contained in $\mathfrak{m}A$, and then $f_1-g_1\in I_m$
 lies in $\mathfrak{m}A$. Hence the assumption $f_1-g_1\in A\setminus\mathfrak{m}A$ implies that
 $f_m-g_m\in A\setminus\mathfrak{m}A$. Now Proposition \ref{prop:int-dom} (ii) tells us that
 $A/I_m$ is flat over $k^\circ$.
\end{prf}

\begin{prop}\label{prop:abc-completion}
 Let $\widehat{A}=\varprojlim_m A/\varpi^mA$ denote the $\varpi$-adic completion of $A$.
 If $A$ satisfies the property (a) (resp.\ (b), resp.\ (c)), so does $\widehat{A}$.
\end{prop}

\begin{prf}
 Assume that $A$ satisfies the property (a). Then, for every integer $m\ge 1$, the map
 $A/\varpi^mA\xrightarrow{\times\varpi}A/\varpi^{m+1}A$ is an injection. By taking projective limit,
 we obtain the injectivity of $\widehat{A}\xrightarrow{\times\varpi}\widehat{A}$.
 Hence $\widehat{A}$ satisfies (a).
 
 Next assume that $A$ satisfies (b).
 By \cite[Chapter 0, Lemma 7.2.8, Corollary 7.2.9, Proposition 7.2.15]{Fujiwara-Kato}, we have
 $\widehat{A}/\varpi\widehat{A}\cong A/\varpi A$. 
 By taking base change $(-)\otimes_{k^\circ/\varpi k^\circ}k^\circ/\mathfrak{m}$, we have
 $\widehat{A}/\mathfrak{m}\widehat{A}\cong A/\mathfrak{m}A$.
 Therefore $\widehat{A}/\mathfrak{m}\widehat{A}$ is also an integral domain, and $\widehat{A}$ satisfies
 the property (b).

 Finally assume that $A$ satisfies (c). We write $i$ for the natural homomorphism $A\to \widehat{A}$.
 Take an arbitrary non-zero element $a\in \widehat{A}$.
 By \cite[Chapter 0, Proposition 7.2.15]{Fujiwara-Kato}, $\widehat{A}$ is $\varpi$-adically complete, hence
 $\varpi$-adically separated. Thus there exists an integer $m\ge 0$ such that $a\notin \varpi^m\widehat{A}$.
 Since $\widehat{A}/\varpi^m\widehat{A}\cong A/\varpi^mA$, we can find $f\in A$ such that
 $a-i(f)\in \varpi^m\widehat{A}$.
 Since $a\notin \varpi^m\widehat{A}$, we have $i(f)\notin \varpi^m\widehat{A}$.
 In particular $f\neq 0$, therefore, by the property (c) for $A$, we have a decomposition $f=c_ff'$
 with $c_f\in k^\circ$ and $f'\in A\setminus \mathfrak{m}A$.
 As $i(f)\notin \varpi^m\widehat{A}$, $c_f$ does not lie in $\varpi^mk^\circ$.
 Hence we have $\varpi^m/c_f\in \mathfrak{m}$.

 If we write $a-i(f)=\varpi^ma_1$ with $a_1\in\widehat{A}$, then
 \[
  a=i(f)+\varpi^ma_1=c_f\bigl(i(f')+(\varpi^m/c_f)a_1\bigr).
 \]
 Since $A/\mathfrak{m}A\cong\widehat{A}/\mathfrak{m}\widehat{A}$,
 we have $i(f')\in \widehat{A}\setminus \mathfrak{m}\widehat{A}$. 
 Together with $(\varpi^m/c_f)a_1\in \mathfrak{m}\widehat{A}$,
 we obtain $i(f')+(\varpi^m/c_f)a_1\in \widehat{A}\setminus \mathfrak{m}\widehat{A}$.
 Therefore $a=c_f(i(f')+(\varpi^m/c_f)a_1)$ gives the desired decomposition.
\end{prf}

\begin{cor}\label{cor:flatness-complete}
 Assume that $A$ is $\varpi$-adically complete and satisfies the properties (a), (b) and (c).
 Fix an integer $q>1$ and let $(f_m)_{m\ge 1}$ and $(g_m)_{m\ge 1}$ be as in Corollary \ref{cor:flatness}.
 We denote by $\overline{I}$ the closure in $A$ of the ideal $I=(f_m-g_m\mid m\ge 1)$.
 Then, $A/\overline{I}$ is flat over $k^\circ$.
\end{cor}

\begin{prf}
 By Corollary \ref{cor:flatness}, $A/I$ is $\varpi$-torsion free.
 By \cite[Chapter 0, Proposition 7.4.5]{Fujiwara-Kato}, the $\varpi$-adic completion of $A/I$ equals
 $A/\overline{I}$. Hence Proposition \ref{prop:abc-completion} tells us that $A/\overline{I}$ is also
 $\varpi$-torsion free, that is, flat over $k^\circ$.
\end{prf}

\begin{exa}\label{exa:perfection}
 Let $n\ge 0$ be an integer.
 \begin{enumerate}
  \item The polynomial ring $k^\circ[T_1,\ldots,T_n]$ obviously satisfies the properties (a), (b), (c).
  \item Suppose that the characteristic of $\kappa=k^\circ/\mathfrak{m}$ is $p>0$. Let $q$ be a power of $p$
	and put $k^\circ[T_1^{q^{-\infty}},\ldots,T_n^{q^{-\infty}}]=\varinjlim_{T_i\mapsto T_i^q}k^\circ[T_1,\ldots,T_n]$. Then, $k^\circ[T_1^{q^{-\infty}},\ldots,T_n^{q^{-\infty}}]$ satisfies the properties (a), (b), (c).
	Indeed, (a) and (b) follow from (i), since these are preserved by filtered inductive limits.
	The property (c) follows from that on $k^\circ[T_1,\ldots,T_n]$ and the fact that
	$\kappa[T_1,\ldots,T_n]\to \kappa[T_1,\ldots,T_n]$; $T_i\mapsto T_i^q$ is an injection.
  \item Under the same assumption as in (ii), 
	let $k^\circ\langle T_1^{q^{-\infty}},\ldots,T_n^{q^{-\infty}}\rangle$ denote the $\varpi$-adic
	completion of $k^\circ[T_1^{q^{-\infty}},\ldots,T_n^{q^{-\infty}}]$.
	Then, Proposition \ref{prop:abc-completion} tells us that 
	$k^\circ\langle T_1^{q^{-\infty}},\ldots,T_n^{q^{-\infty}}\rangle$ satisfies
	the properties (a), (b), (c). 
 \end{enumerate}
\end{exa}

\begin{cor}\label{cor:flatness-dimension}
 Let $A=k^\circ\langle T_1^{q^{-\infty}},\ldots,T_n^{q^{-\infty}}\rangle$ be
 as in Example \ref{exa:perfection} (iii), and $(f_m)_{m\ge 1}$, $(g_m)_{m\ge 1}$ sequences of elements in 
 $A$ satisfying $f_{m+1}^q=f_m$ and $g_{m+1}^q=g_m$.
 We denote by $\overline{I}$ the closure in $A$ of the ideal $I=(f_m-g_m\mid m\ge 1)$.
 Assume that $\dim (A/\overline{I}\otimes_{k^\circ}\kappa)<n$.
 Then, $A/\overline{I}$ is flat over $k^\circ$.
\end{cor}

\begin{prf}
 Since $A/\overline{I}$ is the $\varpi$-adic completion of $A/I$, as in the proof of
 Proposition \ref{prop:abc-completion} we have 
 $A/\overline{I}\otimes_{k^\circ}\kappa=A/I\otimes_{k^\circ}\kappa=\kappa[T_1^{q^{-\infty}},\ldots,T_n^{q^{-\infty}}]/(\overline{f_m}-\overline{g_m}\mid m\ge 1)$, where $\overline{f_m}$ and $\overline{g_m}$ denote the images of
 $f_m$ and $g_m$ in $\kappa[T_1^{q^{-\infty}},\ldots,T_n^{q^{-\infty}}]$, respectively. Since
 $\kappa[T_1^{q^{-\infty}},\ldots,T_n^{q^{-\infty}}]$ is $n$-dimensional,
 the condition $\dim (A/\overline{I}\otimes_{k^\circ}\kappa)<n$ implies that
 $\overline{f_m}-\overline{g_m}\neq 0$ for some $m\ge 1$. 
 Since $\kappa[T_1^{q^{-\infty}},\ldots,T_n^{q^{-\infty}}]$ is an integral domain,
 we have $\overline{f_1}-\overline{g_1}=(\overline{f_m}-\overline{g_m})^{q^{m-1}}\neq 0$,
 in other words, $f_1-g_1\notin \mathfrak{m}A$.
 Hence Corollary \ref{cor:flatness-complete} tells us that $A/\overline{I}$ is flat over $k^\circ$,
 as desired.
\end{prf}

\section{Cohomology of Artin-Schreier sheaves}\label{sec:Artin-Schreier}
Let $k$ be an algebraically closed field of characteristic $p>0$.
Fix a non-trivial additive character $\psi\colon \F_p\to \C^\times$.
We write $\mathcal{L}_\psi$ for the rank $1$ sheaf defined by $\psi$ and the $\F_p$-torsor $a^p-a=t$
over $\A^1=\Spec k[t]$.

\begin{rem}\label{rem:m=1}
 For $m\ge 1$, let $\psi_m$ be the additive character $\psi\circ\Tr_{\F_{p^m}/\F_p}$ of $\F_{p^m}$.
 Then, $\mathcal{L}_\psi$ is the rank $1$ sheaf defined by $\psi_m$ and the $\F_{p^m}$-torsor $a^{p^m}-a=t$
 over $\A^1$.
\end{rem}

\begin{lem}\label{lem:AS-A^1}
 Let $d\ge 1$ be an integer prime to $p$, and $\phi$ denote the morphism $\A^1\to \A^1$; $t\mapsto t^d$.
 Then the following hold.
 \begin{enumerate}
  \item We have $\dim H^1_c(\A^1,\phi^*\mathcal{L}_\psi)=d-1$ and
	$H^i_c(\A^1,\phi^*\mathcal{L}_\psi)=H^i(\A^1,\phi^*\mathcal{L}_\psi)=0$ for $i\neq 1$.
  \item The canonical map $H^1_c(\A^1,\phi^*\mathcal{L}_\psi)\to H^1(\A^1,\phi^*\mathcal{L}_\psi)$ is
	an isomorphism.
 \end{enumerate}
\end{lem}

\begin{prf}
 In \cite[(8.11)]{MR0340258}, (i) and the following is proved:
 \begin{enumerate}
  \item[($\text{ii}'$)] the cup product pairing
	\[
	 H^1_c(\A^1,\phi^*\mathcal{L}_\psi)\times H^1_c(\A^1,\phi^*\mathcal{L}_{\psi^{-1}})\to H^2_c(\A^1,\overline{\Q}_\ell)=\overline{\Q}_\ell(-1)
	\]
	is perfect.
 \end{enumerate}
 The assertion (ii) follows from ($\text{ii}'$) and the Poincar\'e duality.
\end{prf}

\begin{lem}\label{lem:AS-quadratic}
 Fix an integer $n\ge 0$ and consider a non-degenerate quadratic form with $n$ variables
 $Q\colon \A^n\to \A^1$.
 If $p=2$, we assume that $n$ is even. Then, we have $\dim H^n_c(\A^n,Q^*\mathcal{L}_\psi)=1$ and
 the canonical map
 $H^n_c(\A^n,Q^*\mathcal{L}_\psi)\to H^n(\A^n,Q^*\mathcal{L}_\psi)$ is bijective.
\end{lem}

\begin{prf}
 First assume that $p\neq 2$. Then, $Q$ is diagonalizable,
 hence the claim follows from Lemma \ref{lem:AS-A^1} and the K\"unneth formula.

 Next consider the case $p=2$. We may assume that $n>0$.
 Let $X\subset \A^{n+1}$ be the affine variety defined by $y^2-y=Q(x_1,\ldots,x_n)$.
 We have
 \begin{align*}
  H^n_c(X,\overline{\Q}_\ell)&=H^n_c(\A^n,\overline{\Q}_\ell)\oplus H^n_c(\A^n,Q^*\mathcal{L}_\psi)=H^n_c(\A^n,Q^*\mathcal{L}_\psi),\\
  H^n(X,\overline{\Q}_\ell)&=H^n(\A^n,\overline{\Q}_\ell)\oplus H^n(\A^n,Q^*\mathcal{L}_\psi)=H^n(\A^n,Q^*\mathcal{L}_\psi).
 \end{align*}
 Therefore, it suffices to show that $\dim H^n_c(X,\overline{\Q}_\ell)=1$ and
 the map $H^n_c(X,\overline{\Q}_\ell)\to H^n(X,\overline{\Q}_\ell)$ is
 a bijection.
 Consider the quadric $\overline{X}$ in $\P^{n+1}$ defined by the homogeneous equation
 $y^2-yz=Q(x_1,\ldots,x_n)$, which includes $X$ as the complement of the hyperplane section $z=0$.
 We will show that $\overline{X}$ is smooth over $k$. Since $Q$ is non-degenerate and $n$ is even,
 it is ordinary in the sense of \cite[Expos\'e XII]{SGA7}. Therefore, by changing variables, we may
 assume that $Q(x_1,\ldots,x_n)=\sum_{i=1}^mx_ix_{i+m}$, where we put $m=n/2$
 (see \cite[Expos\'e XII, Proposition 1.2]{SGA7}).
 In this case, the quadratic form $y(y-z)-\sum_{i=1}^m x_ix_{i+m}$ 
 is obviously ordinary; in other words, the quadric $\overline{X}$ is smooth over $k$.
 Now \cite[Expos\'e XII, \S 3.6 and Table 3.7]{SGA7} tells us that
 $\dim H^n_c(X,\overline{\Q}_\ell)=1$ and $H^n_c(X,\overline{\Q}_\ell)\to H^n(X,\overline{\Q}_\ell)$
 is an isomorphism.
 This concludes the proof.
\end{prf}

\begin{exa}\label{exa:non-degenerate}
 For an integer $n\ge 0$ and $a\in k^\times$, the quadratic form
 \[
  \nu(x_1,\ldots,x_n)=a\sum_{1\le i\le j\le n}x_ix_j
 \]
 is non-degenerate if $p\neq 2$ or $n$ is even. 
\end{exa}

\begin{prf}
 If $p\neq 2$, this is proved in \cite[Lemma 4.1]{Imai-Tsushima-tame}.
 Assume that $p=2$ and $n$ is even. For $x=(x_1,\ldots,x_n)$ and $y=(y_1,\ldots,y_n)$,
 we have
 \[
  \Phi(x,y):=\nu(x+y)-\nu(x)-\nu(y)=a\sum_{i\neq j}x_iy_j.
 \]
 This satisfies $\Phi(e_i,e_i)=0$ and $\Phi(e_i,e_j)=a$ for $i\neq j$,
 where $e_1,\ldots,e_n$ denote the standard basis.
 It suffices to show that $\Phi$ is non-degenerate.
 Take a non-zero element $x=\sum_{i=1}^n a_ie_i$ with $a_i\in k$.
 Put $s=a_1+\cdots +a_n$. If $s-a_i=0$ for every $1\le i\le n$,
 we have $0=\sum_{i=1}^n(s-a_i)=ns-s=s$, $a_i=s=0$ and $x=0$.
 Hence there exists $1\le j\le n$ such that $s-a_j\neq 0$.
 This $j$ satisfies $\Phi(x,e_j)=a(s-a_j)\neq 0$,
 which concludes that $\Phi$ is non-degenerate.
\end{prf}

\def\cftil#1{\ifmmode\setbox7\hbox{$\accent"5E#1$}\else
  \setbox7\hbox{\accent"5E#1}\penalty 10000\relax\fi\raise 1\ht7
  \hbox{\lower1.15ex\hbox to 1\wd7{\hss\accent"7E\hss}}\penalty 10000
  \hskip-1\wd7\penalty 10000\box7}
  \def\cftil#1{\ifmmode\setbox7\hbox{$\accent"5E#1$}\else
  \setbox7\hbox{\accent"5E#1}\penalty 10000\relax\fi\raise 1\ht7
  \hbox{\lower1.15ex\hbox to 1\wd7{\hss\accent"7E\hss}}\penalty 10000
  \hskip-1\wd7\penalty 10000\box7}
  \def\cftil#1{\ifmmode\setbox7\hbox{$\accent"5E#1$}\else
  \setbox7\hbox{\accent"5E#1}\penalty 10000\relax\fi\raise 1\ht7
  \hbox{\lower1.15ex\hbox to 1\wd7{\hss\accent"7E\hss}}\penalty 10000
  \hskip-1\wd7\penalty 10000\box7}
  \def\cftil#1{\ifmmode\setbox7\hbox{$\accent"5E#1$}\else
  \setbox7\hbox{\accent"5E#1}\penalty 10000\relax\fi\raise 1\ht7
  \hbox{\lower1.15ex\hbox to 1\wd7{\hss\accent"7E\hss}}\penalty 10000
  \hskip-1\wd7\penalty 10000\box7} \def\cprime{$'$} \def\cprime{$'$}
  \newcommand{\dummy}[1]{}
\providecommand{\bysame}{\leavevmode\hbox to3em{\hrulefill}\thinspace}
\providecommand{\MR}{\relax\ifhmode\unskip\space\fi MR }
\providecommand{\MRhref}[2]{%
  \href{http://www.ams.org/mathscinet-getitem?mr=#1}{#2}
}
\providecommand{\href}[2]{#2}

\end{document}